\documentclass[ijoc,nonblindrev]{informs3} %

\OneAndAHalfSpacedXII %
\usepackage{natbib}
 \bibpunct[, ]{(}{)}{,}{a}{}{,}%
 \def\bibfont{\small}%

 \usepackage[inline]{enumitem}
\newlist{inlinelist}{enumerate*}{1}
\setlist[inlinelist]{label=(\roman*)}

\usepackage{amsfonts}
\usepackage{amsmath}
\usepackage{array}
\usepackage{mathtools}
\usepackage{xfrac}
\usepackage{IEEEtrantools}
\usepackage{hyperref}

\usepackage{booktabs}
\usepackage{multirow}

\DeclarePairedDelimiter\Paren{\lparen}{\rparen}
\DeclarePairedDelimiter\Brac{\lbrace}{\rbrace}
\DeclarePairedDelimiter\Brak{\lbrack}{\rbrack}
\DeclarePairedDelimiter\Card{\lvert}{\rvert}

\DeclarePairedDelimiter\Ceil{\lceil}{\rceil}
\newcommand{\paren}[1]{\Paren*{#1}}
\newcommand{\brac}[1]{\Brac*{#1}}
\newcommand{\brak}[1]{\Brak*{#1}}
\newcommand{\card}[1]{\Card*{#1}}

\newcommand{\ceil}[1]{\Ceil*{#1}}

\TheoremsNumberedThrough     %
\EquationsNumberedThrough    %
\begin{document}
\RUNAUTHOR{Farham, Rostami, and Haughton}

\RUNTITLE{Vehicle utilization in the SAHLP}

\TITLE{Vehicle Utilization in Hub Network Design: Exploiting Economies of Scale in Transportation}

\ARTICLEAUTHORS{%
    \AUTHOR{Mohammad Saleh Farham}
    \AFF{Lazaridis School of Business \& Economics, Wilfrid Laurier University \EMAIL{mfarham@wlu.ca}}%
    \AUTHOR{Borzou Rostami}
    \AFF{Alberta School of Business, University of Alberta, \EMAIL{borzou@ualberta.ca}}%
    \AUTHOR{Michael Haughton}
    \AFF{Lazaridis School of Business \& Economics, Wilfrid Laurier University \EMAIL{mhaughton@wlu.ca}}%
} %

\ABSTRACT{%
    We study a vehicle-based hub network design problem (HNDPv) with the main applications in freight distribution and parcel delivery systems, where the economies of scale stem from the effective utilization of vehicles that move consolidated freight. The HNDPv is a generalization of the classical single allocation hub location problem, in which the transportation costs are stepwise functions of the number (and type) of vehicles that move the demand. We present the quadratic mixed-integer programming formulation of the problem and its linear reformulation.
    Exploiting the special structures of the linearized model, we develop a branch-and-cut method based on Benders decomposition with solely feasibility subproblems. We derive closed-form solutions for the extreme rays of the feasibility subproblems that improve the efficiency of the proposed algorithm through generating stronger feasibility cuts.
    We also address the HNDPv under demand uncertainty and show the flexibility of our solution methodology in handling the stochastic variant of the problem. To evaluate the efficiency of our models and solution approaches, we perform extensive computational experiments on  uncapacitated and capacitated instances of the problem derived from the classical Australian Post dataset. The results show a considerable advantage of using HNDPv compared to the classical HLP with constant discount factors in terms of  vehicle utilization and total transportation costs. Our computational experiments also demonstrate the efficiency of our proposed solution method in solving large-scale problem instances.
}%

\KEYWORDS{Hub location problem; vehicle utilization; economies of scale; demand uncertainty; Benders decomposition}
\HISTORY{}

\maketitle

\section{Introduction}
\label{sec:intro}

Consolidation-based freight transportation is used for systems where several freight loads of different demand nodes are aggregated to be transported in a less costly manner. When the direct shipments between the origin and destination of demands are not economically justifiable or even feasible, such systems provide a profitable balance between economy-of-scale-based costs and high service quality.

Postal and parcel delivery companies, less-than-truckload (LTL) motor carriers, railroads or maritime liner navigation companies, and air/land or water/land-based intermodal carriers use consolidation centers, called hub facilities, to centralize commodity handling and sorting operations, reduce setup costs, and achieve economies of scale on transportation costs by consolidating flows. This builds a hierarchical network called a \textit{hub network}, where at the \textit{access-level}, individual demand nodes connect to the hubs, and at the \textit{hub-level}, interconnected hubs send and receive consolidated flows.
The objective is to minimize the cost of locating hubs and transporting origin-destination (OD) flows on access and hub-level links.

Economies of scale in freight transportation networks are directly related to the level at which hub/vehicle capacities are utilized to handle/transport large volumes of loads.
At the hub-level, loads are consolidated and can be moved in bulks. Therefore, vehicles that travel on the inter-hub links are commonly large and are made to transport high-volume of loads over long distances efficiently (e.g., cargo jets, trains with freight wagons, etc.).
Vehicles that are operated at the access-level are different.
In a postal delivery system, for instance, small or medium-size trucks are used to transport parcels from hubs (e.g., airports) to non-hubs (e.g., regional collection/distribution centers).
Although the cost of using inter-hub vehicles might be larger, they are capable of moving bulks of consolidated freight more efficiently.
Therefore, properly utilizing the capacity of inter-hub vehicles loads to a smaller unit transportation cost on hub arcs compared to the access arcs.
This enables transportation carriers to exploit economies of scale and achieve lower transportation costs by consolidating loads into larger vehicles.

In this paper, we consider vehicle-related decisions and utilization costs in the hub network design problem to determine the resources required to route the flow through the network. More specifically, we study the vehicle-based hub network design problem (HNDPv) with the single-assignment strategy in which each demand node is assigned to precisely one hub. The objective is to select hub nodes, assign demand nodes to the selected hubs, and determine the number and type of vehicles that travel on hub-level and access-level networks to minimize the total location and vehicle cost. We present mathematical programming models of the HNDPv and solve large-scale instances of the problem by developing an exact branch-and-cut solution method based on Benders decomposition. We show the advantages of using HNDPv compared to the classical HLP with constant discount factors in terms of vehicle utilization and total transportation costs.

The HNDPv is a strategic problem where the long-term location and allocation decisions determine the underlying network structure and the fleet size decisions determine the investment that should be made to move the flow through the network. However, since selecting the optimal hub location and an appropriate vehicle fleet management are directly related to the amount of anticipated OD shipping volumes, a poor demand estimation may lead to an unsustainable or infeasible solution under the actual demand. To address the demand uncertainty, we extend the mathematical programming models of HNDPv in a stochastic environment and show how to adjust our branch-and-cut solution method to handle the demand uncertainty.


\subsection{Related Literature}
\label{sec:literature}

Despite the important economic impact of vehicle utilization in hub network design problems, it has not gained enough attention in the literature.
The large majority of in the hub location problems (HLPs) in the literature model the transportation cost as a linear function of volume transported on the network links. To address economies of scale, the transportation cost on the inter-hub links are multiplied by a fixed constant $\alpha$ called the \textit{discount factor} \citep[see][for a recent review]{Alumur+Campbell-EurJOpeRes-2021}.
While this approach leads to a tractable (linear) cost function to minimize, it does not provide an adequate representation of economies of scale in consolidation-based transportation systems. Considering a fixed discount factor on all inter-hub links may lead to an oversimplified modeling of economies of scale and produce a poor estimation of the real savings.
This simplification may result in solutions that grant discounted transportation on the inter-hub links, even though the flow on these links are poorly consolidated \citep{Kimms--2006, Real+Contreras-TraResParE:LogandTraRev-2021}.
In addition, it requires advanced knowledge about the technologies used in the system to estimate potential $\alpha$ values in advance.

To reflect economies of scale more realistically, one can consider either a decreasing unit cost for increasing transport volumes or a function that reflects the actual cost of vehicle utilization on the inter-hub links \citep{Alumur+Campbell-EurJOpeRes-2021}.
The former leads to a nonlinear concave cost function of the flow on an inter-hub link \citep{Horner+OKelly-JTraGeo-2001}.
Such an approach, however, has two drawbacks. First, it does not consider the technology (i.e., the means of transport) used on the inter-hub links. Hence, vehicle utilization cost and capacity are ignored.
Second, considering concave functions in mathematical modeling brings new computational challenges.
As a result, researchers often consider representing such functions by their linear approximation \citep{OKelly+Bryan-TraResParB:Met-1998,Racunica+Wynter-TraResParB:Met-2005,deCamargo+deMiranda-TraSci-2009,Rostami+Chitsaz-OpeRes-2022}.

\begin{figure}
    \FIGURE
    {\includegraphics[width=0.42\textwidth]{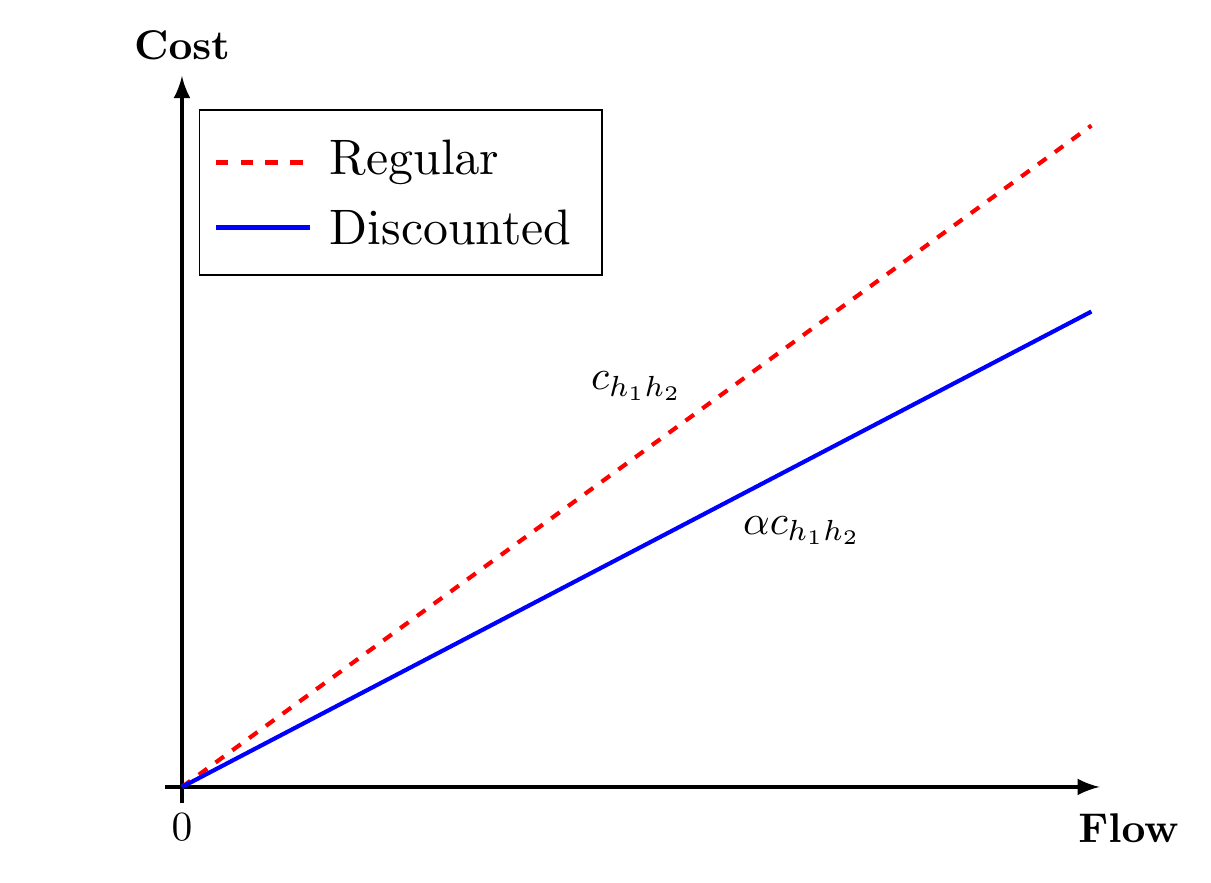} \quad \includegraphics[width=0.42\textwidth]{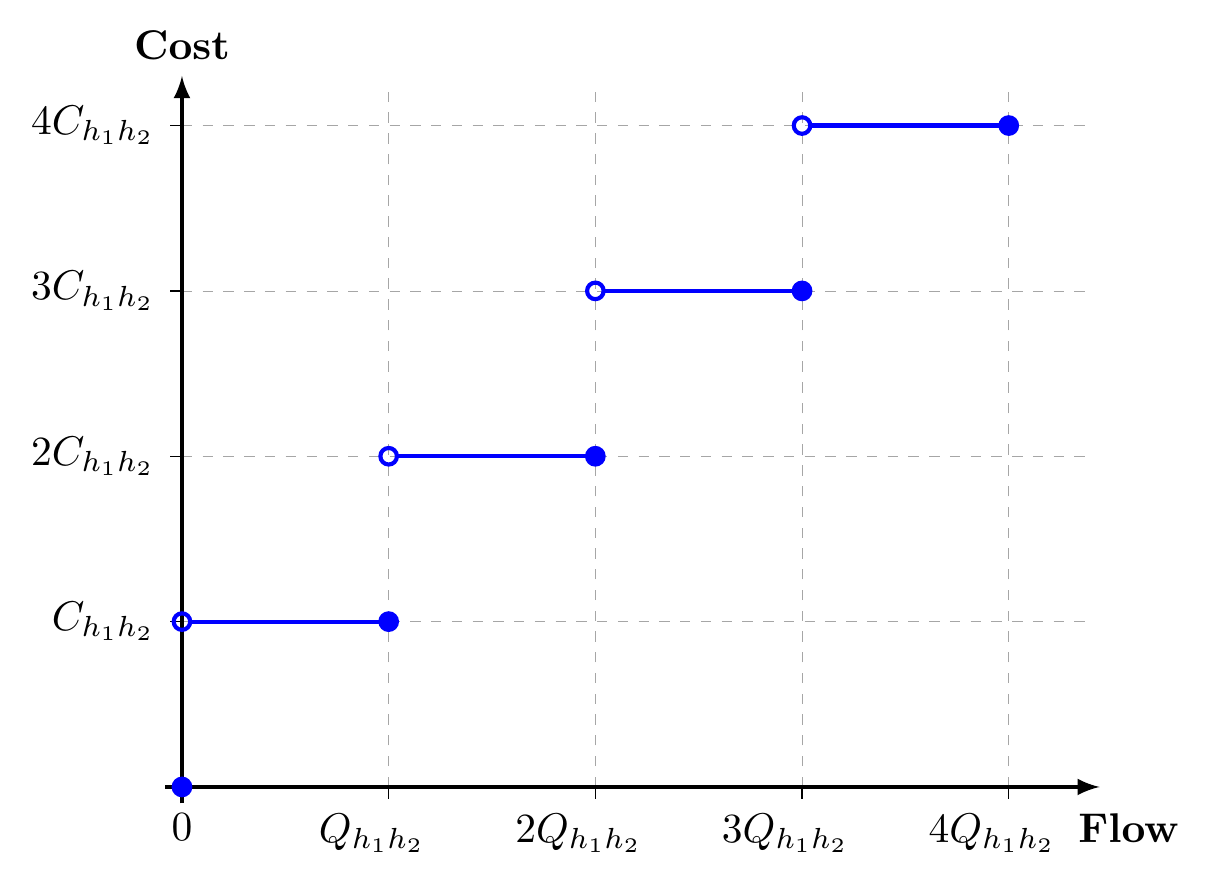}}
    {Transportation cost function over an inter-hub arc.\label{fig:cost-func}}
    {Left: Linear discount-based cost function. Right: Stepwise vehicle-based cost function.}
\end{figure}

Models with vehicle-based cost functions, on the other hand, take the vehicle characteristics into account by measuring the transportation costs per vehicle (as in freight transportation networks) or per capacitated link (as in telecommunication networks). Note that in the context of telecommunication, vehicle capacities translate to link capacities and the corresponding problem is called HLP with modular links (HLPm), which   decides the number of capacitated links to build between each pair of nodes in the network. Such cost functions often lead to a stepwise linear cost function. \autoref{fig:cost-func} plots a conventional discount-based and a vehicle-based cost function that calculates the transportation cost on an inter-hub arc $(h_1,h_2)$.
In the constant-discount models, the transportation cost per flow unit, say $c_{h_1,h_2}$, is reduced by a factor in the range $(0, 1)$. The cost grows linearly with the amount of flow on the arc.
The vehicle-based function calculates the transportation cost according to how vehicles are utilized.
If each vehicle on arc $(h_1,h_2)$ has a capacity $Q_{h_1,h_2}$ and a fixed utilization cost $C_{h_1,h_2}$, then the cost of transporting loads on that arc is represented by a step-wise function of the number of vehicles in use.

The main difficulty of the HNDPv (and the HLPm), in addition to the natural complexity of the design problem, is dealing with the integer number of vehicles (capacitated links) in the mathematical programming models. 
Different exact and approximation methodologies have been presented in the literature to address this challenge.

In the steam of exact solution approaches, \cite{Yaman+Carello-ComOpeRes-2005}  presents a branch-and-cut 
 method for the HLPm. They solve a set of modified Australian Post (AP) instances with up to 20 nodes to optimality and provide feasible solutions for larger instances with up to 50 nodes using a heuristic. \cite{Rostami+Buchheim-OptOnl-2015} formulate the uncapacitated $p$-Hub median problem with vehicle-based transportation costs and develop a branch-and-bound algorithm where lower bounds are computed using a Lagrangian relaxation algorithm. The authors report the exact solutions for modified AP instances with up to 50 nodes. \cite{Tanash+Contreras-ComOpeRes-2017} develop a Lagrangian-based branch-and-bound algorithm for the HLPm  with an incomplete inter-hub network structure. In such problems, the OD paths are allowed to visit more than two hubs, and the design of the inter-hub network is part of the decision process. They solve AP instances with up to 40 nodes to optimality. 

In the stream of heuristics, different solution methods based on local search, iterated greedy search combined with memory strategies, and variable neighborhood search have been proposed in the literature  \cite[see][]{Carello+DellaCroce-Net-2004,Corberan+Peiro-JHeu-2016,Hoff+Peiro-ComOpeRes-2017,Serper+Alumur-TraResParB:Met-2016,KeshvariFard+Alfandari-ComOpeRes-2019}. In particular, \cite{KeshvariFard+Alfandari-ComOpeRes-2019} proposes an approximation method to transform the stepwise vehicle-based cost function to a linear function of the flow. The authors claim that by choosing the proper intercept and slope parameters of the linear function, one can obtain a solution close to the one found for the original stepwise cost function. The authors provide the inexact to a set of problem test instances from CAB, AP, and Turkish datasets with up to 50 nodes. While this approach is more efficient, its solution quality highly depends on the estimated parameters of the generalized linear function and the capacity of the inter-hub vehicles. This approach may lead to inferior solutions when the flow on some inter-hub links is small or the vehicle capacities are large.

\autoref{tab:liter} lists recent research in the literature on HNDPv and the HLPm with the single-allocation strategy. For each study, the table specifies the considered inter-hub network structure, whether the demand uncertainty is addressed, whether all vehicles (including the access-level vehicles) are capacitated, the proposed solution approach, and the largest instance size which could be solved exactly.
Although there are a few other studies that consider HLP variants with step-wise cost functions \citep[e.g.,][]{Kimms--2006, Sender+Clausen--2011,Baumung+Gunduz--2015}, we do not list them in this table as they only introduce problem formulations, and do not present any problem-specific solution algorithms.
The terms that are made bold in \autoref{tab:liter} share similarities with our problem.

The available solution approaches for the HDNPv are either approximation methods or exact algorithms that are only capable of solving small-size problem instances.
Moreover, despite the fact that the demand uncertainty directly affects the flow through the network, hence vehicle utilization, no research has been carried out to investigate the HNDPv with stochastic demand.
Previous HLPs studies only consider sources of uncertainty in the discount-based models \citep[see, for example,][]{Alumur+Nickel-TraResParB:Met-2012,Qin+Gao-JIntMan-2017,Tran+OHanley-TraSci-2017,Wang+Chen-TraSci-2020,Hu+Hu-TraResParB:Met-2021,Rostami+Kammerling-EurJOpeRes-2021}, and leave the stochastic HNDPv unexplored.

\begin{table}
    \TABLE
    {Related literature to the vehicle-based HLP.\label{tab:liter}}
    {\footnotesize \begin{tabular}{@{}llllll@{}}
            \toprule
            \multicolumn{1}{@{}c}{Article}               & \begin{tabular}[c]{@{}c@{}}Inter-hub\\network structure\end{tabular} & \begin{tabular}[c]{@{}c@{}}Demand\\uncertainty\end{tabular} & \begin{tabular}[c]{@{}c@{}}Capacitated\\vehicles\end{tabular} & \multicolumn{1}{c}{Solution approach} & \begin{tabular}[c]{c@{}}Problem type\\(optimal size)\end{tabular} \\ \midrule
            \cite{Carello+DellaCroce-Net-2004}           & \textbf{Complete}                                                    & No                                                          & \textbf{Yes}                                                  & Heuristic                             & CHLP                                                                                       \\
            \cite{Yaman+Carello-ComOpeRes-2005}          & \textbf{Complete}                                                    & No                                                          & \textbf{Yes}                                                  & Branch-and-cut                        & CHLP(20)                                                                                   \\
            \cite{Rostami+Buchheim-OptOnl-2015}          & \textbf{Complete}                                                    & No                                                          & \textbf{Yes}                                                  & Branch-and-bound                      & $p$-HLP(50)                                                                                \\
            \cite{Corberan+Peiro-JHeu-2016}              & \textbf{Complete}                                                    & No                                                          & \textbf{Yes}                                                  & Heuristic                             & CHLP                                                                                       \\
            \cite{Serper+Alumur-TraResParB:Met-2016}     & incomplete                                                           & No                                                          & \textbf{Yes}                                                  & Heuristic                             & CHLP                                                                                       \\
            \cite{Hoff+Peiro-ComOpeRes-2017}             & \textbf{Complete}                                                    & No                                                          & No                                                            & Heuristic                             & CHLP                                                                                       \\
            \cite{Tanash+Contreras-ComOpeRes-2017}       & incomplete                                                           & No                                                          & \textbf{Yes}                                                  & Branch-and-bound                      & HLP(40)                                                                                    \\
            \cite{KeshvariFard+Alfandari-ComOpeRes-2019} & \textbf{Complete}                                                    & No                                                          & \textbf{Yes}                                                  & Approximation                         & $p$-CHLP                                                                                   \\
            This study                                   & \textbf{Complete}                                                    & \textbf{Yes}                                                & \textbf{Yes}                                                  & Branch-and-cut                        & CHLP(200\textsuperscript{$\dagger$}, 75\textsuperscript{$\ddagger$})                       \\ \bottomrule
        \end{tabular}%
    }
    {CHLP: capacitated HLP, $p$-HLP: $p$-hub median problem, $p$-CHLP: capacitated $p$-HLP. \\
        \indent {\textsuperscript{$\dagger$}} Deterministic instance size.\\
        \indent {\textsuperscript{$\ddagger$}} Stochastic instance size.}
\end{table}


\subsection{Our Contribution}

While vehicle-based cost functions provide an adequate description of scale economies, they make the already challenging hub network design problem even more difficult to solve due to the additional integer variables that determine the number of inter-hub vehicles.
Incorporating demand variability, which is necessary to make reliable location/allocation and fleet-size decisions, also adds another layer of complexity.
In this study,
\begin{inlinelist}
    \item we develop an exact solution method based on Benders decomposition. While a natural approach to tackling the problem's difficulty is to project out the integer vehicle variables, we propose an alternative where all the main decisions are made in the master problem. This leads to a branch-and-cut algorithm with feasibility-checking subproblems.

    \item To generate feasibility cuts more efficiently, we derive the extreme rays of the feasibility subproblem in an analytical way that, in addition to preventing solving many linear programs within the search tree, also provides an opportunity to add multi cuts in each call to the subproblem.

    \item We address the problem with demand uncertainty under the stochastic programming framework and show the potential of our solution methodology in solving the deterministic equivalent formulation of the problem.

    \item Extensive computational experiments are conducted to evaluate the potential, robustness, and efficiency of our models and solution methodologies on uncapacitated and capacitated instances derived from the classical Australian Post dataset. The results show a considerable advantage of using HNDPv compared to the classical HLP with constant discount factors in terms of  vehicle utilization and total transportation costs. 
    The proposed solution algorithm is able to solve large-scale deterministic HNDPv instances with up to 200 nodes for the first time in the literature. We also show the capability of the proposed solution methodology in solving problems with uncertain demands and general (or incomplete) inter-hub network structures.
\end{inlinelist}

The remainder of the paper is organized as follows.
\autoref{sec:problem} introduces mathematical formulations for the deterministic HNDPv and presents some valid inequalities.
In \autoref{sec:sol}, we propose our solution approaches for the HNDPv.
The HNDPv under demand uncertainty and its solution method is discussed in \autoref{sec:uncertainty}.
Our computational study in \autoref{sec:experiments} investigates the performance of our solution algorithm in solving a set of benchmark problem test instances and provides managerial insights.
In this section, we also explain how our solution approach can be adopted to solve problems with general inter-hub network structure.
\autoref{sec:conclusion} concludes the paper and highlights future research directions.
\section{Problem Statement and Formulation}
\label{sec:problem}

The HNDPv is defined over a many-to-many network $G=\paren{N,A}$, where $N$ represents the set of demand points (including hubs) and $A$ is the set of (directed) arcs.
The set of candidate hub nodes is denoted as $H \subset N$. Each node in $N$ can potentially be the sender/receiver of specific OD flows. That is, there exists a flow amount of $w_{ij}$ of a single commodity between each pair of $i$ and $j$ in $N$.
We use $O_i=\sum_{j\in N} w_{ij}$ and $D_i=\sum_{j\in N} w_{ji}$ to denote the total amount of flow originated and destined at demand node~$i$, respectively.
The amount of flow that can be consolidated at a hub is usually restricted by hub capacity. Associated with each hub $h\in H$, we consider a limited capacity $U_h$ (that is set to a big number if the hub is uncapacitated) and a fixed setup cost $F_h$.

Node-to-node connections are established through vehicle movements. At the hub level, loads are consolidated and are transported by high-capacity vehicles called \textit{primary} vehicles. At the access-level network, smaller vehicles, called \textit{secondary} vehicles, are utilized to pickup/deliver demands from/to non-hub nodes.
Vehicles that are available at each hub may differ in terms of capacity and cost factors (e.g., fixed utilization cost).
Therefore, we identify the fleet of secondary vehicles at any hub $h\in H$ by capacity $q_h$, a fixed utilization cost $g_h$, and unit traveling cost $b_h$. We define $c_{hi} = g_h + b_h\, d(h,i)$  as the one-way transportation cost on access arc $(h,i)\in A$, where $d(i,j)$ is the distance between nodes $i$ and $j$. Therefore, the cost of serving node~$i\in N$ by a secondary vehicle from hub~$h\in H$ is $c^\pm_{hi} = c_{hi} + b_h\, d(i,h)$.
Similarly, associated with each primary vehicle connecting a pair of hubs $h, k\in H$ are capacity $Q_{hk}$, fixed utilization cost $G_{hk}$, and unit traveling cost $B_{hk}$. Therefore, we let $C_{hk} = G_{hk} + B_{hk}\, d(h,k)$ be the cost of using a primary vehicle on each hub-hub connection $(h,k)$.

Note that in many strategic problems, calculating vehicle costs are based on the utilization cost over a given distance, while the filling quota plays a minor role. Therefore, the vehicle utilization cost is defined as a function of a fixed cost associated with using a vehicle on a specific link, and the cost of traversing that link. The traveling cost is assumed as a function of distance and not the load.
In this way, dispatching vehicles with high load factors leads to a lower total transportation cost.

We assume that the fixed cost, capacity, and unit traveling costs of a primary vehicle are strictly greater than those of a secondary vehicle, and $\sfrac{C_{hk}}{Q_{hk}} < \sfrac{c_{hk}}{q_h}$ holds for $(h,k)\in A$.
This ensures that the unit transportation cost is less on inter-hub arcs compared to the access arcs when vehicles are adequately utilized. Consequently, the economies of scale is enhanced through consolidating freight into primary vehicles.

\subsection{Mathematical Formulation}
\label{sec:model}

The HNDPv aims to
\begin{inlinelist}
    \item select a set of hub nodes,
    \item assign non-hub nodes to the selected hubs by satisfying the single-assignment property, and
    \item determine the number and type of primary and secondary vehicles,
\end{inlinelist}
such that the overall hub location and vehicle utilization costs are minimized.
Therefore, the major decisions in the HNDPv involve hub location, demand node allocation, and vehicle fleet management. We define a binary decision variable $x_{hi}$ to indicate whether node $i\in N$ is assigned to hub $h\in H$. Variable $x_{hh}$ indicates whether node $h\in H$ is selected as a hub. Moreover, we define $y_{ij}$ as an integer variable that determines the number of vehicles needed on arc $(i,j)\in A$.

One of the main component of the objective cost function is the cost corresponding to the number of primary and secondary vehicles.
Due to the single-assignment property, we have full information on the amount of flow on each selected access arc traversed in direct shipments. For example, when node $i$ is assigned to hub $h$, the total flow on arc $(h,i)$ is equal to $D_i$ and the total flow on arc $(i,h)$ is equal to $O_i$.
We can calculate the number of vehicles required for delivering and picking up the demand of node $i$ as $n^+_{hi} = \ceil{D_i / q_h}$ and $n^-_{ih} = \ceil{O_i / q_h}$, respectively, where $\ceil{\,\cdot\,}$ is the ceiling function.
Therefore, the total number of secondary vehicles $y_{hi}$ that need to dispatch from hub $h$ to serve demand node $i\in N$ via direct shipment can be prepossessed and set to $n^\pm_{hi}=\max\brac{n^+_{hi}, n^-_{hi}}$.
This number is based on the fact that vehicles are available at the hub locations and start and end their trip at their hosting hubs. In the HLPs with modular links, $n^\pm_{hi}$ translates into the number of links with capacity $q_h$ that need to be installed in order to serve demand node~$i$ from hub~$h$ \citep[see][]{Yaman+Carello-ComOpeRes-2005, Corberan+Peiro-JHeu-2016}.
Therefore, we define the following cost function at the access level network
\begin{IEEEeqnarray}{rCl}
    \mathit{DC}(\BFx) &=& \sum_{h\in H} \sum_{i\in N : i\neq h} n^\pm_{hi}\, c^\pm_{hi}\, x_{hi}. \label{eq:vcost-dir}
\end{IEEEeqnarray}

However, the decisions about the number of primary vehicles in the hub-level network can not be prepossessed and must explicitly be addressed by $y$ variables. The following function calculates the total transportation cost at the hub level.
\begin{IEEEeqnarray}{l}
    \mathit{HC}(\BFy) = \sum_{h\in H} \sum_{k\in H: k \neq h} C_{hk}\, y_{hk}.
\end{IEEEeqnarray}

Finally, the total hub location cost is defined as:
\begin{IEEEeqnarray}{l}
    \mathit{LC}(\BFx) = \sum_{h\in H} F_h\, x_{hh}.
\end{IEEEeqnarray}

Using the defined decisions and their associated costs, we model the HNDPv as the mixed-integer quadratic formulation given below.
\begin{IEEEeqnarray}{lll}
    \IEEEyesnumber\label{eq:hlp}
    \IEEEyessubnumber*
    \text{Minimize}\ & \IEEEeqnarraymulticol{2}{l}{\mathit{LC}(\BFx) + \mathit{DC}(\BFx) + \mathit{HC}(\BFy)} \label{eq:hlp.0}\\ [1ex]
    \text{subject to:}\ & \sum_{i\in N}\sum_{j\in N} w_{ij}\,x_{hi} x_{kj} \leq Q_{hk}\, y_{hk},\  & h, k\in H, \label{eq:hlp.y}\\
    & \BFx \in \BFcalX, \label{eq:hlp.x} \\
    & \BFy \in \mathbb{Z}^{\card{H}\times\card{H}}_{\geq 0}, \label{eq:hlp.00}
\end{IEEEeqnarray}
where $\BFcalX$ is the set of constraints that ensure a feasible assignment.
Objective function \eqref{eq:hlp.0} minimizes the total cost of locating hubs and utilizing vehicles to transport flow through the network.
Constraint \eqref{eq:hlp.y} relates assignment and vehicle variables.
It calculates the total flow on each hub arc and ensures that a correct number of vehicles travel on that arc.
Constraint \eqref{eq:hlp.00} restricts $y$ to nonnegative integer values. For simplicity, we denote $\mathbb{Z}^{\card{H}\times\card{H}}_{\geq 0}$ by $\BFcalY$.

Set $\BFcalX$ is well-defined in the HLP literature \citep[see][for a review]{Farahani+Hekmatfar-ComIndEng-2013}.
A standard feasible set $\BFcalX$ is given as:
\begin{IEEEeqnarray}{l}
    \BFcalX = \brac{x_{hi}\in\brac{0, 1},\ h\in H, i\in N \left\lvert
        \begin{array}{ll}
            \sum_{h\in H} x_{hi} = 1,\  & i\in N,         \\
            x_{hi} \leq x_{hh},\        & h\in H, i\in N.
        \end{array}
        \right.},
\end{IEEEeqnarray}
where the first constraint assigns each node to exactly one hub, and the second constraint restricts assignments to the selected hubs.
If a $p$-hub median problem is targeted, equality constraint $\sum_{h\in H} x_{hh} = p$ is added to $\BFcalX$ to ensure that exactly $p$ hubs are selected.
When hub $h$ has a limited capacity $U_h$, $\BFcalX$ also contains the following constraint to ensure that the total outgoing demand from a hub does not exceed its capacity.
\begin{IEEEeqnarray}{l}
    \sum_{i\in N} O_i x_{hi} \leq U_h\, x_{hh},\ h \in H. \label{eq:hlp.cap}
\end{IEEEeqnarray}

\subsection{A Linear Reformulation}

The HNDPv formulation \eqref{eq:hlp} is a constrained binary quadratic program, which is intractable for many standard solvers.
One can use the ``path-based'' formulation of \cite{Skorin-Kapov+Skorin-Kapov-EurJOpeRes-1996} or the ``flow-based'' formulation of \cite{Ernst+Krishnamoorthy-LocSci-1996} to obtain an equivalent MIP formulation.
Although the flow-based formulation is widely considered as the most effective model for the classical single allocation HLPs, a crucial assumption for its validity is that the triangle inequality for the transportation costs holds \citep{Correia+Nickel-EurJOpeRes-2010}.
As the inter-hub transportation costs are vehicle-dependent in our application, the triangular inequality condition does not generally hold and, therefore, the classical flow-based technique cannot be applied.
Here, we use a modified version of the flow-based linearization of \cite{Rostami+Chitsaz-OpeRes-2022} that always provides a valid linearization regardless of the underlying cost structure.
Consider a non-negative variable $z_{ihk}$ that determines the amount of node $i$'s demand transported from hub $h$ to hub $k$, i.e., $z_{ihk}=x_{hi}\sum_{j} w_{ij} x_{kj}, \forall i\in N, h,k\in H$.
Then, constraint \eqref{eq:hlp.y} can be replaced by
\begin{IEEEeqnarray}{l}
    \IEEEyesnumber\label{eq:hlp.lin}
    \IEEEyessubnumber*
    \sum_{i\in N} z_{ihk} \leq Q_{hk}\, y_{hk},\ h, k\in H, \label{eq:hlp.lin.0}
\end{IEEEeqnarray}
where variable $\boldsymbol{z}$ is determined by the following set of constraints.
\begin{IEEEeqnarray}{ll}
    \IEEEyessubnumber*
    \sum_{k\in H} z_{ihk} = O_i\, x_{hi},\  & h\in H, i \in N, \label{eq:linflow1} \\
    \sum_{h\in H} z_{ihk} = \sum_{j\in N} w_{ij}\, x_{kj},\  & k\in H, i \in N, \label{eq:linflow2} \\
    \BFz \in \mathbb{R}^{\card{N}\times\card{H}\times\card{H}}_{\geq 0}. \label{eq:hlp.lin.00}
\end{IEEEeqnarray}

Replacing the quadratic constraints \eqref{eq:hlp.y} by, new set of constraints in \eqref{eq:hlp.lin} we obtain the following mixed-integer linear programming reformulation
\begin{IEEEeqnarray}{l}
    \label{model:SM}
    \text{P}:\quad \min_{\BFx,\BFy,\BFz}\brac{ \mathit{LC}(\BFx) + \mathit{DC}(\BFx) + \mathit{HC}(\BFy)\left\lvert \eqref{eq:hlp.lin}, ~\BFx \in \BFcalX, \,
        \BFy \in \BFcalY\right.}.
\end{IEEEeqnarray}

The validity of this linearization follows from the equivalence between constraints \eqref{eq:linflow1} and \eqref{eq:linflow2} and the mathematical definition of the flow variables \citep[see Remark 3 in][]{Rostami+Chitsaz-OpeRes-2022}.
\section{A Benders Decomposition-based Solution Algorithm}
\label{sec:sol}

Model P can be solved by state-of-the-art MIP solvers. However, it is still very challenging due to a large number of variables and linking constraints.
A natural way to handle such a difficulty is to apply a Benders decomposition (BD) to project out integral variables $\BFy$ and deal with them in a subproblem.
However, since the subproblem involves integral variables, it requires special treatments, such as the integer L-shaped method \citep{Laporte+Louveaux-OpeResLet-1993}, to generate optimality cuts to the Benders master problem. Given that integer L-shaped cuts are usually not strong, one can solve the linear programming (LP) relaxation of the resulting subproblems and add Benders optimality cuts to improve the global lower bound on the value of the subproblems \citep[see][for a recent implementation on the integer L-shaped method to solve an HLP]{Rostami+Chitsaz-OpeRes-2022}. Our preliminary experiments, however, showed that such a treatment is time-consuming and negatively affects the performance of the BD algorithm. Therefore, we do not provide details of the integer L-shaped method  here.

In what follows, we describe an alternative Benders decomposition in \autoref{sec:Feasibility} in which the flow variables are projected out and will be handled through feasibility cuts. In \autoref{sec:closed-form}, we show how to compute such feasibility cuts efficiently. Finally, in \autoref{sec:valid-inequality}, we present some valid inequalities to enhance the algorithm.

\subsection{Benders Feasibility Cuts}
\label{sec:Feasibility}

Consider P in \ref{model:SM} and project out the flow variables from the model, while keeping the location/allocation and integer vehicle variables. This leads to an integer Benders master problem (BMP) and a linear program (LP) subproblem.
The Benders subproblem is defined for a given $\paren{\overline{\BFx}, \overline{\BFy}} \in \BFcalX\times \BFcalY$ in order to check whether the decided vehicle variables yield a feasible flow on hub-hub connections. This problem is referred to as the BSP$(\overline{\BFx},\overline{\BFy})$ given by
\begin{IEEEeqnarray}{llll}
    \IEEEyesnumber\label{eq:bsp}
    \IEEEyessubnumber*
    \text{BSP}\paren{\overline{\BFx}, \overline{\BFy}}:\qquad & \min\ & \IEEEeqnarraymulticol{2}{l}{0} \label{eq:psp.0} \\ [1ex]
    &\text{s.t.}\ & \sum_{i\in N} z_{ihk} \leq Q_{hk} \overline{y}_{hk},\  & h, k\in H,  \\
    && \sum_{k\in H} z_{ihk} = O_i \overline{x}_{hi},\  & h\in H, i \in N, \label{eq:psp.1}\\
    && \sum_{h\in H} z_{ihk} = \sum_{j\in N} w_{ij} \overline{x}_{kj},\  & k\in H, i \in N, \label{eq:psp.2}\\
    && \BFz \in \mathbb{R}^{\card{N}\times\card{H}\times\card{H}}_{\geq 0}. \label{eq:psp.00}
\end{IEEEeqnarray}

By defining dual variables $\BFlambda, \BFmu, \BFnu$, we can write the  dual of the BSP$\paren{\overline{\BFx}, \overline{\BFy}}$, called DSP$\paren{\overline{\BFx}, \overline{\BFy}}$, as the following LP formulation:
\begin{IEEEeqnarray}{lll}
    \IEEEyesnumber\label{eq:dsp}
    \IEEEyessubnumber*
    \text{DSP}\paren{\overline{\BFx}, \overline{\BFy}}:\quad & \max\ & \sum_{h \in H}\sum_{k \in H} Q_{hk}\overline{y}_{hk}\, \lambda_{hk}
    + \sum_{h \in H}\sum_{i \in N} \paren{ O_i \overline{x}_{hi}\, \mu_{hi} + \sum_{j \in N} w_{ij} \overline{x}_{hj}\, \nu_{hi}} \label{eq:dsp.0} \IEEEeqnarraynumspace \\ [1ex]
    & \text{s.t.}\ & \lambda_{hk} + \mu_{hi} + \nu_{ki} \leq 0, \qquad h, k\in H, i \in N, \label{eq:dsp.1} \\
    && \BFlambda \in \mathbb{R}^{\card{H}\times\card{H}}_{\leq 0}, \BFmu, \BFnu \in \mathbb{R}^{\card{H}\times\card{N}}.  \label{eq:dsp.00}
\end{IEEEeqnarray}

The feasible set of $\text{DSP}\paren{\overline{\BFx}, \overline{\BFy}}$ is independent of the choice of $\paren{\overline{\BFx}, \overline{\BFy}}$. Therefore, if it is not empty, the DSP$\paren{\overline{\BFx}, \overline{\BFy}}$ becomes either feasible or unbounded for any arbitrary choice of $\paren{\overline{\BFx}, \overline{\BFy}}$. In the former case, no further action is required and $\paren{\overline{\BFx}, \overline{\BFy}}$ is feasible and hence optimal for the problem. In the latter case, given the set of extreme rays $\mathcal{R}$ of the set $\brac{\paren{\overline{\BFlambda}, \overline{\BFmu}, \overline{\BFnu}} \colon \eqref{eq:dsp.1}}$, there is an unbounded ray $(\overline{\BFlambda}^r, \overline{\BFmu}^r, \overline{\BFnu}^r)$, $r\in \mathcal{R}$ for which
\begin{IEEEeqnarray}{l}
    \Omega^r =\sum_{h \in H}\sum_{k \in H} Q_{hk}\overline{\lambda}_{hk}^r\,y_{hk} + \sum_{h \in H}\sum_{i \in N} \paren{ O_i \overline{\mu}_{hi}^r\, x_{hi}  + \sum_{j \in N} w_{ij}\overline{\nu}_{hi}^r\, x_{hj}} > 0. \label{eq:fcut}
\end{IEEEeqnarray}

We must cut solution $\paren{\overline{\BFx}, \overline{\BFy}}$ to restrict movement in this direction. This will result in the following reformulation of model P refer to master problem
\begin{IEEEeqnarray}{lll}
    \IEEEyesnumber\label{eq:bmp}
    \IEEEyessubnumber*
    \text{BMP}:\qquad & \min\ & \mathit{LC}(\BFx) + \mathit{DC}(\BFx) + \mathit{HC}(\BFy) \\
    & \text{s.t.}\ & \Omega^r \leq 0,\ r\in R \label{eq:feasibilityCut}\\
    && \BFx \in \BFcalX,\ \BFy\in \BFcalY.
\end{IEEEeqnarray}

In our implementation, which is evaluated in \autoref{sec:experiments}, we solve BMP using a branch-and-cut framework of a state-of-the-art optimization solver. The feasibility cuts are incorporated into the master problem by using callbacks, allowing to add the cutting planes step-by-step. A callback is executed whenever an optimal solution of the LP-relaxation is found at the root node of the branch-and-bound-tree or an incumbent solution at any node of the branch-and-bound-tree is found. For the current choice of variables $\paren{\overline{\BFx}, \overline{\BFy}}$, if this solution satisfies the following conditions, then it is also feasible to the original problem \eqref{eq:hlp} and no further action is required.
\begin{IEEEeqnarray}{l}
    \begin{cases}
        z_{ihk} =  \overline{x}_{hi}\sum_{j} w_{ij} \overline{x}_{kj},\  & \forall i\in N, h, k\in H, \\
        \sum_{i\in N} z_{ihk}  \leq Q_{hk}\, \overline{y}_{hk},\         & h, k\in H.
    \end{cases}
    \label{eq:fcond}
\end{IEEEeqnarray}
Otherwise, i.e., if condition \eqref{eq:fcond} is violated, the  feasibility cut \eqref{eq:feasibilityCut} is added to the BMP to cut off the current $\paren{\overline{\BFx}, \overline{\BFy}}$.  To find the feasibility cuts, one can solve the BSP, or its dual, using a standard LP solver.

\begin{remark}
    The Benders decomposition approach developed here can also be applied to solve the HNDPv with the general hub network structures with small modifications. The general hub network structure allows the OD pairs' demands to go through more than two hubs if needed \cite{Tanash+Contreras-ComOpeRes-2017}. Therefore, the quadratic constraint \eqref{eq:hlp.y} and inequality \eqref{eq:valid3} are no longer valid. The Benders subproblem should incorporate a flow conservation constraint for each hub node to ensure that the flows are correctly routed through the network. That is, we have  
    
    \begin{IEEEeqnarray}{llll}
        \IEEEyesnumber\label{eq:gbsp}
        \IEEEyessubnumber*
        \text{G-BSP}\paren{\overline{\BFx}, \overline{\BFy}}:\ \ & \min\ & \IEEEeqnarraymulticol{2}{l}{0} \label{eq:gpsp.0} \\ [1ex]
        & \text{s.t.}\ & \eqref{eq:hlp.lin.0}, \eqref{eq:hlp.lin.00}, \nonumber \\
        && \sum_{k\in H} z_{ihk} - \sum_{k\in H} z_{ikh} = O_i \overline{x}_{hi} - \sum_{j\in N} w_{ij} \overline{x}_{hj},\  & h\in H, i \in N. \IEEEeqnarraynumspace \label{eq:gpsp.00}
    \end{IEEEeqnarray}

    Constraints \eqref{eq:gpsp.00} are the flow conservation, which can also be obtained by subtracting  \eqref{eq:psp.2} from  \eqref{eq:psp.1} in the BSP for complete inter-hub networks. In the online companion, we provide more details on the method and its computational performance. 
\end{remark}

\subsection{Cut Generation Improvement}
\label{sec:closed-form}

There are two main issues with the cut generation procedure described in \autoref{sec:Feasibility}. First, while condition \eqref{eq:fcond} might be violated for more than one hub-hub connection $(h,k)$, we only add one feasibility cut. This is because BSP and DSP can not be decomposed on inter-hub link $(h,k)$. Moreover, solving LP subproblems can be time-consuming, as many of these suproblems must be solved within the search tree (see \autoref{sec:exp.alg}).
To overcome these challenges, in the following theorem, we show how to exploit the structure of the DSP to obtain an unbounded ray $\paren{\overline{\BFlambda}, \overline{\BFmu}, \overline{\BFnu}}$ for each pair of hubs for which condition \eqref{eq:fcond} is violated.

\begin{theorem}
    \label{thm:fcut}
    Let $\paren{\overline{\BFx}, \overline{\BFy}}$ be a feasible solution to the BMP.
    Let $\hat{h}, \hat{k} \in H$ be an arbitrary pair of hubs for which condition \eqref{eq:fcond} is violated.
    Then, given a constant $\varGamma > 0$, the vector $\paren{\overline{\lambda}, \overline{\mu}, \overline{\nu}}$ with
    \begin{IEEEeqnarray}{rCll}
        \IEEEyesnumber\label{eq:fcut.ray}
        \IEEEyessubnumber*
        \overline{\lambda}_{hk} &=&\begin{cases}
            -\varGamma & \text{if } h=\hat{h} \text{ and } k=\hat{k}, \\
            0          & \text{otherwise},
        \end{cases}\ & h,k \in H, \label{eq:fcut.ray.1}\\
        \overline{\mu}_{hi} &=& \sum_{h'\in H}\sum_{k\in H} \lambda_{h'k}\,\overline{x}_{h'i} - \sum_{k\in H} \lambda_{hk}\,\overline{x}_{hi}
        \ &  h\in H, i \in N, \IEEEeqnarraynumspace\label{eq:fcut.ray.2}\\
        \overline{\nu}_{ki} &=& - \sum_{h \in H} \lambda_{hk}\,\overline{x}_{hi} \ & k\in H, i \in N, \label{eq:fcut.ray.3}
    \end{IEEEeqnarray}
    is an unbounded ray for DSP$\paren{\overline{\BFx}, \overline{\BFy}}$.
\end{theorem}
\proof{Proof.}
Given $\hat{h}, \hat{k} \in H$, we first show that the vector $\paren{\overline{\lambda}, \overline{\mu}, \overline{\nu}}$ is feasible.
Substituting the values in left-hand side of constraint \eqref{eq:dsp.1}, we get

\begin{IEEEeqnarray*}{l}
    \lambda_{hk} + \mu_{hi} + \nu_{ki} = \begin{cases}
        -\varGamma\,\paren{1-\overline{x}_{\hat{h}i}} & \text{if } h=\hat{h} \text{ and } k=\hat{k},       \\
        -\varGamma\,\overline{x}_{\hat{h}i}           & \text{if } h\neq\hat{h} \text{ and } k\neq\hat{k}, \\
        0                                             & \text{otherwise},
    \end{cases}\  h, k\in H, i \in N,
\end{IEEEeqnarray*}
which is always lest than or equal to 0 as $-\varGamma  < 0$ and $\overline{x}_{hi} \leq 1$ holds for any $h, k\in H, i \in N$.

We now show that the objective function is unbounded over $\paren{\overline{\lambda}, \overline{\mu}, \overline{\nu}}$. Substituting the values in expression \eqref{eq:dsp.0} yields
\begin{IEEEeqnarray*}{l}
    \sum_{h \in H}\sum_{k \in H} Q_{hk}\overline{y}_{hk}\, \lambda_{hk} + \sum_{h \in H}\sum_{i \in N} \paren{ O_i \overline{x}_{hi}\, \mu_{hi} + \sum_{j \in N} w_{ij} \overline{x}_{hj}\, \nu_{hi}} \\
    \qquad = -\varGamma\,Q_{\hat{h}\hat{k}}\overline{y}_{\hat{h}\hat{k}} - \sum_{h:h\neq\hat{h}}\sum_{i\in N} \varGamma\, \overline{x}_{\hat{h}i} O_i \overline{x}_{hi} + \sum_{i\in N} \sum_{j\in N} \varGamma\, w_{ij} \overline{x}_{\hat{h}i}\overline{x}_{\hat{k}j} \\
    \qquad = -\varGamma\,Q_{\hat{h}\hat{k}}\overline{y}_{\hat{h}\hat{k}} - 0 + \varGamma\, \sum_{i\in N} \overline{x}_{\hat{h}i} \sum_{j\in N} w_{ij}\overline{x}_{\hat{k}j} \\
    \qquad = \varGamma\,\paren{\sum_{i \in N}z_{i\hat{h}\hat{k}} - Q_{\hat{h}\hat{k}}\overline{y}_{\hat{h}\hat{k}}}.
\end{IEEEeqnarray*}
Since $\hat{h}, \hat{k} \in H$ violate condition \eqref{eq:fcond}, we have $\sum_{i \in N}z_{i\hat{h}\hat{k}} - Q_{\hat{h}\hat{k}}\overline{y}_{\hat{h}\hat{k}} > 0$.
Therefore, the maximization problem \eqref{eq:dsp} becomes unbounded as $\varGamma \to \infty$.
\Halmos
\endproof

Using \autoref{thm:fcut}, we can generate the feasibility cut \eqref{eq:feasibilityCut} without solving the DSP subproblems. More importantly, it implies that we can generate a feasibility cut whenever we detect a pair of hubs $(\hat{h}, \hat{k})$ on which there exists an insufficient number of primary vehicles.
If there are more than one infeasible inter-hub link, one can choose an arbitrary pair or select an $(\hat{h}, \hat{k})$ equal to $\arg\max_{(h,k)\in H\times H} \paren{\sum_{i\in N} \sum_{j\in N} w_{ij}\overline{x}_{\hat{h}i}\overline{x}_{\hat{k}j} - Q_{hk}\overline{y}_{hk}}$, i.e., the pair for which the highest amount of capacity violation is observed.
Therefore, we consider a modified branch-and-cut framework, where at each node, we find all $(\hat{h}, \hat{k})$  with infeasible flows, calculate the extreme rays  $\paren{\overline{\BFlambda}, \overline{\BFmu}, \overline{\BFnu}}$ using \autoref{thm:fcut}, and add the resulting cuts to prune that node (if needed).
Our computational results show that the multi-cut approach outperforms the single-cut version, specially for the problem under multiple demand scenarios (see \autoref{sec:uncertainty}).

\subsection{Valid Inequalities}
\label{sec:valid-inequality}
\label{sec:model.valid}

Initially, the BMP has poor information on the flows over the inter-hub links as it includes $\BFy$ variables, but not their relation to the flow variables $\BFz$. Therefore, the initial bounds are usually loose, and the algorithm may go through many iterations to obtain some information through feasibility cuts.
In a desire to give the algorithm a better warm start and improve its linear relaxation, we can exploit the model's structure to generate some valid inequalities.
In particular, we can set a lower bound on the total number of vehicles that arrive to and dispatch from a hub based on the total incoming and outgoing demand assigned to that hub.
Let $Q^\text{in}_h = \max_k \brac{Q_{kh}}$ and $Q^\text{out}_h = \max_k \brac{Q_{hk}}$. Then, the following set of inequalities provide the aforementioned bounds.
\begin{IEEEeqnarray}{rCll}
    \IEEEyesnumber\label{eq:valids1}
    \IEEEyessubnumber*
    \frac{1}{Q^\text{in}_h}\, \sum_{i \in N} D_i\, x_{hi} &\leq& \sum_{k\in H} y_{kh},\ & h\in H \label{eq:valid1} \\
    \frac{1}{Q^\text{out}_h} \sum_{i\in N} O_i\, x_{hi} &\leq& \sum_{k\in H} y_{hk},\ & h\in H. \label{eq:valid2}
\end{IEEEeqnarray}

Moreover, when there is a flow between two nodes and both nodes are selected as hubs, they must be connected by at least one primary vehicle. The following valid inequality calculates the minimum number of vehicles required to travel between two hubs based on their shipment volume and primary vehicle capacity.
\begin{IEEEeqnarray}{l}
    \IEEEyessubnumber*
    \ceil{\frac{w_{hk}}{Q_{hk}}}\,\paren{x_{hh} + x_{kk} - 1} \leq y_{hk},\ h, k\in H. \label{eq:valid3}
\end{IEEEeqnarray}

\section{Demand Uncertainty}
\label{sec:uncertainty}

The HNDPv problem presented in Section \ref{sec:problem} assumes that  the OD demands are known in the planning stage. In reality, however, shipment volumes are stochastic, and long-term deterministic forecasts are unreliable. In this section, we address the demand uncertainty in HNDPv under a stochastic programming framework. The goal is to account for demand uncertainty in the design phase of the network in order to maintain the operational reliability of the network when the actual demand is realized.

For each $i,j \in N$, let random variable $w_{ij}(\xi)$ represent the  flow that needs to be sent from node $i$ to node $j$, where $\xi\in \Xi$ for a given support $\Xi$.
Define $O_i(\xi)=\sum_{j\in N}w_{ij}(\xi)$ and $D_i(\xi)=\sum_{j\in N}w_{ji}(\xi)$ as random variables representing the total outgoing flow from and the total incoming flow to node $i$, respectively.
We consider a two-stage stochastic program with recourse in which the location and allocation variables are dealt with in the first stage, while the flow variables and the required number of vehicles are determined in the second stage. The two-stage stochastic formulation of HNDPv is given as 
\begin{IEEEeqnarray}{llll}
     & \min_{\BFx \in \BFcalX}\ & \IEEEeqnarraymulticol{2}{l}{\mathit{LC}(\BFx) + \mathbb{E}_{\xi}\brak{P(\BFx,\xi)}} \label{eq:2shlp.1}
\end{IEEEeqnarray}
where $\mathbb{E}_{\xi}$ denotes the mathematical expectation with respect to $\xi \in \Xi$, 
\begin{IEEEeqnarray}{rlll}
    P(\BFx, \xi) = \min_{\BFy\in \BFcalY} \brac{\mathit{DC}_{\xi}(\BFx) + \mathit{HC}(\BFy) \left\lvert \sum_{i\in N}\sum_{j\in N} w_{ij}(\xi)\,x_{hi} x_{kj} \leq Q_{hk}\, y_{hk},\ h, k\in H \right.},
\end{IEEEeqnarray}
and $\mathit{DC}_{\xi}(\BFx)$ is the cost of direct access to non-hub node~$i$ which is defined as a function of random variables calculated as:
\begin{IEEEeqnarray}{rCl}
    \mathit{DC}_\xi(x) &=& \sum_{h\in H} \sum_{i\in N} \max\brac{\ceil{\frac{D_i(\xi)}{q_h}}, \ceil{\frac{O_i(\xi)}{q_h}}}\, c^\pm_{hi}\, x_{hi}. \label{eq:vcost-dir2}
\end{IEEEeqnarray}

Evaluating $\mathbb{E}_{\xi}\brak{P(\BFx,\xi)}$ in \eqref{eq:2shlp.1} is difficult and makes the optimization problem intractable. Therefore, following other works on stochastic hub location \citep[see, for example, ][]{Alumur+Nickel-TraResParB:Met-2012,Rostami+Kammerling-EurJOpeRes-2021}, we assume that the random variable $\xi$ follows a discrete distribution with finite support $S =\{s_1,\ldots, s_m\}$, where each event $s\in S$ occurs with probability $P(\xi=s)= p_s$. Accordingly, we use $w_{ij}^s$ to denote the amount of flow from node $i$ to node $j$ for each scenario $s\in S$. Therefore, $O_i^s=\sum_{j\in N}w_{ij}^{s}$ and $D_i^s=\sum_{j\in N}w_{ji}^{s}$ represent the total outgoing flow from and the total incoming flow to node~$i$, respectively.
Since vehicle selection decisions are to be done once scenario $s$ is realized, we redefine $\BFy$ variables as $y^{s}_{hk}$ to indicate the number of primary vehicles traveling on hub arc $(h,k)\in H\times H$ under scenario $s$. Moreover, for each scenario $s\in S$, we redefine flow variables as $z_{ihk}^s$ to determine the amount of node $i$'s demand transported from hub $h$ to hub $k$, i.e., $z_{ihk}^s=x_{hi}\sum_{j} w_{ij}^s x_{kj}, \forall i\in N, h,k\in H$.
Then, the deterministic equivalent formulation of $\text{P}$ is stated as:
\begin{IEEEeqnarray}{llll}
    \IEEEyesnumber\label{eq:shlp}
    \IEEEyessubnumber*
    \text{DEP}:\quad & \min \ & \IEEEeqnarraymulticol{2}{l}{\mathit{LC}(\BFx) + \sum_{s\in S} p_s \paren{\mathit{DC}_{s}(\BFx) + \mathit{HC}_{s}(\BFy)}} \label{eq:shlp.0} \\ [1ex]
    & \text{s.t.}\ & \sum_{i\in N} z^s_{ihk} \leq Q_{hk}\, {y}^s_{hk},\  & h, k\in H, s\in S\IEEEeqnarraynumspace \label{eq:shlp.y} \\
    && \sum_{k\in H} z^s_{ihk} = O^s_i\, {x}_{hi},\  & h\in H, i \in N, s\in S\label{eq:shlp.2}\\
    && \sum_{h\in H} z^s_{ihk} = \sum_{j\in N} w^s_{ij}\, {x}_{kj},\  & k\in H, i \in N, s\in S\label{eq:shlp.3}\\
    && \BFx \in \BFcalX,\ \BFz^s \in\mathbb{R}^{\card{N}\times\card{H}\times\card{H}}_{\geq 0},\ \BFy^s\in \BFcalY, s\in S \label{eq:shlp.00}
\end{IEEEeqnarray}
where,
\begin{IEEEeqnarray}{rCl}
    \mathit{DC}_s(\BFx) &=& \sum_{h\in H} \sum_{i\in N} \max\brac{\ceil{\frac{D_i^s}{q_h}}, \ceil{\frac{O_i^s}{q_h}}}\, c^\pm_{hi}\, x_{hi}, \label{eq:vcost-dir3} \\
    \mathit{HC}_s(\BFy) &=& \sum_{h\in H} \sum_{k\in H} C_{hk}\, y^s_{hk}.
\end{IEEEeqnarray}

\subsection{Solving the HNDPv Under Demand Uncertainty}
\label{sec:uncertainty.bd}

The Benders decomposition approach presented in \autoref{sec:sol} can be easily adjusted to solve DEP.
For each scenario $s\in S$, and for any feasible solution $\paren{\overline{\BFx}, \overline{\BFy}}$ to the master problem at a given iteration of the algorithm, we define one scenario-based Benders subproblem (S-BSP) as follows.
\begin{IEEEeqnarray}{llll}
    \IEEEyesnumber\label{eq:spsp}
    \IEEEyessubnumber*
    \text{S-BSP}_s\paren{\overline{\BFx}, \overline{\BFy}}:\qquad & \min\ & \IEEEeqnarraymulticol{2}{l}{0} \label{eq:spsp.0} \\ [1ex]
    & \text{s.t.}\ & \sum_{i\in N} z^s_{ihk} \leq Q_{hk}\, \overline{y}^s_{hk},\  & h, k\in H, \label{eq:spsp.1}\\
    && \sum_{k\in H} z^s_{ihk} = O^s_i\, \overline{x}_{hi},\  & h\in H, i \in N, \label{eq:spsp.2}\\
    && \sum_{h\in H} z^s_{ihk} = \sum_{j\in N} w^s_{ij}\, \overline{x}_{kj},\  & k\in H, i \in N, \label{eq:spsp.3}\\
    && \BFz^s \in\mathbb{R}^{\card{N}\times\card{H}\times\card{H}}_{\geq 0}.\label{eq:spsp.00}
\end{IEEEeqnarray}

Following the approach described in \autoref{sec:sol}, one can solve the dual of $\text{S-BSP}_s\paren{\overline{\BFx}, \overline{\BFy}}$ to generate feasibility cuts and apply them to the master problem whenever an infeasible S-BSP is observed. The Benders master problem corresponding to the S-HNDPv (i.e., the S-BMP) is formulated as:
\begin{IEEEeqnarray}{llll}
    \IEEEyesnumber\label{eq:sbmp}
    \IEEEyessubnumber*
    \text{S-BMP}:\qquad & \min\ & \IEEEeqnarraymulticol{2}{l}{\mathit{LC}(\BFx) + \sum_{s\in S} p_s \paren{\mathit{DC}_{s}(\BFx) + \mathit{HC}_{s}(\BFy)}} \label{eq:sbmp.0} \\ [1ex]
    & \text{s.t.}\ & \Omega^{sr} \leq 0,\ r \in R^s,  s\in S,\\
    & & \BFx \in \BFcalX,\ \BFy^s\in \BFcalY, s\in S, \label{eq:sbmp.00}
\end{IEEEeqnarray}
where
\begin{IEEEeqnarray}{l}
    \Omega^{sr} = \sum_{h \in H}\sum_{k \in H} Q_{hk}\overline{\lambda}^{sr}_{hk}\,y^s_{hk} + \sum_{h \in H}\sum_{i \in N} \paren{ O^s_i \overline{\mu}^{sr}_{hi}\, x_{hi}  + \sum_{j \in N} w^s_{ij}\overline{\nu}^{sr}_{hi}\, x_{hj}}, \IEEEeqnarraynumspace\label{eq:sfcut}
\end{IEEEeqnarray}
and $\overline{\lambda}^{s\tau}_{hk}$, $\overline{\mu}^{s\tau}_{hi}$, and $\overline{\mu}^{s\tau}_{hi}$ being the unbounded rays in the set extreme rays $\mathcal{R}^s$ corresponding to constraint \eqref{eq:spsp.1}, \eqref{eq:spsp.2}, and \eqref{eq:spsp.3}, respectively, in the dual space.




Note that the results of \autoref{thm:fcut} are still valid for S-BSP$_s\paren{\overline{\BFx}, \overline{\BFy}}$ for a given scenario $s\in S$. Therefore, similar to the deterministic case, we may add multiple cuts for any scenario $s$ with infeasible inter-hub flows. 
\section{Computational Experiments}
\label{sec:experiments}

In this section, we present the numerical results evaluating models and solution algorithms.
The algorithms are coded in Python. We used Guropi optimizer v9.5 \citep{gurobi} and its callback features to solve our optimization models. Experiments are conducted on a Linux laptop with Intel\textsuperscript{\textregistered} Core\textsuperscript{\texttrademark} i9-11900H CPU @ 2.50GHz and 32GB of RAM using up to 14 threads. All experiments are done within a time limit of 12,000 seconds.

In the following sections, we first present the test instances and then provide the experimental results of solving deterministic and stochastic HNDPv test instances.

\subsection{Problem Test Instances and Experimental Design}
\label{sec:inst}

To perform our experiments, we use the Australian Post (AP) dataset used by \cite{Contreras+Diaz-ORSpe-2009}.
We consider location cost and capacity values provided in the AP dataset. The location costs are all assumed to be tight, while both tight (T) and loose (L) values are tested for hub capacity levels. An uncapacitated case (U) is also considered where the hubs have unrestricted capacity.
We assume fixed capacities $Q_{hk} = Q, \forall h,k\in H$, and $q_h = q, \forall h\in H$, and fixed unit transportation costs $b_{ij} = b$ and $B_{ij} = B$, $\forall (i,j)\in A$, for all primary and secondary vehicles, respectively.
Inspired by \cite{Tanash+Contreras-ComOpeRes-2017}, we set all vehicle fixed costs to 0 and consider the following configurations for primary and secondary vehicle parameters:
\begin{itemize}
    \item L1$\colon\paren{Q=600, q=100, B=600, b=260} \Rightarrow \frac{B/Q}{b/q} \approx 0.38$
    \item L2$\colon\paren{Q=600, q=150, B=600, b=300} \Rightarrow \frac{B/Q}{b/q} = 0.50$
    \item L3$\colon\paren{Q=320, q=100, B=500, b=260} \Rightarrow \frac{B/Q}{b/q} \approx 0.60$
    \item L4$\colon\paren{Q=320, q=150, B=500, b=300} \Rightarrow \frac{B/Q}{b/q} \approx 0.78$
\end{itemize}
These configurations are chosen such that the unit transportation cost on hub arcs, when fully utilized, is smaller than the unit transportation cost on access arcs.
We also indicate the $\frac{B/Q}{b/q}$ ratio for L1 to L4. This value is considered as the smallest discount factor that can be achieved on hub arcs \citep{Tanash+Contreras-ComOpeRes-2017}.
Instances with $\card{N}$ equal to 20, 25, 40, 50, 100, and 200 are considered for our computational tests.
We refer to each instance by \#\textsubscript{1}\#\textsubscript{2}-\#\textsubscript{3} notation, where \#\textsubscript{1} denotes the number of nodes in the instance (e.g., 50), \#\textsubscript{2} indicates the capacity configuration (i.e., T, L, or U), and \#\textsubscript{3} shows the vehicle configuration (e.g., L1).

\subsection{Algorithmic Efficiency}

\begin{figure}
    \FIGURE
    {\includegraphics[width=0.85\textwidth]{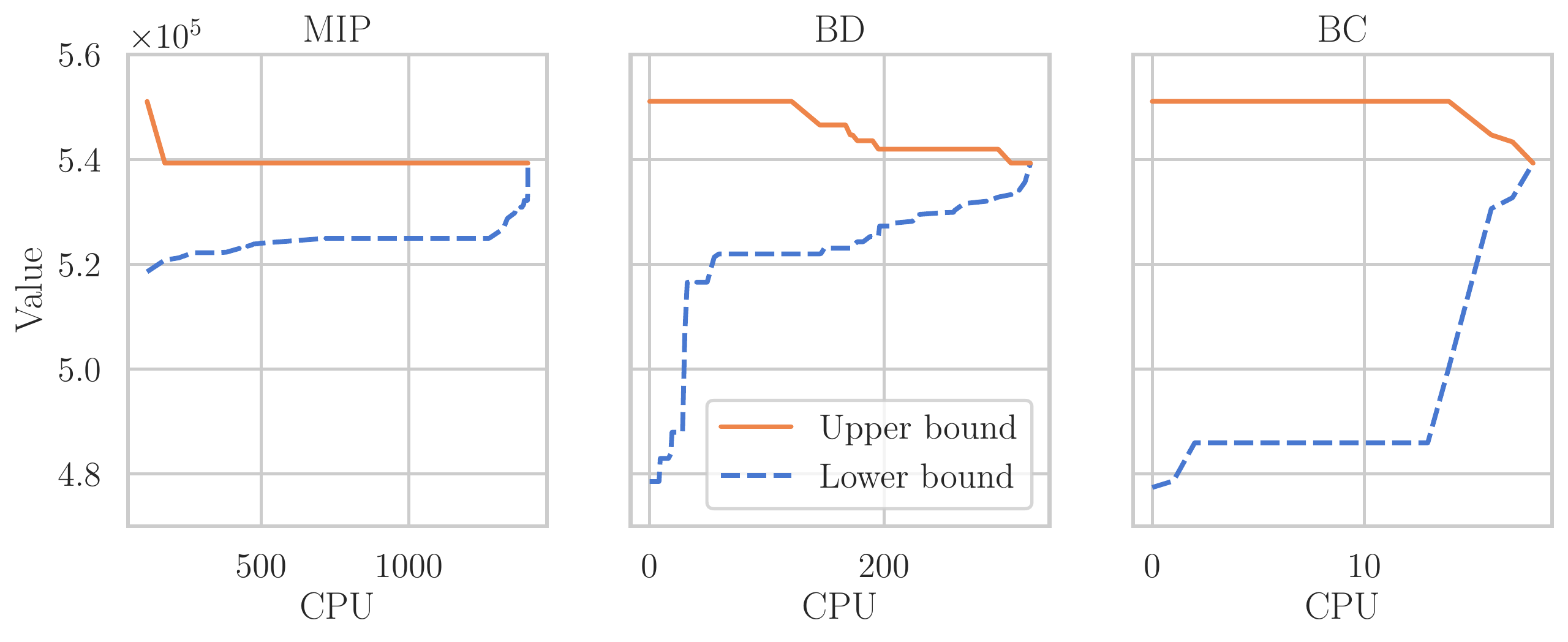}}
    {Closing optimality gap by different approaches (problem instances: 50L-L4).\label{fig:conv-gaps.50}}
    {}
\end{figure}

\begin{figure}
    \FIGURE
    {\includegraphics[width=0.85\textwidth]{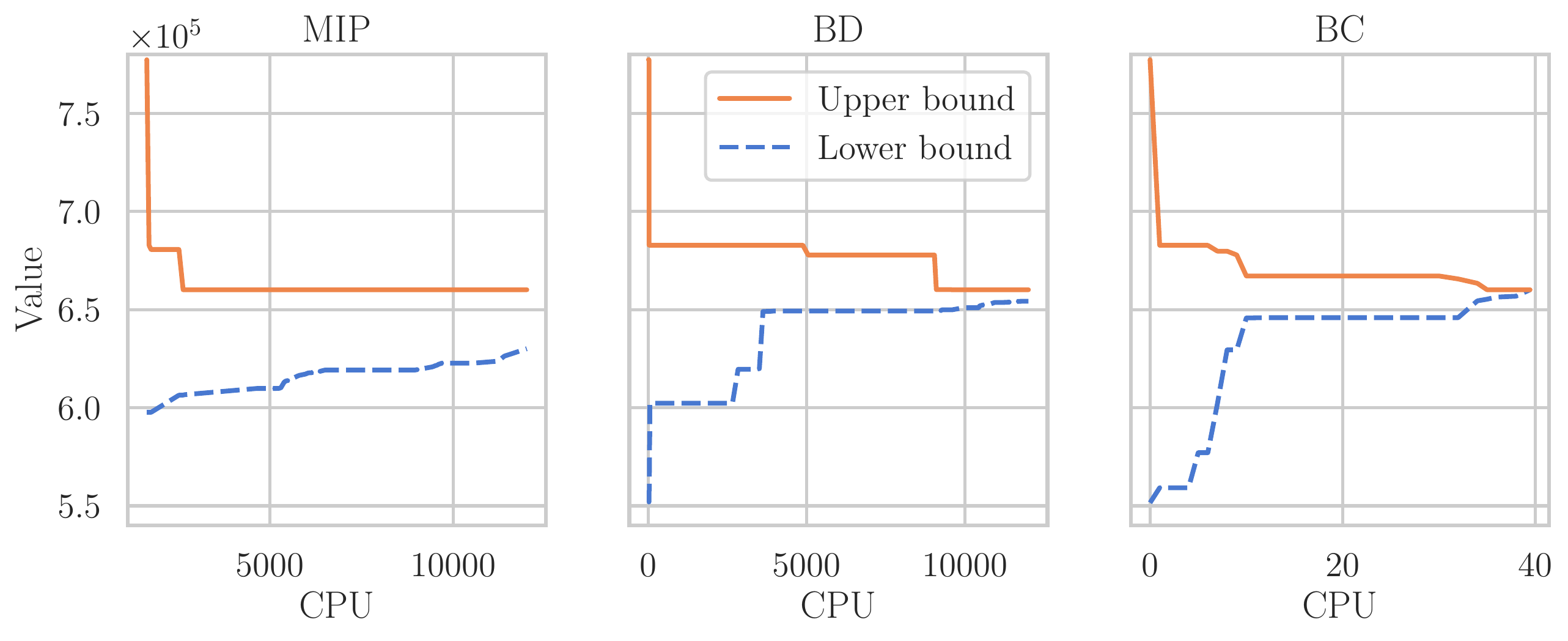}}
    {Closing optimality gap by different approaches (problem instances: 75L-L4).\label{fig:conv-gaps.75}}
    {}
\end{figure}

\begin{figure}
    \FIGURE
    {\includegraphics[width=0.7\textwidth]{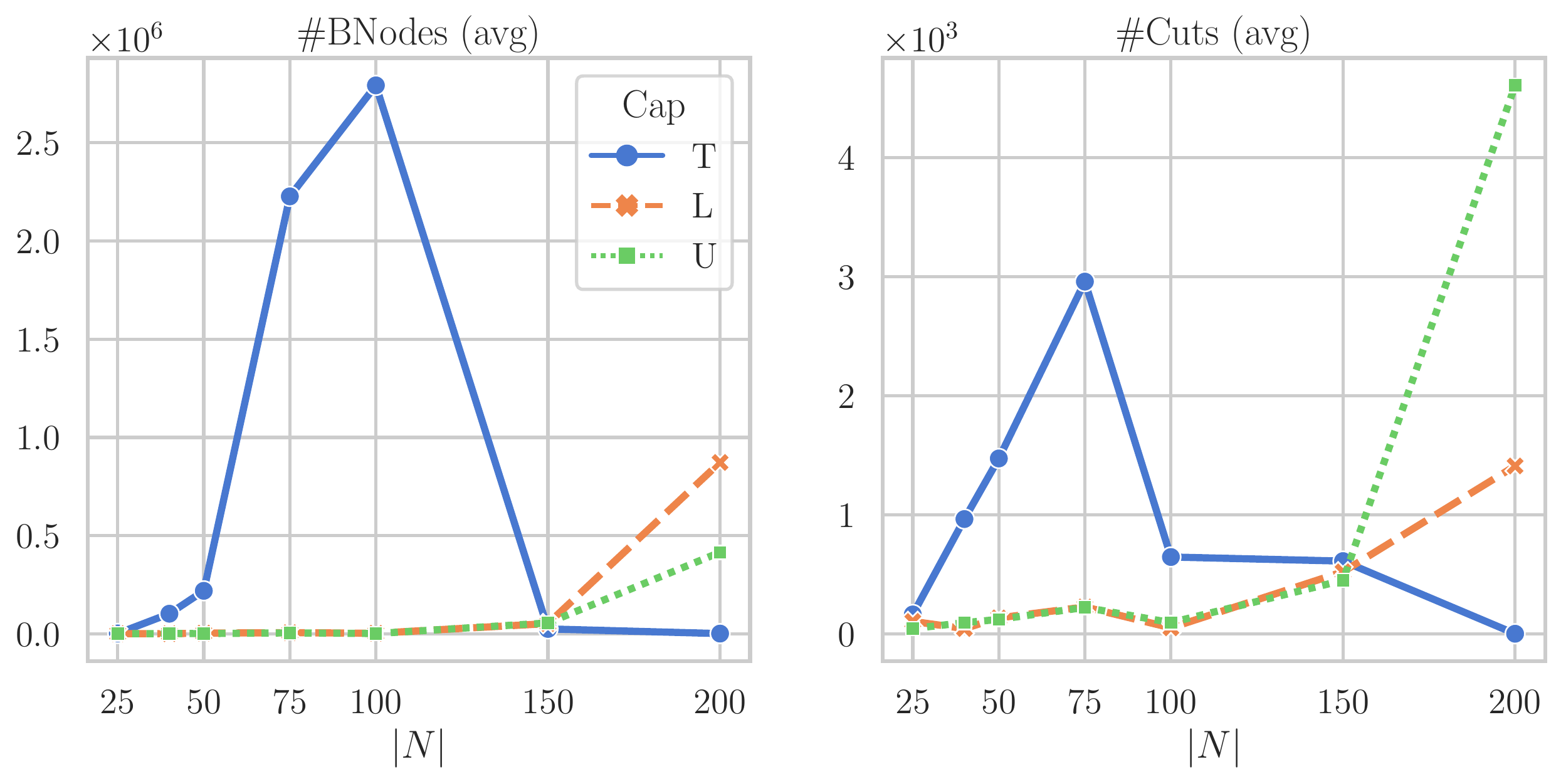}}
    {Effect of hub capacity configurations on the number of explored branch-and-bound nodes and the number of generated feasibility cuts.\label{fig:perf.cap}}
    {}
\end{figure}

\begin{figure}
    \FIGURE
    {\includegraphics[width=0.7\textwidth]{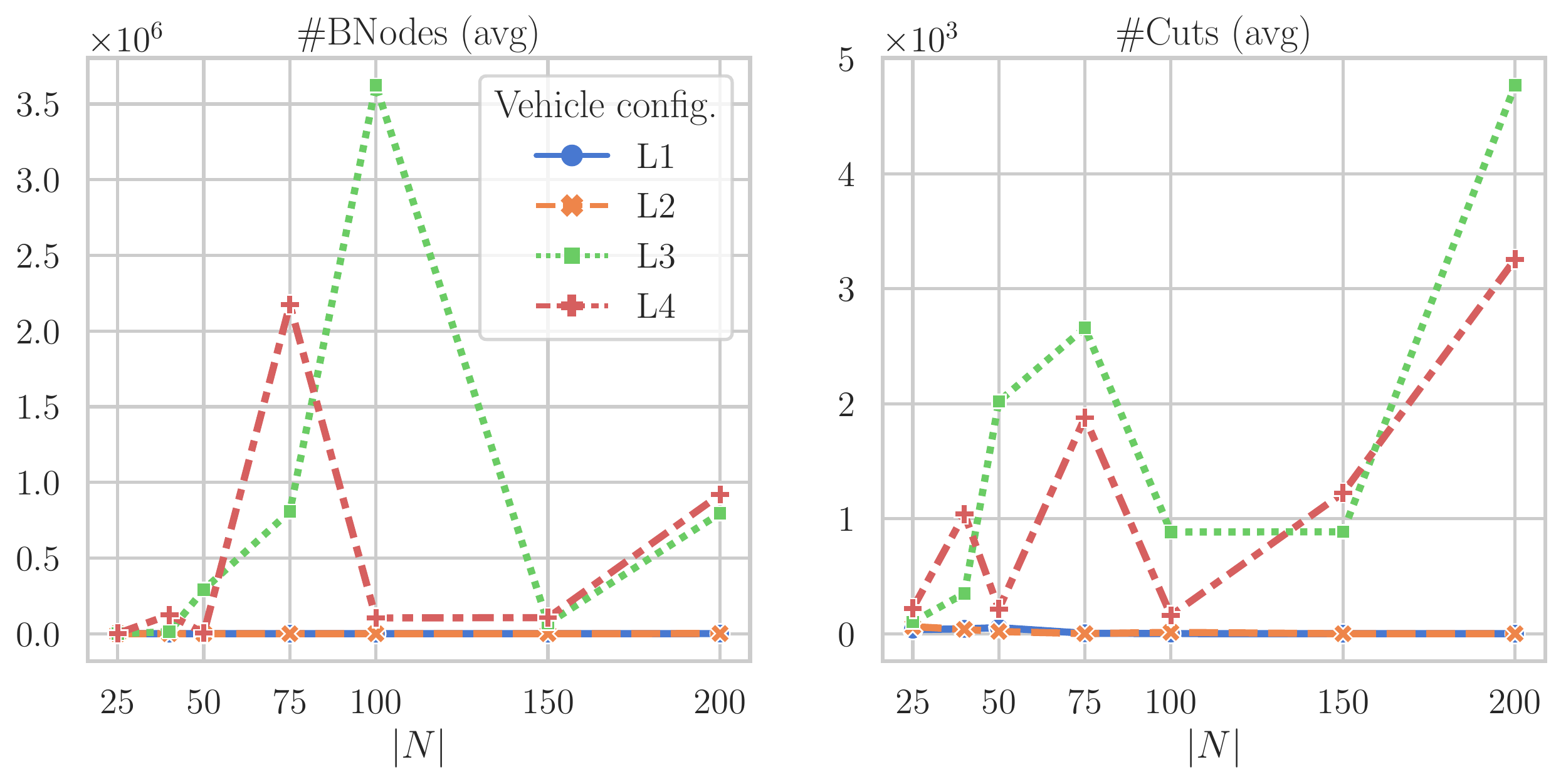}}
    {Effect of vehicle configurations on the number of explored branch-and-bound nodes and the number of generated feasibility cuts.\label{fig:perf.veh}}
    {}
\end{figure}

We evaluate the performance of the proposed algorithms in terms of efficiency and compare them with Gurobi applied directly to solve the MIP formulation P in \eqref{model:SM}. We show by MIP the direct application of Gurobi to solve the MIP model, by BD the Benders-based branch-and-cut algorithm with single feasibility cut, and by BC, the Benders-based branch-and-cut algorithm with multiple feasibility cuts derived from Theorem \ref{thm:fcut}. First, we compare the performance of these algorithms on two different instances in Figures \ref{fig:conv-gaps.50} and \ref{fig:conv-gaps.75}. We then compare MIP and BC on small to medium-size instances in \autoref{tab:perf-mipvbc}, and show the performance of the BC on large-scale instances in \autoref{tab:perf-comp-summary}, and Figures \ref{fig:perf.cap} and \ref{fig:perf.veh}. These tables and figures have been designed to give a full picture of the algorithms' efficiency. The detailed performance of the algorithms can be found in the online supplement of the paper.
\label{sec:exp.alg}


Figures \ref{fig:conv-gaps.50} and \ref{fig:conv-gaps.75} illustrate how the upper bound and the lower bound converge during the procedure of solving instances 50L-L4 and 75L-L4. In each figure, the $x$ access shows the CPU time in seconds, while the $y$ access shows the values for lower and upper bounds.
To solve the 50L-L4, Gurobi takes 23 minutes, the BD takes 5 minutes, while the BC only takes 18 seconds. Although the initial lower bound in MIP is better, and the solver is able to find a good upper bound after some time, it takes a long time to close the gap. BD and BC start with worse lower bounds (since less information is available at the beginning), but they close the gap significantly faster. When the number of nodes is increased to 75 (\autoref{fig:conv-gaps.75}), we see the same behavior for the upper and lower bounds. However, neither MIP nor BD  was able to solve the instance within the given time limit of 12,000 seconds, while the BC was able to find the optimal solution in less than 40 seconds.

\begin{table}
    \TABLE
    {Comparison of the BD and BC performances.\label{tab:bd-vs-bc}}
    {\scriptsize \begin{tabular}{@{}lrrrrrrrrrr@{}}
            \toprule
            \multicolumn{1}{c}{\multirow{2}{*}{Inst}} & \multicolumn{5}{c}{BD}       & \multicolumn{5}{c}{BC}                                                                                                                                                                                                                                                                                              \\ \cmidrule(lr){2-6} \cmidrule(l){7-11}
            \multicolumn{1}{c}{}                      & \multicolumn{1}{c}{\#BNodes} & \multicolumn{1}{c}{\#Cuts} & \multicolumn{1}{c}{\#Calls} & \multicolumn{1}{c}{CPU\textsuperscript{c}} & \multicolumn{1}{c}{CPU (\%Gap)} & \multicolumn{1}{c}{\#BNodes} & \multicolumn{1}{c}{\#Cuts} & \multicolumn{1}{c}{\#Calls} & \multicolumn{1}{c}{CPU\textsuperscript{c}} & \multicolumn{1}{c@{}}{CPU (\%Gap)} \\ \midrule
            50L-L4                                    & {7,561}                      & 118                        & 125                         & 317.53                                     & 324.61 (0.0)                    & {5,002}                      & 259                        & 119                         & 2.76                                       & 17.97 (0.0)                        \\
            75L-L4                                    & {5,239}                      & 168                        & 175                         & \textgreater{}12,000                       & \textgreater{}12,000 (0.88)     & {10,230}                     & 476                        & 220                         & 9.81                                       & 39.47 (0.0)                        \\ \bottomrule
        \end{tabular}}
    {}
\end{table}

\autoref{tab:bd-vs-bc} report more information on the performance of  the BD and BC in solving the same two instances. For each instance, CPU\textsuperscript{c} shows the time spent to solve the feasibility subproblems and generating the associated cuts, \#BNodes, \#Cuts, \#Calls, and \%Gap show the number of explored branch-and-bound nodes, the number of generated feasibility cuts, the number of time the solver has called the subproblem, and the percent relative optimality gap, respectively.
For 50L-L4, while both call the feasibility check subproblems almost the same, the BC generated more than twice the cuts than the BD. However, the BC generates the cuts 114 times faster than the BD.  For the 75L-L4, the number of feasibility cuts generated by BC is larger, but the total time spent to generate these cuts is significantly smaller. We can see that most of the time in the BD is spent in finding the coefficients in the feasibility cuts. \autoref{thm:fcut}, on the other hand, allows us to skip solving an optimization problem and add multiple cuts at once. That is, generating cuts in this way allows us to solve larger problem instances more efficiently. Therefore, in the rest of this section, we do not report the results of the BD.

\begin{table}
    \TABLE
    {CPU values of MIP and BC in solving small and medium-size instances.\label{tab:perf-mipvbc}}
    {\small \begin{tabular}{@{}llrrrrrrrr@{}}
            \toprule
            \multicolumn{1}{c}{\multirow{2}{*}{$\card{N}$}} & \multicolumn{1}{c}{\multirow{2}{*}{Cap\textsuperscript{$\dagger$}}} & \multicolumn{4}{c}{MIP} & \multicolumn{4}{c}{BC}                                                                                                                                                       \\ \cmidrule(lr){3-6} \cmidrule(l){7-10}
            \multicolumn{1}{c}{}                            & \multicolumn{1}{c}{}                                                & \multicolumn{1}{c}{L1}  & \multicolumn{1}{c}{L2} & \multicolumn{1}{c}{L3} & \multicolumn{1}{c}{L4} & \multicolumn{1}{c}{L1} & \multicolumn{1}{c}{L2} & \multicolumn{1}{c}{L3} & \multicolumn{1}{c}{L4} \\ \midrule
            \multirow{3}{*}{20}                             & T                                                                   & 7.27                    & 3.69                   & 2.97                   & 9.18                   & 0.15                   & 0.36                   & 0.29                   & 0.74                   \\
                                                            & L                                                                   & 9.88                    & 6.52                   & 2.76                   & 1.99                   & 0.40                   & 0.27                   & 0.37                   & 0.44                   \\
                                                            & U                                                                   & 8.28                    & 8.11                   & 3.41                   & 1.99                   & 0.20                   & 0.32                   & 0.24                   & 0.33                   \\ \midrule
            \multirow{3}{*}{25}                             & T                                                                   & 16.53                   & 14.84                  & 17.79                  & 12.88                  & 0.27                   & 0.24                   & 2.71                   & 11.28                  \\
                                                            & L                                                                   & 21.01                   & 24.55                  & 10.91                  & 10.18                  & 1.50                   & 1.36                   & 1.41                   & 1.12                   \\
                                                            & U                                                                   & 18.41                   & 9.40                   & 4.72                   & 2.21                   & 0.69                   & 0.68                   & 0.31                   & 1.31                   \\ \midrule
            \multirow{3}{*}{40}                             & T                                                                   & 266.52                  & 219.05                 & 270.53                 & 304.10                 & 0.49                   & 0.53                   & 44.16                  & 474.52                 \\
                                                            & L                                                                   & 98.27                   & 44.47                  & 87.70                  & 42.07                  & 1.08                   & 3.12                   & 1.68                   & 1.70                   \\
                                                            & U                                                                   & 110.98                  & 41.58                  & 71.28                  & 44.68                  & 2.70                   & 1.65                   & 1.39                   & 3.47                   \\ \midrule
            \multirow{3}{*}{50}                             & T                                                                   & 9,765.28                & 6,691.78               & 5,218.47               & 4,488.73               & 2.50                   & 1.41                   & 1,327.72               & 11.68                  \\
                                                            & L                                                                   & 734.93                  & 1,459.41               & 1,286.41               & 1,402.19               & 1.72                   & 0.80                   & 7.49                   & 17.97                  \\
                                                            & U                                                                   & 915.23                  & 988.54                 & 766.35                 & 1,238.07               & 2.82                   & 1.66                   & 8.26                   & 4.74                   \\ \midrule
            \multicolumn{2}{l}{Average}                     & 997.72                                                              & 792.66                  & 645.28                 & 629.86                 & 1.21                   & 1.03                   & 116.34                 & 44.11                                           \\ \bottomrule
        \end{tabular}
    }
    {{\textsuperscript{$\dagger$}} Hub capacity configuration.}
\end{table}

\begin{table}
    \TABLE
    {CPU and (\%Gap) values of BC in solving large instances.\label{tab:perf-comp-summary}}
    {\small \begin{tabular}{@{}llrrrr@{}}
            \toprule
            \multicolumn{1}{c}{\multirow{2}{*}{$\card{N}$}} &
            \multicolumn{1}{c}{\multirow{2}{*}{Cap}}        &
            \multicolumn{4}{c}{Vehicle configuration}                                                                                                      \\ \cmidrule(l){3-6}
            \multicolumn{1}{c}{}                            &
            \multicolumn{1}{c}{}                            &
            \multicolumn{1}{c}{L1}                          &
            \multicolumn{1}{c}{L2}                          &
            \multicolumn{1}{c}{L3}                          &
            \multicolumn{1}{c}{L4}                                                                                                                         \\ \midrule
            \multirow{3}{*}{75}                             & T            & 2.15         & 2.34            & (2.71)\textsuperscript{$\dagger$} & (1.96)   \\
                                                            & L            & 2.21         & 2.19            & 81.24                             & 39.47    \\
                                                            & U            & 1.54         & 2.23            & 50.77                             & 40.26    \\ \midrule
            \multirow{3}{*}{100}                            & T            & 8.31         & 5.57            & (0.60)                            & 143.39   \\
                                                            & L            & 3.14         & 6.26            & 29.75                             & 11.95    \\
                                                            & U            & 2.35         & 2.86            & 76.16                             & 28.71    \\ \midrule
            \multirow{3}{*}{150}                            & T            & 24.81        & 23.49           & 520.12                            & 399.36   \\
                                                            & L            & 9.40         & 14.61           & 722.13                            & 615.61   \\
                                                            & U            & 6.47         & 12.79           & 619.59                            & 421.49   \\ \midrule
            \multirow{3}{*}{200}                            & T            & 109.19       & 106.40          & 80.79                             & 104.07   \\
                                                            & L            & 83.46        & 106.87          & (0.24)                            & 7,395.43 \\
                                                            & U            & 56.82        & 70.47           & 7,689.79                          & (0.09)   \\ \midrule
            \multicolumn{2}{l}{Average}                     & 25.82 (0.00) & 29.67 (0.00) & 3,822.53 (1.18) & 2,766.65 (1.03)                              \\ \bottomrule
        \end{tabular}%
    }
    {{\textsuperscript{$\dagger$}} Values in parentheses represent percentage optimality gaps.}
\end{table}

Next, in \autoref{tab:perf-mipvbc} we compare BC and MIP performances in solving small to medium-size problem instances for which both approaches could solve them to optimality within the time limit. For each instance, with size $\card{N}$, the hub capacity configuration (Cap), and the vehicle configurations L1 to L4, the table reports CPU times for each solution method.  For both methods, the instances with tight capacity are the most difficult to solve. From the vehicle configuration perspective, the L3 and L4 are the most difficult instances for the BC, while this is not the case for the MIP. This will be investigated further \autoref{tab:perf-comp-summary}. Overall, as can be seen, the BC was considerably faster than MIP, particularly when solving larger problems. For instances with $\card{N}=75$ and more, the MIP could not close the gap within the time limit. Therefore, in \autoref{tab:perf-comp-summary}, we only report the performance of the BC in solving larger problem instances.

 \autoref{tab:perf-comp-summary} lists the CPU time for each instance under different hub capacities and vehicle configurations. In cases where the time limit is reached, instead of the time, the \%Gap is reported in parentheses.
Out of 48 instances, the BC was able to find the optimal solution for 43 instances. These instances took, on average, about 20 seconds for $\card{N}=75$, 30 seconds for $\card{N}=100$, 5 minutes for $\card{N}=150$,  and 30 minutes when $\card{N}=200$ to be solved by BC. The optimality gap for unsolved instances is reported as 1.12\%.
We observe that in larger instances, the vehicle configurations have a more significant effect on the BC performance. The last row of \autoref{tab:perf-comp-summary} shows the average CPU time and \%Gap values for different vehicle configurations. Vehicle configurations L3 and L4 are the most difficult settings to deal with. The reason is that the primary vehicle capacities are smaller in these settings.
Hence, more effort should be made to ensure the feasibility of the flows on the inter-hub links, and more feasibility cuts are generated (see the online supplement for more details).
The average CPU time and \%Gap values of the L3 configuration are the highest among all vehicle configurations.
The reason is that the L3 configuration has not only the smallest primary vehicle capacity, but also the smallest secondary capacity. Therefore, more vehicles are required on the access-level network. This increases the cost of assignment decisions, and, as a result, more exploration is required to find the optimal location and allocation decisions.

In Figures \ref{fig:perf.cap} and \ref{fig:perf.veh}, we show the effect of instance characteristics on the number of branch-and-bound nodes explored by the algorithm (\#BNodes) and the number of generated feasibility cuts (\#Cuts).
\autoref{fig:perf.cap} illustrates the average values for different hub capacity configurations and instance sizes. For instances with size 100 or smaller, the algorithm explores significantly more \#BNodes and outputs more \#Cuts when hub capacities are tight.
We observed that in the AP instances with $\card{N} \geq 150$, although there exists a large number of OD pairs, shipment volumes are smaller than those in smaller instances.
Therefore, when hub capacities get large, more assignment options become available for the demand nodes, and more demands can get consolidated at hubs. This leads to more consolidated flows on the inter-hub links, which in turn requires more cuts to ensure feasibility on those links.


We illustrate the effect of different vehicle configurations on \#BNodes and \#Cuts in \autoref{fig:perf.veh}. We can observe that many feasibility cuts are generated for the instances with small primary vehicles (i.e., L3 and L4).
When primary vehicles have larger capacities, the introduced valid inequalities, along with other cuts generated by Gurobi, help to obtain feasible flows on the inter-hub links with no or only a few feasibility cuts. The online companion of this paper provides more details about the BC performance.

\subsection{Managerial Insights}
\label{sec:res.insights}

The main objective of the HNDPv is to optimize vehicle utilization to reduce cost.  In this section, we first investigate how such utilization is compared to conventional HLP with a constant discount factor. Then, we explore the effect of problem instance characteristics, i.e., size, hub capacity, and vehicle configurations on different cost factors and operational decisions.  While solving the conventional HLP with a constant discount factor, we assume that the transfer (discount), collection, and distribution factors are $0.75$, $3$, and $2$, respectively, as given in the original AP dataset. After solving each instance of the HLP and finding the location and allocation decisions, we assign the minimum number of vehicles to transport flows on the resulting network. We use the same vehicle configurations as in the HNDPv.

\autoref{tab:vals-veh-vs-disc} compares the HNDPv and the conventional HLP solutions for instances with 
$\card{N}=25$, $40$, and $50$ in terms of percentage difference of the following factors: the actual total cost (TC), the number of utilized primary and secondary vehicles (\#Veh1 and \#Veh2, respectively), and the average capacity utilization of the primary vehicles in percent (\%VUtil).
We use $100\times\frac{\text{HFlow}}{\text{\#Veh1} \times Q}$ to compute \%VUtil, where HFlow is the total flow on the inter-hub links.
As the HLP solution opens only one hub in uncapacitated instances (and hence no need for inter-hub vehicles), we only report the results for the capacitated instances.

\begin{table}[htbp]
    \TABLE
    {Comparison of the HLP-based solutions and the HNDPv solutions.\label{tab:vals-veh-vs-disc}}
    {\footnotesize \begin{tabular}{@{}lrrrr@{}}
            \toprule
            \multicolumn{1}{@{}c}{Inst} & \multicolumn{1}{c}{TC (\% diff.)} & \multicolumn{1}{c}{\#Veh1 (\% diff.)} & \multicolumn{1}{c}{\#Veh2 (\% diff.)} & \multicolumn{1}{c@{}}{\%VUtil (\% diff.)} \\ \midrule
            25T-L1                      & $ +0.62 $                         & 0                                     & $ +2.00 $                             & $ -0.27 $                                 \\
            25T-L2                      & $ +0.28 $                         & 0                                     & 0                                     & $ -0.27 $                                 \\
            25T-L3                      & $ +1.88 $                         & $ +12.50 $                            & $ +2.00 $                             & $ -10.91 $                                \\
            25T-L4                      & $ +1.65 $                         & $ +12.50 $                            & 0                                     & $ -10.91 $                                \\ \midrule
            25L-L1                      & $ +1.39 $                         & $ +33.33 $                            & $ +1.96 $                             & $ -25.36 $                                \\
            25L-L2                      & $ +1.89 $                         & $ +33.33 $                            & 0                                     & $ -21.93 $                                \\
            25L-L3                      & $ +0.25 $                         & 0                                     & 0                                     & $ -7.28 $                                 \\
            25L-L4                      & $ +0.32 $                         & 0                                     & 0                                     & $ -7.28 $                                 \\ \midrule
            Average ($\card{N}=25$)     & $ +1.04 $                         & $ +11.46 $                            & $ +0.75 $                             & $ -10.53 $                                \\ \midrule
            40T-L1                      & 0                                 & 0                                     & 0                                     & 0                                         \\
            40T-L2                      & $ +0.25 $                         & 0                                     & 0                                     & $ -4.35 $                                 \\
            40T-L3                      & $ +0.70 $                         & $ +12.50 $                            & 0                                     & $ -11.62 $                                \\
            40T-L4                      & $ +1.17 $                         & $ +12.50 $                            & 0                                     & $ -12.20 $                                \\ \midrule
            40L-L1                      & $ +1.00 $                         & $ +33.33 $                            & 0                                     & $ -22.13 $                                \\
            40L-L2                      & $ +1.17 $                         & $ +33.33 $                            & 0                                     & $ -22.79 $                                \\
            40L-L3                      & $ +0.01 $                         & 0                                     & 0                                     & $ -0.76 $                                 \\
            40L-L4                      & $ +0.02 $                         & 0                                     & 0                                     & $ -0.76 $                                 \\ \midrule
            Average ($\card{N}=40$)     & $ +0.54 $                         & $ +11.46 $                            & $ +0.00 $                             & $ -9.33 $                                 \\ \midrule
            50T-L1                      & $ +7.10 $                         & $ +100.00 $                           & $ -1.61 $                             & $ -44.90 $                                \\
            50T-L2                      & $ +7.44 $                         & $ +100.00 $                           & $ -1.79 $                             & $ -40.12 $                                \\
            50T-L3                      & $ +6.51 $                         & $ +87.50 $                            & $ -1.61 $                             & $ -40.10 $                                \\
            50T-L4                      & $ +6.45 $                         & $ +87.50 $                            & $ -1.79 $                             & $ -36.10 $                                \\ \midrule
            50L-L1                      & $ +0.31 $                         & 0                                     & 0                                     & $ -16.13 $                                \\
            50L-L2                      & $ +0.34 $                         & 0                                     & 0                                     & $ -16.13 $                                \\
            50L-L3                      & $ +0.30 $                         & 0                                     & 0                                     & $ -16.13 $                                \\
            50L-L4                      & $ +0.33 $                         & 0                                     & 0                                     & $ -16.13 $                                \\ \midrule
            Average ($\card{N}=50$)     & $ +3.60 $                         & $ +46.88 $                            & $ -0.85 $                             & $ -28.22 $                                \\ \bottomrule
        \end{tabular}
    }
    {}
\end{table}

\begin{figure}
    \FIGURE
    {\includegraphics[width=0.7\textwidth]{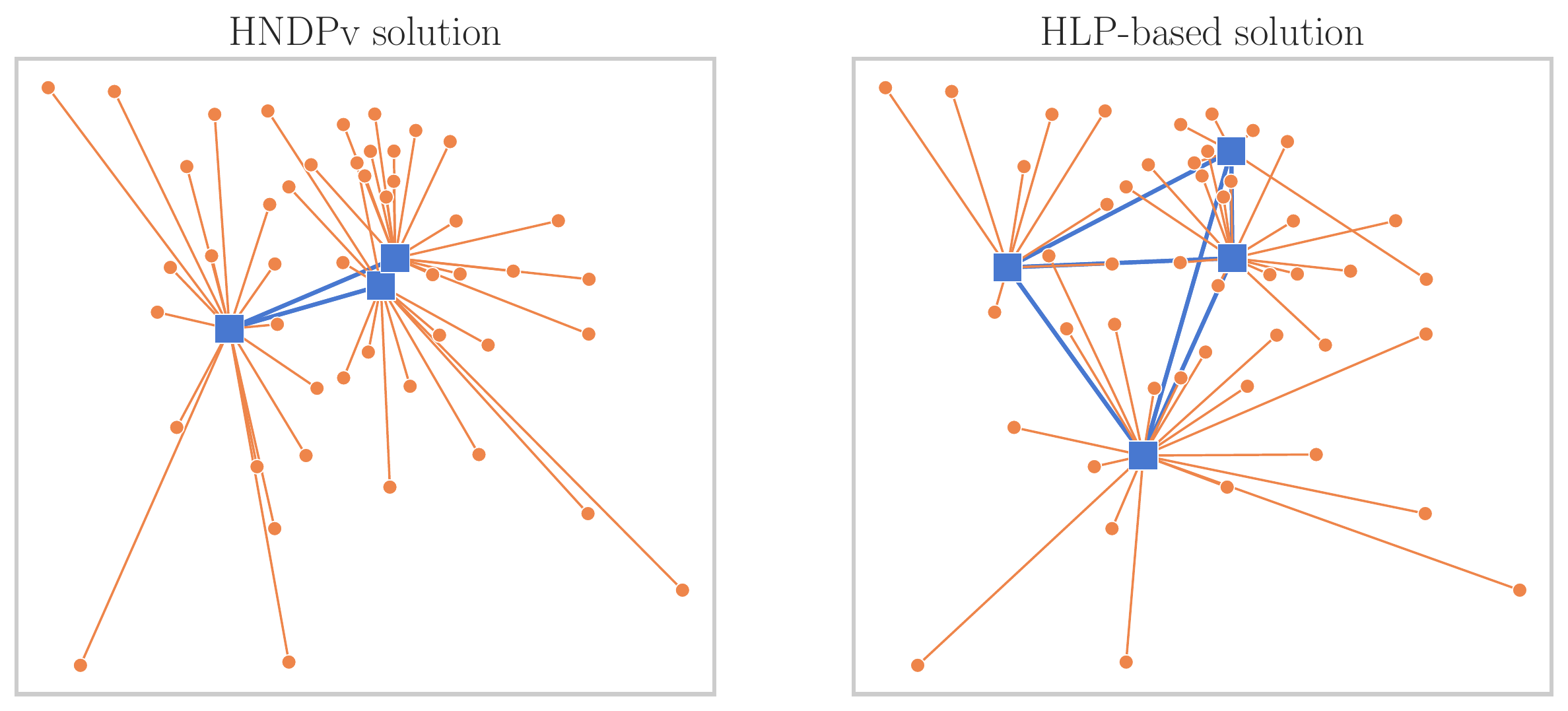}}
    {Illustration of the solutions to problem instance 50T-L1.\label{fig:sol.50t}}
    {Square shapes represent open hubs, small discs show the demand nodes.}
\end{figure}

\begin{figure}
    \FIGURE
    {\includegraphics[width=0.7\textwidth]{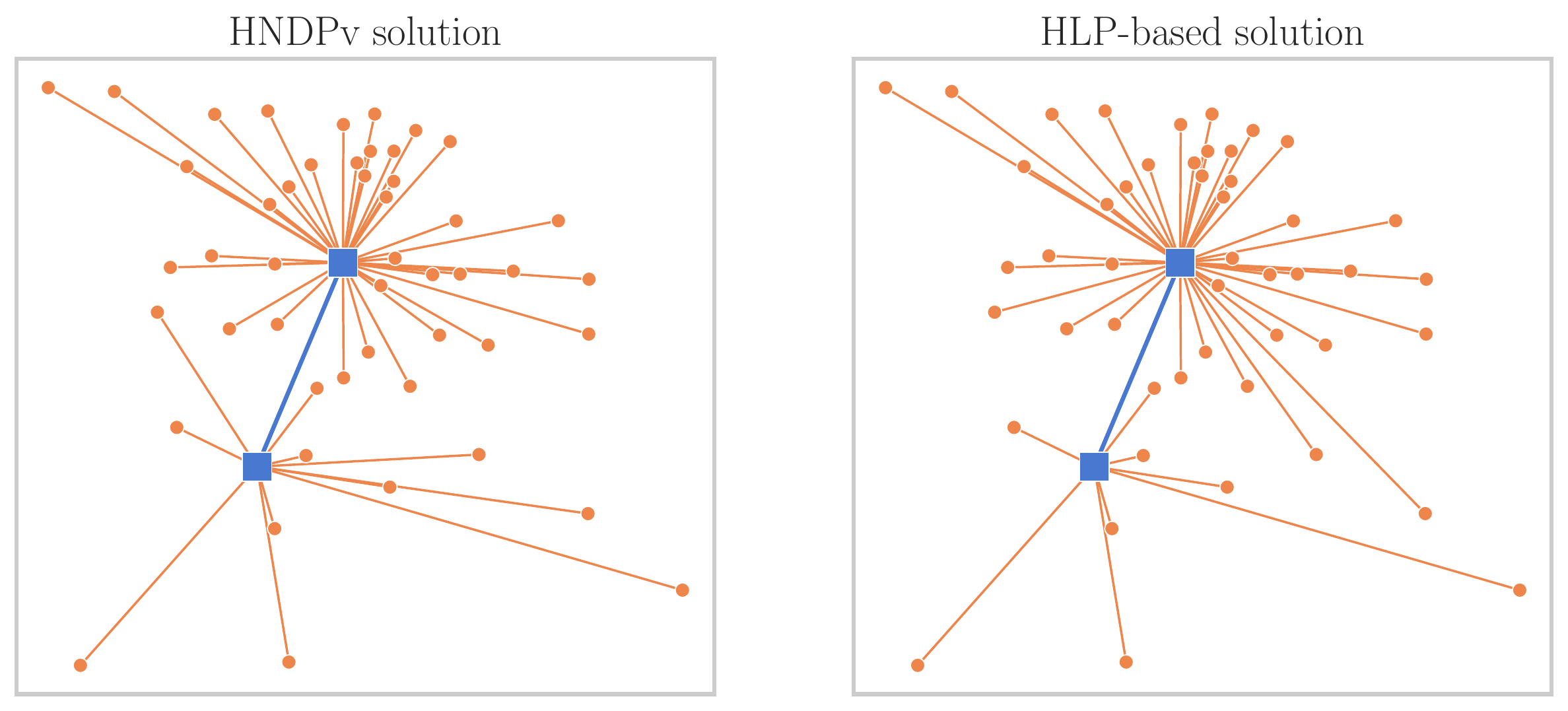}}
    {Illustration of the solutions to problem instance 50L-L1.\label{fig:sol.50l}}
    {}
\end{figure}

As can be seen, there is only one instance (i.e., 40T-L1) where both problems provided the same solution for. On the rest of the instances, the HLP-based solutions lead to, on average,  $1.04\%$, $0.54\%$, and $3.6\%$ higher costs for instances with 25, 40, and 50 nodes, respectively. 
All HLP-based solutions used the same or a larger number of primary vehicles on the network.
In 50-node instances with tight capacities, the HLP required up to double the primary vehicle fleet size. The reason is that the HLP solution opened more hubs than the HNDPv solution. Therefore, although slightly fewer secondary vehicles were used, more primary vehicles were required, leading to poor vehicle utilization.
\autoref{fig:sol.50t} illustrates the HNDPv solution (left) and the HLP-based solution (right). The HNDPv solution opens one less hub and constructs a different hub topology, resulting in a less costly solution with better utilized vehicles at the hub level.
Even when \#Veh1 is the same for both, the HNDPv solution provides a better primary vehicle utilization, as it incorporates vehicle-based decisions in the network design process, which can lead to a different location or allocation decisions (see \autoref{fig:sol.50l}).
The HNDPv solution was able to provide $10.53\%$, $9.33\%$, and $28.22\%$ better primary vehicle utilization, for instances with 25, 40, and 50 nodes, respectively.

\begin{figure}
    \FIGURE
    {\frame{\includegraphics[width=0.48\textwidth]{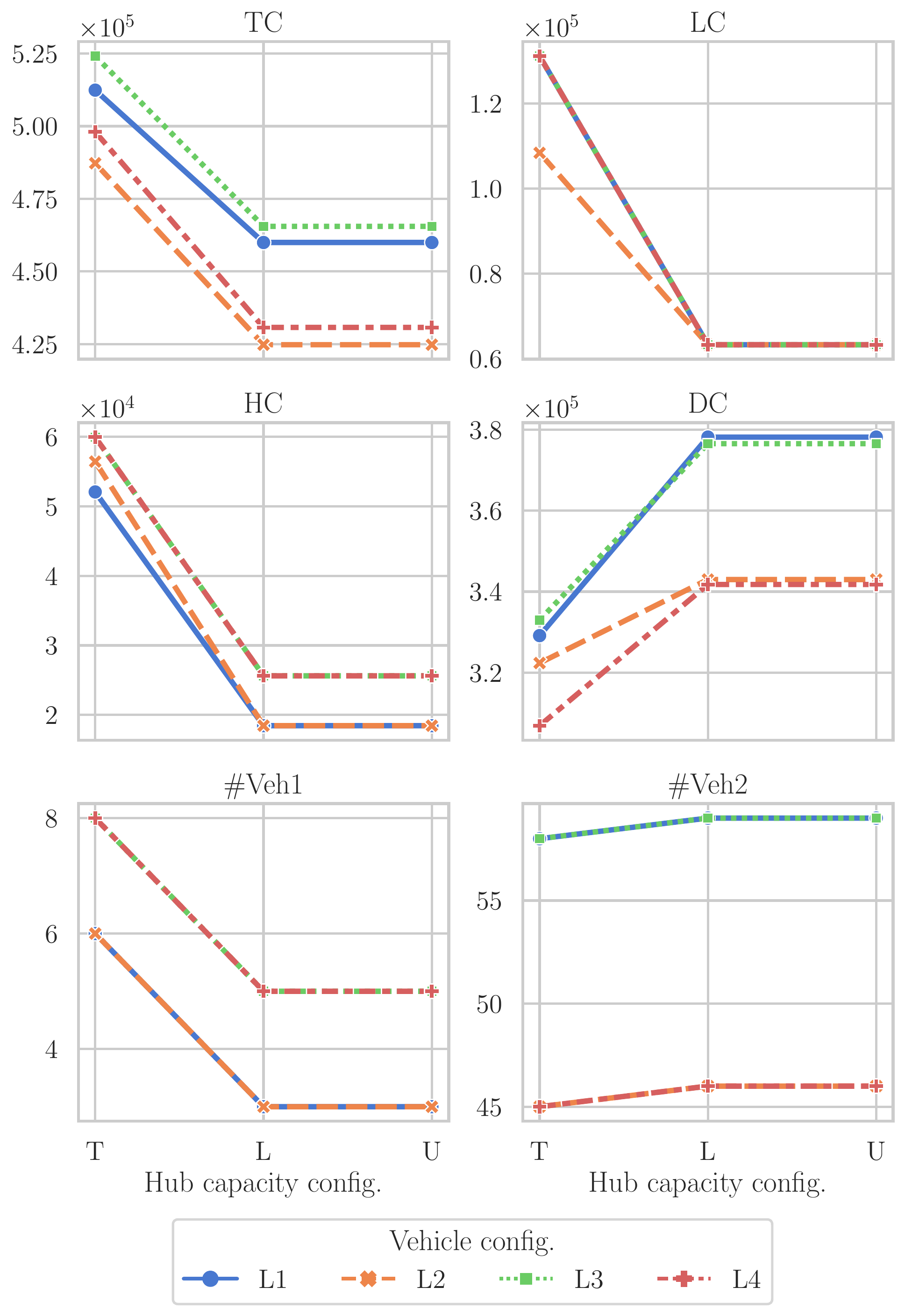}} \frame{\includegraphics[width=0.48\textwidth]{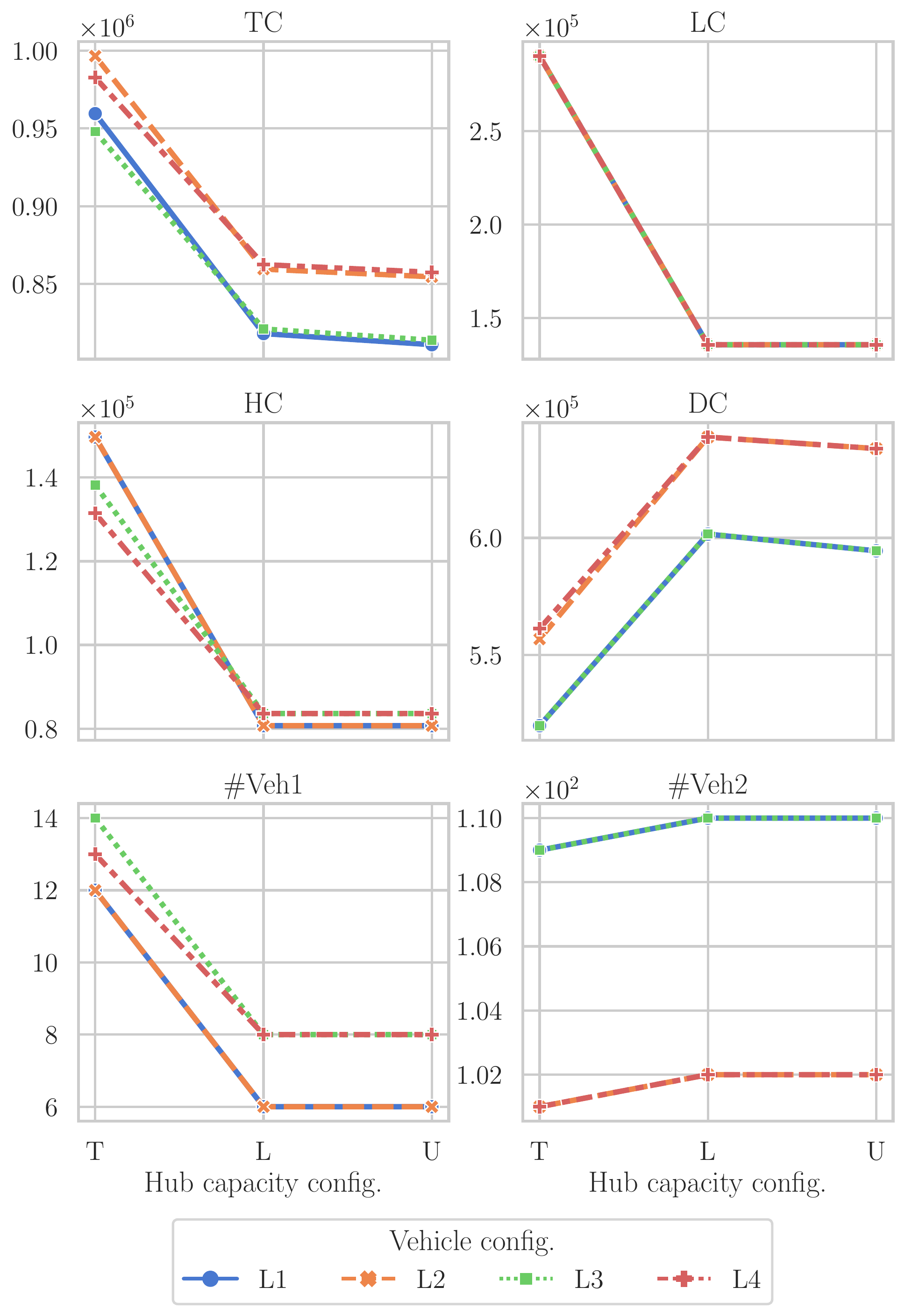}}}
    {Effect of different hub and vehicle configurations on the final solutions (Left: 40-node instances, Right: 100-node instances).\label{fig:configs.40}}
    {}
\end{figure}

We further investigate the effect of problem instance characteristics, i.e., size, hub capacity, and vehicle configurations on different cost factors and operational decisions.
Our aim is to identify which vehicle configuration is preferred under differ conditions.
We consider the criteria in our discussion: The objective function value or the total cost (TC), the total hub location cost (LC), the total transportation cost on the hub-level network (HC), the total transportation cost on the access-level network (DC), the number of selected hubs (\#Hubs), \#Veh1, \#Veh2, and \%VUtil.
Our experiments indicate that the solution is highly sensitive to the number of nodes in the instance.
In the AP dataset, smaller instances are created by aggregating the nodes in larger instances. As the total OD flows remain constant between all instances, smaller instances have larger OD shipment volumes, while larger instances have more fractional $w_{ij}$ values.
To see how such a characteristic may affect the solution, we refer to \autoref{fig:configs.40}, where the value of each criterion is plotted with respect to hub capacity and vehicle configurations.
\autoref{fig:configs.40} (left) illustrates solution characteristics for instances with $\card{N}=40$.
Here, when L2 and L4 are selected, a smaller total cost is incurred under all vehicle capacity (Cap) levels.
The L3 configuration leads to higher HC and DC values, as it offers the combination of the smallest vehicles in both hub and access-level networks. Therefore, more trips are required in both levels, leading to higher operational costs.
For the 40-node instances, using larger vehicles saves costs through consolidation, even though the operational cost for larger vehicles are higher.
Instances with larger set $N$, however, give different results. \autoref{fig:configs.40} (right) depicts similar criteria for instances with $\card{N}=100$.
Larger instances have more OD pairs with less shipment volumes. Therefore, consolidating flows are not as straightforward. TC is larger for vehicle configurations (VehConfs) that offer larger secondary vehicles (i.e., L2 and L4). Larger secondary vehicles are more costly to operate, and since flows are more fractional compared to the $\card{N}=40$ instances, we cannot fully benefit from the excess capacities on these vehicles. Therefore, even though less number of vehicles are used in L2 and L4, larger DC values are obtained.
Larger primary vehicles (as in L1 and L2) also leads to smaller \#Veh1.

For both 40 and 100-node instances, we observe that hub capacities have a direct impact on TC and LC. The higher the Cap level, the lower TC and LC values are, regardless of the VehConf choice.
Larger capacities also result in lower HC values and fewer number of primary vehicles, as more OD pairs can be linked through a single hub. However, DC values might increase as less number of hubs can be opened when the capacities are large.
Overall, tight hub capacities lead to higher total costs, resulting from higher location and inter-hub transportation costs. When capacities are restrictive, either more hubs are selected or larger but more expensive hubs are opened. This may increase the total distance travelled between the hubs, hence more inter-hub transportation cost is incurred.
When hub capacities are nonrestrictive, there is more flexibility in making location and allocation decisions, and this allows us to come up with a less costly solution.

\begin{table}
    \TABLE
    {Effect of instance size and vehicle configuration on the final solution.\label{tab:val-summary-veh}}
    {\small \begin{tabular}{@{}llrrrrrrrr@{}}
            \toprule
            \multicolumn{1}{@{}c}{$\card{N}$}                                  &
            \multicolumn{1}{c}{VehConf\textsuperscript{$\dagger$}}             &
            \multicolumn{1}{c}{TC (avg)}                                       &
            \multicolumn{1}{c}{LC (avg)}                                       &
            \multicolumn{1}{c}{HC (avg)}                                       &
            \multicolumn{1}{c}{DC (avg)}                                       &
            \multicolumn{1}{c}{\begin{tabular}{c} \#Hubs\\(avg) \end{tabular}} &
            \multicolumn{1}{c}{\begin{tabular}{c} \#Veh1\\(avg) \end{tabular}} &
            \multicolumn{1}{c}{\begin{tabular}{c} \#Veh2\\(avg) \end{tabular}} &
            \multicolumn{1}{c@{}}{\begin{tabular}{c@{}} \%VUtil\\(avg) \end{tabular}}                                              \\ \midrule
            \multirow{4}{*}{25}                                                & L1 & 4.42e+5 & 1.21e+5 & 3.57e+4 & 2.85e+5 & 2.33 & 4.00  & 50.67  & 74.37 \\
                                                                               & L2 & 3.93e+5 & 9.69e+4 & 2.41e+4 & 2.72e+5 & 2.00 & 3.00  & 37.67  & 69.07 \\
                                                                               & L3 & 4.47e+5 & 9.69e+4 & 3.04e+4 & 3.19e+5 & 2.00 & 4.33  & 51.67  & 91.35 \\
                                                                               & L4 & 3.98e+5 & 9.69e+4 & 3.04e+4 & 2.70e+5 & 2.00 & 4.33  & 37.67  & 91.35 \\ \hline
            \multirow{4}{*}{40}                                                & L1 & 4.77e+5 & 8.59e+4 & 2.97e+4 & 3.62e+5 & 2.33 & 4.00  & 58.67  & 74.06 \\
                                                                               & L2 & 4.46e+5 & 7.84e+4 & 3.11e+4 & 3.36e+5 & 2.33 & 4.00  & 45.67  & 75.45 \\
                                                                               & L3 & 4.85e+5 & 8.59e+4 & 3.71e+4 & 3.62e+5 & 2.33 & 6.00  & 58.67  & 92.07 \\
                                                                               & L4 & 4.53e+5 & 8.59e+4 & 3.71e+4 & 3.30e+5 & 2.33 & 6.00  & 45.67  & 92.26 \\ \hline
            \multirow{4}{*}{50}                                                & L1 & 5.31e+5 & 1.27e+5 & 3.14e+4 & 3.72e+5 & 2.33 & 3.67  & 62.67  & 77.66 \\
                                                                               & L2 & 5.57e+5 & 1.43e+5 & 3.93e+4 & 3.74e+5 & 2.33 & 3.67  & 56.67  & 75.93 \\
                                                                               & L3 & 5.39e+5 & 1.27e+5 & 4.14e+4 & 3.71e+5 & 2.33 & 5.67  & 62.67  & 89.50 \\
                                                                               & L4 & 5.66e+5 & 1.43e+5 & 4.92e+4 & 3.74e+5 & 2.33 & 5.67  & 56.67  & 87.62 \\ \hline
            \multirow{4}{*}{75}                                                & L1 & 6.18e+5 & 1.24e+5 & 5.81e+4 & 4.36e+5 & 3.00 & 6.00  & 86.00  & 60.25 \\
                                                                               & L2 & 6.64e+5 & 1.24e+5 & 5.81e+4 & 4.82e+5 & 3.00 & 6.00  & 80.00  & 60.25 \\
                                                                               & L3 & 6.29e+5 & 1.16e+5 & 5.91e+4 & 4.54e+5 & 3.00 & 8.67  & 86.00  & 83.03 \\
                                                                               & L4 & 6.74e+5 & 1.24e+5 & 6.63e+4 & 4.84e+5 & 3.00 & 8.67  & 80.00  & 76.22 \\ \hline
            \multirow{4}{*}{100}                                               & L1 & 8.63e+5 & 1.87e+5 & 1.04e+5 & 5.72e+5 & 3.33 & 8.00  & 109.67 & 45.82 \\
                                                                               & L2 & 9.04e+5 & 1.87e+5 & 1.04e+5 & 6.13e+5 & 3.33 & 8.00  & 101.67 & 45.70 \\
                                                                               & L3 & 8.61e+5 & 1.87e+5 & 1.02e+5 & 5.72e+5 & 3.33 & 10.00 & 109.67 & 66.39 \\
                                                                               & L4 & 9.01e+5 & 1.87e+5 & 9.96e+4 & 6.14e+5 & 3.33 & 9.67  & 101.67 & 68.07 \\ \hline
            \multirow{4}{*}{150}                                               & L1 & 9.31e+5 & 1.10e+5 & 1.09e+5 & 7.12e+5 & 4.00 & 12.00 & 158.00 & 34.43 \\
                                                                               & L2 & 1.02e+6 & 1.10e+5 & 1.09e+5 & 7.97e+5 & 4.00 & 12.00 & 151.00 & 34.42 \\
                                                                               & L3 & 9.20e+5 & 1.10e+5 & 9.50e+4 & 7.15e+5 & 4.00 & 12.67 & 158.00 & 60.43 \\
                                                                               & L4 & 1.01e+6 & 1.10e+5 & 9.50e+4 & 8.01e+5 & 4.00 & 12.67 & 151.00 & 60.79 \\ \hline
            \multirow{4}{*}{200}                                               & L1 & 1.22e+6 & 1.09e+5 & 1.56e+5 & 9.51e+5 & 4.33 & 14.67 & 206.67 & 28.89 \\
                                                                               & L2 & 1.34e+6 & 1.11e+5 & 1.61e+5 & 1.06e+6 & 4.33 & 14.67 & 199.67 & 29.76 \\
                                                                               & L3 & 1.19e+6 & 1.09e+5 & 1.32e+5 & 9.53e+5 & 4.33 & 15.00 & 206.67 & 52.63 \\
                                                                               & L4 & 1.31e+6 & 1.30e+5 & 1.75e+5 & 1.01e+6 & 5.00 & 20.00 & 199.00 & 41.28 \\ \hline
        \end{tabular}%
    }
    {{\textsuperscript{$\dagger$}} Vehicle configuration.}
\end{table}

In \autoref{tab:val-summary-veh}, we analyze the effect of vehicle configurations on the solution characteristics in more details.
Here, we see that if shipment volumes are larger, as seen in smaller instances, high-capacity secondary vehicles (i.e., in L2 and L4) helps to save costs, even though they are more expensive to operate.
For larger instances, where more trips of secondary vehicles are required, their operational costs dominate their capacity benefits. Therefore, on average, L2 and L4 become expensive choices.
Larger instances also require more primary vehicles to use, leading to lower \%VUtil values. Therefore, it is always best to go with an option that provides the most vehicle utilization for the operational costs we pay at the hub-level network. That is why in large instances with $\card{N}\geq 100$, the L3 configuration is the cost-efficient choice.
For the smaller instances with more aggregated flows, choosing larger primary vehicles helps to save costs.
VehConf has no significant effect on the number of selected hubs in most of the instances.
We observed that for small-size instances with $\card{N}= 25$ and large hub capacities, the solution tries to open more hubs when larger primary vehicles are available in order to exploit the economies of scale by consolidating demand on the inter-hub links and save costs on the access-level network.
Another exception is observed for the 200-node instance under L and U capacity settings. This problem requires many secondary vehicles. Therefore, when vehicle configuration is L4, one additional hub is opened to save on DC. This did not happen for L2 since the additional hub for L2 would require more primary vehicle trips which would be more expensive due to the higher operational cost of primary vehicles in L2 compared to L4.
More detailed numerical results are provided in the online supplement.

\subsection{HNDPv Instances With Stochastic Demands}
\label{sec:res.stoch}

To evaluate the BC in solving the stochastic version of the HNDPv, we use the dataset introduced in \cite{Rostami+Kammerling-EurJOpeRes-2021} with $\card{N}\in \brac{25,40,50,75}$. The demand value for an OD pair $(i,j)$ in a stochastic scenario is chosen from a Poisson distribution with event rate $\pi_{i} \pi_{j} w_{ij}$, where $w_{ij}$ are the demand values of the underlying AP dataset instance, and $\pi_{i}$  denotes the deviation from the base case being uniformly distributed in the interval $\brak{0.5, 1.5}$. Similar to \cite{Alumur+Nickel-TraResParB:Met-2012} and \cite{Rostami+Kammerling-EurJOpeRes-2021}, we consider five scenarios, each with a probability of occurrence of $0.2$.

\autoref{tab:perf-ss-summary} lists the CPU time for each instance under different hub capacity and vehicle configurations. In cases where the time limit is reached, the \%Gap is reported in parentheses. As expected, stochastic instances are more complex and difficult to solve than deterministic ones. Out of 48 instances, the BC was able to find the exact solution to 35 instances within the time limit.
Like the deterministic case, stochastic instances with tight hub capacities and small primary vehicle capacities are the most difficult ones to solve.
The average optimality gap for unsolved instances is reported as 3.31\%.
The online companion of the paper gives the computational details of solving the stochastic HNDPv instances and a discussion on how the stochastic solution compares to the solution to the expected value problem.

\begin{table}[htbp]
    \TABLE
    {CPU and (\%Gap) values of BC in solving stochastic instances.\label{tab:perf-ss-summary}}
    {\small \begin{tabular}{@{}llrrrr@{}}
            \toprule
            \multicolumn{1}{c}{\multirow{2}{*}{$\card{N}$}} &
            \multicolumn{1}{c}{\multirow{2}{*}{Cap}}        &
            \multicolumn{4}{c}{Vehicle configuration}                                                                                        \\ \cmidrule(l){3-6}
            \multicolumn{1}{c}{}                            &
            \multicolumn{1}{c}{}                            &
            \multicolumn{1}{c}{L1}                          &
            \multicolumn{1}{c}{L2}                          &
            \multicolumn{1}{c}{L3}                          &
            \multicolumn{1}{c}{L4}                                                                                                           \\ \midrule
            \multirow{3}{*}{25}                             & T               & 1.02            & 4.02            & 3.16            & 8.38   \\
                                                            & L               & 0.94            & 0.51            & 2.43            & 2.29   \\
                                                            & U               & 0.89            & 0.38            & 1.18            & 1.09   \\ \midrule
            \multirow{3}{*}{40}                             & T               & 69.26           & 4.06            & (1.98)          & (3.05) \\
                                                            & L               & 24.01           & 12.65           & 54.05           & 29.48  \\
                                                            & U               & 26.20           & 10.00           & 40.90           & 49.37  \\ \midrule
            \multirow{3}{*}{50}                             & T               & 84.02           & 85.17           & (1.41)          & (2.53) \\
                                                            & L               & 16.19           & 17.61           & 615.14          & 282.78 \\
                                                            & U               & 89.45           & 52.17           & 143.78          & 238.97 \\ \midrule
            \multirow{3}{*}{75}                             & T               & (4.46)          & (7.00)          & (5.03)          & (6.20) \\
                                                            & L               & 2,547.60        & (1.17)          & (1.21)          & (3.22) \\
                                                            & U               & 1,483.88        & 979.86          & (1.57)          & (1.81) \\ \midrule
            \multicolumn{2}{l}{Average}                     & 1,361.96 (4.46) & 2,097.20 (4.09) & 5,071.72 (2.24) & 5,051.03 (3.36)          \\ \bottomrule
        \end{tabular}
    }
    {}
\end{table}

\begin{figure}
    \FIGURE
    {\includegraphics[width=0.72\textwidth]{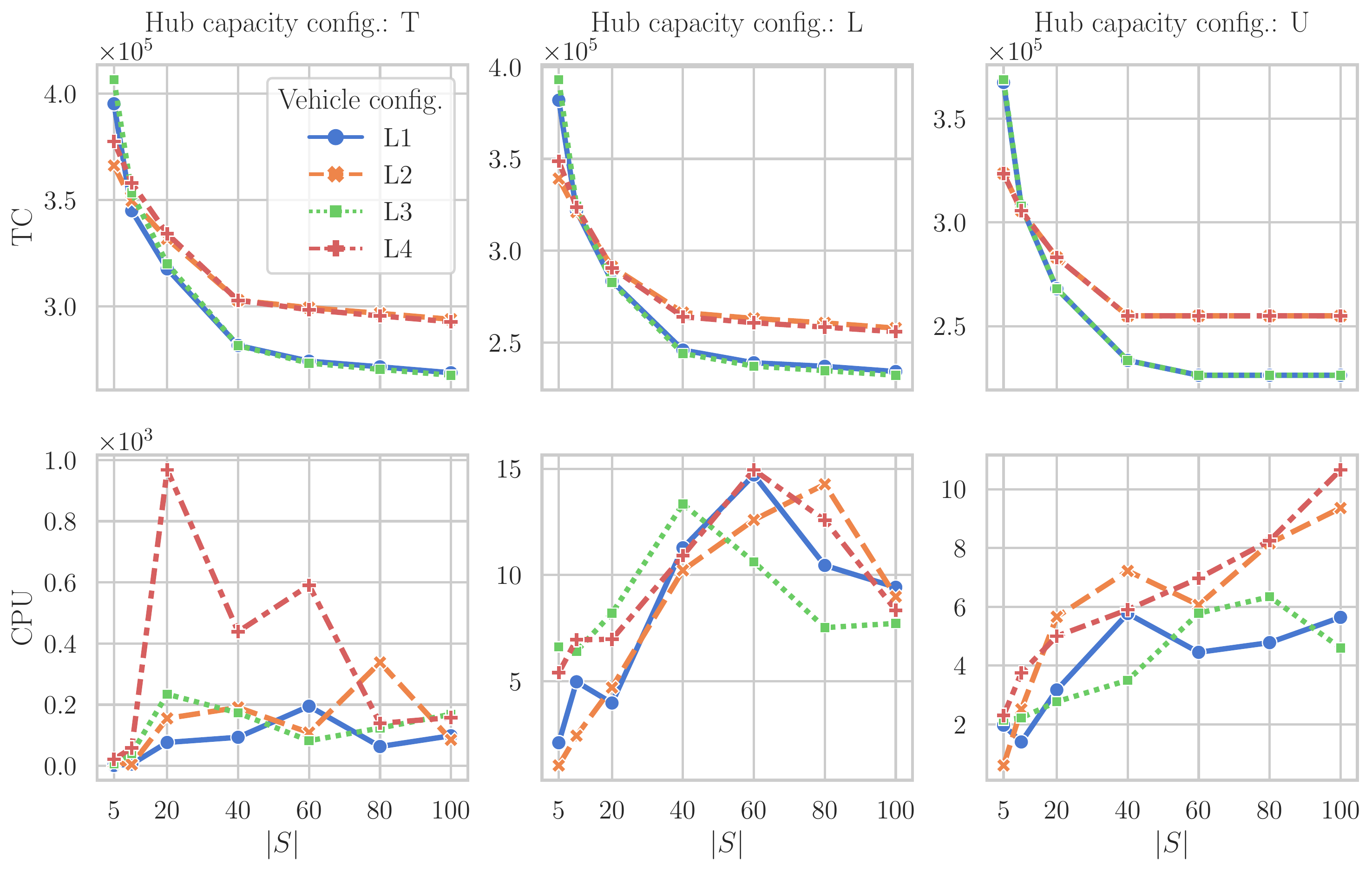}}
    {Effect of the number of scenarios on TC and CPU ($\card{N}=25$).\label{fig:ss-nscen}}
    {}
\end{figure}

In \autoref{fig:ss-nscen}, we illustrate the effect of increasing the number of scenarios on the objective function value and computational time.
Since the computational time increases very quickly as the number of scenarios increases for larger problems, we only consider the 25-node instances under 5 to 100 scenarios.
\autoref{fig:ss-nscen} shows that the increase in the number of scenarios leads to a lower TC value. When $\card{S}$ increases beyond 60, no significant change in TC is observed.
We observe that the effect of size $S$ on the CPU time  depends on the instance characteristics. In uncapacitated instances, increasing $\card{S}$ increases the required computational time.
This might not be the case when hubs are capacitated, specially with tight capacity levels, after $\card{S}$ passes a threshold.
We may explain this situation as follows. The more scenarios we have, the more variables we need to handle in our problem, and the more feasibility cuts are generated. When the number of demand scenarios gets very large for a particular node or when the hubs have very limited capacities, the options that provide a feasible allocation to a capacitated hub become less. Therefore, the algorithm may start with better bounds and fix variables more efficiently, hence its faster termination.

Overall, the stochastic solution may provide a better estimation of the total cost when the number of scenarios increases. However, it is expected that the capability of solving the instances with a large set $S$ becomes prohibitively limited.
\section{Conclusion}
\label{sec:conclusion}

Incorporating decisions about the number and type of vehicles to use adds more complexity to the hub network design problem. 
Therefore, finding the optimal solution to large-scale problem instances remains an issue in this area.
This study presented an efficient solution algorithm which is the first to solve large-scale benchmark instances with up to 200 nodes.
Our solution method relies on Benders decomposition with feasibility subproblems where the extreme rays have been derived in a closed-form solution resulting in a multiple-cut generation approach. Our computational experiments showed the superiority of this approach over the conventional Benders decomposition algorithm. To cope with more realistic situations, we addressed the HNDPv under demand uncertainty and showed the flexibility of our solution methodology in handling the stochastic variant of the problem.

While the benefits of vehicle-based hub network design problems are highlighted in this paper, several extensions can be investigated in the future to address different decisions at the tactical/operational levels. As vehicles are utilized to perform pickup/deliveries from/to the demand nodes, one may employ vehicle routes, instead of direct shipments, at the access-level network to reduce the number of trips and save costs.
Furthermore, to deal with demand uncertainty, one possible research direction is to consider some criteria (e.g., in a risk-averse manner) that strike a balance between the transportation cost and the risk of not having enough resources to meet the demand.

\ACKNOWLEDGMENT{%
    The second author acknowledges the financial support of the Natural Sciences and Engineering Research Council of Canada under Discovery Grant RGPIN- 2020-05395.
}%

\bibliographystyle{informs2014} %
\bibliography{references.bib} %

\end{document}


\maketitle

    

\section{Introduction}
\label{app:tables}
This online companion provides more detailed numerical results of our computational experiments addressed in the original paper.
The computational results of solving deterministic and stochastic HNDPv test instances with complete inter-hub network structure are provided in \autoref{sec:det} and \autoref{sec:stoch}, respectively.
\autoref{sec:gen} presents the results of solving deterministic problem instances with general network structure.

We use the following notation to show the performance of the proposed algorithms and the characteristics of the solution throughout the document.
\begin{inlinelist}
    \item Inst: instance name corresponding to the number of nodes in the instance, capacity type, and the vehicle configuration;
    \item Cap: hub capacity configuration;
    \item TC: total cost, i.e., the objective function value;
    \item LC: total hub location cost;
    \item HC: total transportation cost at the hub level;
    \item DC: total cost of direct shipments at the access level
    \item HFlow: total inter-hub flows;
    \item \#Hubs: number of open hubs;
    \item \#Veh1: number of primary vehicles used at the hub level;
    \item \#Veh2: number of secondary vehicles used at the access level;
    \item \%VUtil: average primary vehicle utilization (in percent)
    \item \#BNodes: number of branch-and-bound nodes explored by the solver;
    \item \#Cuts: number of Benders feasibility cuts generated by the proposed BCP algorithm;
    \item \%Gap: relative optimality gap (in percent);
    \item CPU\textsuperscript{r}: time spent to solve the root node (in seconds);
    \item CPU\textsuperscript{c}: time spent to generate feasibility cuts (in seconds); and
    \item CPU: total computational time (in seconds).
\end{inlinelist}

\section{Deterministic Instances}
\label{sec:det}

\autoref{tab:perf-comp-cap-small} and \autoref{tab:perf-comp-u-small} show the detailed performance measures of MIP and BC solution approaches for small and medium-size deterministic capacitated and uncapacitated instances, respectively.
The numbers regarding BC performance in solving large-size instances is provided in \autoref{tab:perf-comp-large}.
Values of different cost components and strategic/operational decisions of the solutions are listed in Tables \ref{tab:val-comp-t}--\ref{tab:val-comp-u}.

\begin{table}
    \centering
    \caption{BC performance in solving small and medium-size deterministic capacitated instances.\label{tab:perf-comp-cap-small}}
    {\footnotesize \begin{tabular}{@{}lrrrrrrrr@{}}
            \toprule
            \multicolumn{1}{c}{\multirow{2}{*}{Inst}}  & \multicolumn{3}{c}{MIP} & \multicolumn{5}{c}{BC}                                                          \\ \cmidrule(lr){2-4} \cmidrule(l){5-9}
            \multicolumn{1}{c}{}                       &
            \multicolumn{1}{c}{\#BNodes}               &
            \multicolumn{1}{c}{CPU\textsuperscript{r}} &
            \multicolumn{1}{c}{CPU (\%Gap)}            &
            \multicolumn{1}{c}{\#BNodes}               &
            \multicolumn{1}{c}{\#Cuts}                 &
            \multicolumn{1}{c}{CPU\textsuperscript{r}} &
            \multicolumn{1}{c}{CPU\textsuperscript{c}} &
            \multicolumn{1}{c@{}}{CPU (\%Gap)}                                                                                                                     \\ \midrule
            20T-L1                                     & 1.00e+0                 & 6.97                   & 7.27     & 1.69e+2 & 21    & 0.09   & 0.02  & 0.15     \\
            20T-L2                                     & 1.00e+0                 & 3.68                   & 3.69     & 2.13e+2 & 17    & 0.09   & 0.02  & 0.36     \\
            20T-L3                                     & 1.00e+0                 & 2.96                   & 2.97     & 3.25e+2 & 80    & 0.13   & 0.09  & 0.29     \\
            20T-L4                                     & 4.00e+1                 & 8.98                   & 9.18     & 7.59e+2 & 141   & 0.34   & 0.17  & 0.74     \\ \midrule
            20L-L1                                     & 1.00e+0                 & 9.87                   & 9.88     & 1.80e+2 & 41    & 0.21   & 0.05  & 0.40     \\
            20L-L2                                     & 1.00e+0                 & 6.52                   & 6.52     & 2.18e+2 & 47    & 0.15   & 0.06  & 0.27     \\
            20L-L3                                     & 1.00e+0                 & 2.68                   & 2.76     & 1.38e+2 & 62    & 0.28   & 0.07  & 0.37     \\
            20L-L4                                     & 1.00e+0                 & 0.00                   & 1.99     & 2.12e+2 & 78    & 0.26   & 0.09  & 0.44     \\ \midrule
            25T-L1                                     & 2.47e+2                 & 11.02                  & 16.53    & 8.90e+1 & 5     & 0.22   & 0.01  & 0.27     \\
            25T-L2                                     & 2.47e+2                 & 10.12                  & 14.84    & 1.13e+2 & 5     & 0.16   & 0.02  & 0.24     \\
            25T-L3                                     & 1.10e+2                 & 15.21                  & 17.79    & 6.67e+2 & 167   & 2.02   & 0.26  & 2.71     \\
            25T-L4                                     & 2.10e+2                 & 9.06                   & 12.88    & 6.07e+3 & 478   & 9.89   & 0.82  & 11.28    \\ \midrule
            25L-L1                                     & 7.20e+1                 & 18.74                  & 21.01    & 2.02e+2 & 62    & 1.25   & 0.11  & 1.50     \\
            25L-L2                                     & 4.10e+1                 & 24.00                  & 24.55    & 8.53e+2 & 131   & 0.81   & 0.27  & 1.36     \\
            25L-L3                                     & 1.00e+0                 & 10.82                  & 10.91    & 6.93e+2 & 113   & 1.21   & 0.17  & 1.41     \\
            25L-L4                                     & 1.00e+0                 & 9.99                   & 10.18    & 4.69e+2 & 114   & 0.80   & 0.17  & 1.12     \\ \midrule
            40T-L1                                     & 1.12e+2                 & 230.14                 & 266.52   & 5.00e+1 & 0     & 0.40   & 0.01  & 0.49     \\
            40T-L2                                     & 1.51e+3                 & 146.10                 & 219.05   & 1.13e+2 & 2     & 0.37   & 0.03  & 0.53     \\
            40T-L3                                     & 2.68e+3                 & 210.92                 & 270.53   & 3.79e+4 & 934   & 16.45  & 4.15  & 44.16    \\
            40T-L4                                     & 7.20e+3                 & 220.47                 & 304.10   & 3.70e+5 & 2,923 & 8.47   & 13.15 & 474.52   \\ \midrule
            40L-L1                                     & 1.00e+0                 & 98.26                  & 98.27    & 2.13e+2 & 28    & 0.88   & 0.14  & 1.08     \\
            40L-L2                                     & 1.00e+0                 & 44.45                  & 44.47    & 5.28e+2 & 77    & 2.61   & 0.35  & 3.12     \\
            40L-L3                                     & 1.00e+0                 & 87.67                  & 87.70    & 5.25e+2 & 47    & 0.89   & 0.23  & 1.68     \\
            40L-L4                                     & 1.00e+0                 & 42.06                  & 42.07    & 2.16e+2 & 28    & 1.49   & 0.09  & 1.70     \\ \midrule
            50T-L1                                     & 1.09e+4                 & 1,168.49               & 9,765.28 & 1.06e+2 & 18    & 1.91   & 0.13  & 2.50     \\
            50T-L2                                     & 4.92e+3                 & 552.46                 & 6,691.78 & 7.50e+1 & 3     & 1.09   & 0.03  & 1.41     \\
            50T-L3                                     & 2.45e+4                 & 1,425.87               & 5,218.47 & 8.67e+5 & 5,688 & 2.95   & 40.42 & 1,327.72 \\
            50T-L4                                     & 6.32e+3                 & 1,099.91               & 4,488.73 & 7.54e+3 & 185   & 9.64   & 1.19  & 11.68    \\ \midrule
            50L-L1                                     & 1.00e+0                 & 734.07                 & 734.93   & 1.64e+2 & 35    & 1.47   & 0.26  & 1.72     \\
            50L-L2                                     & 2.90e+1                 & 1,405.72               & 1,459.41 & 2.70e+1 & 2     & 0.69   & 0.03  & 0.80     \\
            50L-L3                                     & 1.14e+3                 & 847.77                 & 1,286.41 & 7.76e+3 & 243   & 3.95   & 1.66  & 7.49     \\
            50L-L4                                     & 2.52e+2                 & 1,209.23               & 1,402.19 & 5.00e+3 & 259   & 13.50  & 2.76  & 17.97    \\ \midrule
            75T-L1                                     & 1.41e+2                 & 5,186.28               & 5,808.71 & 1.00e+0 & 0     & 2.15   & 0.02  & 2.15     \\
            75T-L2                                     & 1.08e+2                 & 7,183.66               & 8,190.83 & 2.50e+1 & 0     & 2.15   & 0.03  & 2.34     \\
            75T-L3                                     & 2.07e+3                 & 8,048.18               & (3.19)   & 2.40e+6 & 7,210 & 93.52  & 95.54 & (2.71)   \\
            75T-L4                                     & 2.80e+1                 & 7,169.27               & (3.84)   & 6.51e+6 & 4,624 & 121.09 & 63.82 & (1.96)   \\ \midrule
            75L-L1                                     & 2.00e+1                 & 10,842.34              & (6.05)   & 5.50e+1 & 4     & 1.88   & 0.11  & 2.21     \\
            75L-L2                                     & 1.00e+0                 & 11,993.85              & (8.18)   & 4.50e+1 & 0     & 1.97   & 0.03  & 2.19     \\
            75L-L3                                     & 3.00e+0                 & 11,652.09              & (3.75)   & 1.18e+4 & 427   & 69.63  & 5.91  & 81.24    \\
            75L-L4                                     & 2.60e+1                 & 10,265.19              & (4.55)   & 1.02e+4 & 476   & 6.30   & 9.81  & 39.47    \\ \bottomrule
        \end{tabular}%
    }
\end{table}
\begin{table}
    \centering
    \caption{BC performance in solving small and medium-size deterministic uncapacitated instances.\label{tab:perf-comp-u-small}}
    {\footnotesize \begin{tabular}{@{}lrrrrrrrr@{}}
            \toprule
            \multicolumn{1}{c}{\multirow{2}{*}{Inst}}  & \multicolumn{3}{c}{MIP} & \multicolumn{5}{c}{BC}                                                   \\ \cmidrule(lr){2-4} \cmidrule(l){5-9}
            \multicolumn{1}{c}{}                       &
            \multicolumn{1}{c}{\#BNodes}               &
            \multicolumn{1}{c}{CPU\textsuperscript{r}} &
            \multicolumn{1}{c}{CPU (\%Gap)}            &
            \multicolumn{1}{c}{\#BNodes}               &
            \multicolumn{1}{c}{\#Cuts}                 &
            \multicolumn{1}{c}{CPU\textsuperscript{r}} &
            \multicolumn{1}{c}{CPU\textsuperscript{c}} &
            \multicolumn{1}{c@{}}{CPU (\%Gap)}                                                                                                              \\ \midrule
            20U-L1                                     & 1.00e+0                 & 8.11                   & 8.28     & 8.80e+1 & 22  & 0.15  & 0.02 & 0.20  \\
            20U-L2                                     & 1.00e+0                 & 8.02                   & 8.11     & 6.60e+1 & 20  & 0.27  & 0.02 & 0.32  \\
            20U-L3                                     & 1.00e+0                 & 3.02                   & 3.41     & 1.44e+2 & 54  & 0.16  & 0.06 & 0.24  \\
            20U-L4                                     & 1.00e+0                 & 0.00                   & 1.99     & 4.95e+2 & 108 & 0.13  & 0.12 & 0.33  \\ \midrule
            25U-L1                                     & 7.10e+1                 & 16.28                  & 18.41    & 1.69e+2 & 38  & 0.58  & 0.06 & 0.69  \\
            25U-L2                                     & 1.00e+0                 & 9.38                   & 9.40     & 2.97e+2 & 49  & 0.54  & 0.09 & 0.68  \\
            25U-L3                                     & 1.00e+0                 & 0.01                   & 4.72     & 8.60e+1 & 26  & 0.16  & 0.05 & 0.31  \\
            25U-L4                                     & 1.00e+0                 & 0.01                   & 2.21     & 8.70e+1 & 70  & 1.23  & 0.11 & 1.31  \\ \midrule
            40U-L1                                     & 1.00e+0                 & 110.74                 & 110.98   & 5.28e+2 & 100 & 2.18  & 0.46 & 2.70  \\
            40U-L2                                     & 1.00e+0                 & 41.52                  & 41.58    & 9.70e+1 & 35  & 1.38  & 0.14 & 1.65  \\
            40U-L3                                     & 1.00e+0                 & 71.23                  & 71.28    & 5.68e+2 & 72  & 0.91  & 0.35 & 1.39  \\
            40U-L4                                     & 1.00e+0                 & 44.66                  & 44.68    & 2.84e+3 & 171 & 2.32  & 0.80 & 3.47  \\ \midrule
            50U-L1                                     & 9.00e+1                 & 879.42                 & 915.23   & 1.06e+3 & 106 & 1.83  & 0.88 & 2.82  \\
            50U-L2                                     & 2.71e+2                 & 705.03                 & 988.54   & 5.09e+2 & 53  & 1.02  & 0.46 & 1.66  \\
            50U-L3                                     & 2.50e+1                 & 762.35                 & 766.35   & 9.39e+2 & 129 & 7.38  & 0.70 & 8.26  \\
            50U-L4                                     & 1.69e+2                 & 775.46                 & 1,238.07 & 1.19e+3 & 201 & 3.09  & 1.43 & 4.74  \\ \midrule
            75U-L1                                     & 3.09e+3                 & 5,184.36               & 9,324.58 & 5.30e+1 & 8   & 1.33  & 0.12 & 1.54  \\
            75U-L2                                     & 8.25e+2                 & 6,583.23               & (4.31)   & 5.20e+1 & 3   & 2.07  & 0.06 & 2.23  \\
            75U-L3                                     & 2.80e+1                 & 9,493.26               & (2.88)   & 1.09e+4 & 344 & 40.12 & 4.97 & 50.77 \\
            75U-L4                                     & 1.15e+3                 & 4,093.06               & 8,895.40 & 6.94e+3 & 536 & 5.57  & 7.62 & 40.26 \\ \bottomrule
        \end{tabular}%
    }
\end{table}
\begin{table}
    \centering
    \caption{BC performance in solving large deterministic instances.\label{tab:perf-comp-large}}
    {\footnotesize \begin{tabular}{@{}lrrrrr@{}}
            \toprule
            \multicolumn{1}{c}{Inst}                   &
            \multicolumn{1}{c}{\#BNodes}               &
            \multicolumn{1}{c}{\#Cuts}                 &
            \multicolumn{1}{c}{CPU\textsuperscript{r}} &
            \multicolumn{1}{c}{CPU\textsuperscript{c}} &
            \multicolumn{1}{c@{}}{CPU (\%Gap)}                                                        \\ \midrule
            100T-L1                                    & 2.43e+2 & 0     & 5.97   & 0.05   & 8.31     \\
            100T-L2                                    & 2.15e+2 & 0     & 3.66   & 0.09   & 5.57     \\
            100T-L3                                    & 1.09e+7 & 2,240 & 5.13   & 62.14  & (0.60)   \\
            100T-L4                                    & 3.01e+5 & 344   & 4.62   & 10.93  & 143.39   \\ \midrule
            100L-L1                                    & 3.20e+1 & 0     & 2.97   & 0.04   & 3.14     \\
            100L-L2                                    & 2.10e+1 & 36    & 6.14   & 0.81   & 6.26     \\
            100L-L3                                    & 3.42e+3 & 130   & 28.54  & 2.72   & 29.75    \\
            100L-L4                                    & 7.01e+3 & 39    & 7.69   & 0.97   & 11.95    \\ \midrule
            100U-L1                                    & 2.30e+1 & 4     & 2.25   & 0.15   & 2.35     \\
            100U-L2                                    & 1.70e+1 & 0     & 2.76   & 0.04   & 2.86     \\
            100U-L3                                    & 1.48e+2 & 282   & 75.83  & 6.20   & 76.16    \\
            100U-L4                                    & 4.67e+3 & 95    & 26.55  & 2.25   & 28.71    \\ \midrule
            150T-L1                                    & 2.73e+2 & 0     & 10.98  & 0.28   & 24.81    \\
            150T-L2                                    & 2.48e+2 & 0     & 9.27   & 0.32   & 23.49    \\
            150T-L3                                    & 2.52e+4 & 791   & 386.20 & 58.59  & 520.12   \\
            150T-L4                                    & 7.28e+4 & 1,656 & 11.26  & 116.60 & 399.36   \\ \midrule
            150L-L1                                    & 5.30e+1 & 0     & 8.57   & 0.11   & 9.40     \\
            150L-L2                                    & 8.50e+1 & 0     & 11.67  & 0.11   & 14.61    \\
            150L-L3                                    & 9.08e+4 & 1,001 & 514.39 & 56.53  & 722.13   \\
            150L-L4                                    & 1.16e+5 & 1,087 & 198.24 & 66.28  & 615.61   \\ \midrule
            150U-L1                                    & 8.30e+1 & 0     & 5.17   & 0.09   & 6.47     \\
            150U-L2                                    & 1.04e+2 & 0     & 9.65   & 0.13   & 12.79    \\
            150U-L3                                    & 8.97e+4 & 868   & 210.26 & 52.09  & 619.59   \\
            150U-L4                                    & 1.30e+5 & 932   & 31.35  & 60.20  & 421.49   \\ \midrule
            200T-L1                                    & 8.30e+2 & 0     & 39.09  & 0.55   & 109.19   \\
            200T-L2                                    & 1.19e+3 & 0     & 37.31  & 0.63   & 106.40   \\
            200T-L3                                    & 6.89e+2 & 6     & 36.59  & 0.85   & 80.79    \\
            200T-L4                                    & 1.04e+3 & 1     & 31.83  & 0.75   & 104.07   \\ \midrule
            200L-L1                                    & 1.08e+3 & 0     & 28.88  & 0.56   & 83.46    \\
            200L-L2                                    & 1.46e+3 & 0     & 39.48  & 0.52   & 106.87   \\
            200L-L3                                    & 1.87e+6 & 4,400 & 73.97  & 454.45 & (0.24)   \\
            200L-L4                                    & 1.62e+6 & 1,242 & 90.06  & 134.70 & 7,395.43 \\ \midrule
            200U-L1                                    & 1.10e+3 & 0     & 29.05  & 0.64   & 56.82    \\
            200U-L2                                    & 1.48e+3 & 0     & 28.74  & 0.66   & 70.47    \\
            200U-L3                                    & 5.18e+5 & 9,900 & 451.13 & 935.23 & 7,689.79 \\
            200U-L4                                    & 1.14e+6 & 8,525 & 31.50  & 816.18 & (0.09)   \\ \bottomrule
        \end{tabular}%
    }
\end{table}

\begin{table}
    \centering
    \caption{Solution details of the instances with tight capacity configuration.\label{tab:val-comp-t}}
    {\footnotesize \begin{tabular}{@{}lrrrrrrrrr@{}}
            \toprule
            \multicolumn{1}{c}{Inst}   &
            \multicolumn{1}{c}{TC}     &
            \multicolumn{1}{c}{LC}     &
            \multicolumn{1}{c}{HC}     &
            \multicolumn{1}{c}{DC}     &
            \multicolumn{1}{c}{HFlow}  &
            \multicolumn{1}{c}{\#Hubs} &
            \multicolumn{1}{c}{\#Veh1} &
            \multicolumn{1}{c}{\#Veh2} &
            \multicolumn{1}{c@{}}{\%VUtil}                                                                                       \\ \midrule
            20T-L1                     & 393,933.67   & 75,763.01  & 28,123.55  & 290,047.10   & 1,760.61 & 2 & 4  & 48  & 73.36 \\
            20T-L2                     & 358,634.95   & 75,763.01  & 28,123.55  & 254,748.38   & 1,760.61 & 2 & 4  & 36  & 73.36 \\
            20T-L3                     & 400,964.56   & 75,763.01  & 35,154.44  & 290,047.10   & 1,760.61 & 2 & 6  & 48  & 91.70 \\
            20T-L4                     & 365,665.84   & 75,763.01  & 35,154.44  & 254,748.38   & 1,760.61 & 2 & 6  & 36  & 91.70 \\ \midrule
            25T-L1                     & 466,619.56   & 135,626.35 & 39,049.37  & 291,943.84   & 2,172.08 & 3 & 6  & 50  & 60.34 \\
            25T-L2                     & 421,901.84   & 135,626.35 & 39,049.37  & 247,226.13   & 2,172.08 & 3 & 6  & 37  & 60.34 \\
            25T-L3                     & 473,230.43   & 135,626.35 & 45,053.60  & 292,550.48   & 2,161.53 & 3 & 8  & 50  & 84.43 \\
            25T-L4                     & 428,606.04   & 135,626.35 & 45,053.60  & 247,926.09   & 2,161.53 & 3 & 8  & 37  & 84.43 \\ \midrule
            40T-L1                     & 512,392.47   & 131,166.19 & 52,074.69  & 329,151.59   & 2,231.72 & 3 & 6  & 58  & 61.99 \\
            40T-L2                     & 487,239.06   & 108,415.85 & 56,441.18  & 322,382.03   & 2,333.33 & 3 & 6  & 45  & 64.81 \\
            40T-L3                     & 524,122.65   & 131,166.19 & 59,967.56  & 332,988.90   & 2,244.60 & 3 & 8  & 58  & 87.68 \\
            40T-L4                     & 498,033.20   & 131,166.19 & 59,967.56  & 306,899.45   & 2,259.29 & 3 & 8  & 45  & 88.25 \\ \midrule
            50T-L1                     & 606,952.39   & 190,994.69 & 37,034.70  & 378,923.01   & 2,351.07 & 3 & 6  & 62  & 65.31 \\
            50T-L2                     & 625,720.88   & 239,836.08 & 60,576.55  & 325,308.25   & 2,163.50 & 3 & 6  & 56  & 60.10 \\
            50T-L3                     & 605,982.07   & 189,332.90 & 38,561.64  & 378,087.53   & 2,306.82 & 3 & 8  & 62  & 90.11 \\
            50T-L4                     & 627,201.17   & 239,836.08 & 61,915.06  & 325,450.03   & 2,162.33 & 3 & 8  & 56  & 84.47 \\ \midrule
            75T-L1                     & 638,638.40   & 123,463.64 & 57,693.27  & 457,481.49   & 2,362.14 & 3 & 6  & 86  & 65.61 \\
            75T-L2                     & 682,970.06   & 123,463.64 & 57,693.27  & 501,813.15   & 2,362.14 & 3 & 6  & 80  & 65.61 \\
            75T-L3                     & 660,172.05   & 123,463.64 & 74,167.93  & 462,540.48   & 2,270.51 & 3 & 10 & 86  & 70.95 \\
            75T-L4                     & 703,216.95   & 123,463.64 & 74,167.93  & 505,585.38   & 2,270.51 & 3 & 10 & 80  & 70.95 \\ \midrule
            100T-L1                    & 959,502.13   & 290,007.91 & 149,656.59 & 519,837.62   & 2,508.70 & 4 & 12 & 109 & 34.84 \\
            100T-L2                    & 996,488.07   & 290,007.91 & 149,656.59 & 556,823.56   & 2,518.63 & 4 & 12 & 101 & 34.98 \\
            100T-L3                    & 948,021.18   & 290,007.91 & 138,175.66 & 519,837.62   & 2,508.70 & 4 & 14 & 109 & 56.00 \\
            100T-L4                    & 982,736.28   & 290,007.91 & 131,444.74 & 561,283.63   & 2,520.77 & 4 & 13 & 101 & 60.60 \\ \midrule
            150T-L1                    & 989,945.89   & 119,093.32 & 119,861.86 & 750,990.70   & 2,470.96 & 4 & 12 & 158 & 34.32 \\
            150T-L2                    & 1,061,402.32 & 119,093.32 & 119,861.86 & 822,447.14   & 2,470.96 & 4 & 12 & 151 & 34.32 \\
            150T-L3                    & 972,945.76   & 119,093.32 & 99,884.88  & 753,967.56   & 2,436.72 & 4 & 12 & 158 & 63.46 \\
            150T-L4                    & 1,046,745.20 & 119,093.32 & 99,884.88  & 827,767.00   & 2,473.54 & 4 & 12 & 151 & 64.42 \\ \midrule
            200T-L1                    & 1,237,882.56 & 120,274.44 & 215,456.20 & 902,151.92   & 2,748.19 & 5 & 20 & 206 & 22.90 \\
            200T-L2                    & 1,352,076.99 & 120,274.44 & 215,456.20 & 1,016,346.35 & 2,748.19 & 5 & 20 & 199 & 22.90 \\
            200T-L3                    & 1,201,973.19 & 120,274.44 & 179,546.83 & 902,151.92   & 2,748.19 & 5 & 20 & 206 & 42.94 \\
            200T-L4                    & 1,316,167.62 & 120,274.44 & 179,546.83 & 1,016,346.35 & 2,748.19 & 5 & 20 & 199 & 42.94 \\ \bottomrule
        \end{tabular}%
    }
\end{table}
\begin{table}
    \centering
    \caption{Solution details of the instances with loose capacity configuration.\label{tab:val-comp-l}}
    {\footnotesize \begin{tabular}{@{}lrrrrrrrrr@{}}
            \toprule
            \multicolumn{1}{c}{Inst}   &
            \multicolumn{1}{c}{TC}     &
            \multicolumn{1}{c}{LC}     &
            \multicolumn{1}{c}{HC}     &
            \multicolumn{1}{c}{DC}     &
            \multicolumn{1}{c}{HFlow}  &
            \multicolumn{1}{c}{\#Hubs} &
            \multicolumn{1}{c}{\#Veh1} &
            \multicolumn{1}{c}{\#Veh2} &
            \multicolumn{1}{c@{}}{\%VUtil}                                                                                       \\ \midrule
            20L-L1                     & 370,221.96   & 77,781.69  & 33,217.39  & 259,222.88   & 1,423.66 & 2 & 3  & 48  & 79.09 \\
            20L-L2                     & 334,821.03   & 77,781.69  & 33,217.39  & 223,821.95   & 1,423.66 & 2 & 3  & 36  & 79.09 \\
            20L-L3                     & 375,980.02   & 28,436.87  & 0.00       & 347,543.15   & 0.00     & 1 & 0  & 50  &       \\
            20L-L4                     & 336,616.68   & 28,436.87  & 0.00       & 308,179.81   & 0.00     & 1 & 0  & 38  &       \\ \midrule
            25L-L1                     & 434,368.05   & 128,271.77 & 32,736.50  & 273,359.79   & 1,464.97 & 2 & 3  & 51  & 81.39 \\
            25L-L2                     & 387,027.79   & 114,600.55 & 33,305.03  & 239,122.22   & 1,400.52 & 2 & 3  & 37  & 77.81 \\
            25L-L3                     & 441,156.96   & 114,600.55 & 46,256.98  & 280,299.42   & 1,572.29 & 2 & 5  & 52  & 98.27 \\
            25L-L4                     & 394,911.22   & 114,600.55 & 46,256.98  & 234,053.69   & 1,572.29 & 2 & 5  & 37  & 98.27 \\ \midrule
            40L-L1                     & 459,954.90   & 63,320.02  & 18,452.88  & 378,182.00   & 1,441.66 & 2 & 3  & 59  & 80.09 \\
            40L-L2                     & 424,777.22   & 63,320.02  & 18,452.88  & 343,004.33   & 1,453.92 & 2 & 3  & 46  & 80.77 \\
            40L-L3                     & 465,508.65   & 63,320.02  & 25,629.00  & 376,559.63   & 1,508.31 & 2 & 5  & 59  & 94.27 \\
            40L-L4                     & 430,719.72   & 63,320.02  & 25,629.00  & 341,770.71   & 1,508.31 & 2 & 5  & 46  & 94.27 \\ \midrule
            50L-L1                     & 501,301.26   & 93,859.16  & 22,081.68  & 385,360.42   & 1,087.24 & 2 & 2  & 63  & 90.60 \\
            50L-L2                     & 524,568.40   & 93,859.16  & 22,081.68  & 408,627.56   & 1,087.24 & 2 & 2  & 57  & 90.60 \\
            50L-L3                     & 516,022.38   & 93,859.16  & 36,802.80  & 385,360.42   & 1,087.24 & 2 & 4  & 63  & 84.94 \\
            50L-L4                     & 539,289.52   & 93,859.16  & 36,802.80  & 408,627.56   & 1,087.24 & 2 & 4  & 57  & 84.94 \\ \midrule
            75L-L1                     & 608,275.69   & 124,218.94 & 58,375.88  & 425,680.87   & 2,072.25 & 3 & 6  & 86  & 57.56 \\
            75L-L2                     & 655,180.78   & 124,218.94 & 58,375.88  & 472,585.95   & 2,072.25 & 3 & 6  & 80  & 57.56 \\
            75L-L3                     & 612,805.00   & 111,886.08 & 51,520.45  & 449,398.46   & 2,279.96 & 3 & 8  & 86  & 89.06 \\
            75L-L4                     & 660,107.60   & 124,218.94 & 62,412.50  & 473,476.15   & 2,018.48 & 3 & 8  & 80  & 78.85 \\ \midrule
            100L-L1                    & 818,129.60   & 135,801.12 & 80,709.13  & 601,619.35   & 1,944.14 & 3 & 6  & 110 & 54.00 \\
            100L-L2                    & 859,558.42   & 135,801.12 & 80,709.13  & 643,048.17   & 1,926.45 & 3 & 6  & 102 & 53.51 \\
            100L-L3                    & 821,113.09   & 135,801.12 & 83,609.62  & 601,702.35   & 1,915.22 & 3 & 8  & 110 & 74.81 \\
            100L-L4                    & 862,458.91   & 135,801.12 & 83,609.62  & 643,048.17   & 1,926.45 & 3 & 8  & 102 & 75.25 \\ \midrule
            150L-L1                    & 900,848.18   & 105,159.67 & 103,668.22 & 692,020.29   & 2,484.61 & 4 & 12 & 158 & 34.51 \\
            150L-L2                    & 992,952.58   & 105,159.67 & 103,668.22 & 784,124.69   & 2,484.61 & 4 & 12 & 151 & 34.51 \\
            150L-L3                    & 892,931.33   & 105,159.67 & 92,527.64  & 695,244.01   & 2,450.68 & 4 & 13 & 158 & 58.91 \\
            150L-L4                    & 985,531.68   & 105,159.67 & 92,527.64  & 787,844.37   & 2,456.14 & 4 & 13 & 151 & 59.04 \\ \midrule
            200L-L1                    & 1,207,345.18 & 104,030.28 & 125,968.59 & 977,346.31   & 2,261.76 & 4 & 12 & 207 & 31.41 \\
            200L-L2                    & 1,330,584.32 & 106,407.66 & 133,326.50 & 1,090,850.16 & 2,357.71 & 4 & 12 & 200 & 32.75 \\
            200L-L3                    & 1,194,673.29 & 104,030.28 & 112,098.72 & 978,544.29   & 2,318.17 & 4 & 13 & 207 & 55.73 \\
            200L-L4                    & 1,316,087.39 & 126,894.38 & 174,748.34 & 1,014,444.67 & 2,541.35 & 5 & 20 & 199 & 39.71 \\ \bottomrule
        \end{tabular}%
    }
\end{table}
\begin{table}
    \centering
    \caption{Solution details of the ucapacitated instances.\label{tab:val-comp-u}}
    {\footnotesize \begin{tabular}{@{}lrrrrrrrrr@{}}
            \toprule
            \multicolumn{1}{c}{Inst}   &
            \multicolumn{1}{c}{TC}     &
            \multicolumn{1}{c}{LC}     &
            \multicolumn{1}{c}{HC}     &
            \multicolumn{1}{c}{DC}     &
            \multicolumn{1}{c}{HFlow}  &
            \multicolumn{1}{c}{\#Hubs} &
            \multicolumn{1}{c}{\#Veh1} &
            \multicolumn{1}{c}{\#Veh2} &
            \multicolumn{1}{c@{}}{\%VUtil}                                                                                       \\ \midrule
            20U-L1                     & 370,221.96   & 77,781.69  & 33,217.39  & 259,222.88   & 1,423.66 & 2 & 3  & 48  & 79.09 \\
            20U-L2                     & 334,821.03   & 77,781.69  & 33,217.39  & 223,821.95   & 1,423.66 & 2 & 3  & 36  & 79.09 \\
            20U-L3                     & 375,980.02   & 28,436.87  & 0.00       & 347,543.15   & 0.00     & 1 & 0  & 50  &       \\
            20U-L4                     & 336,616.68   & 28,436.87  & 0.00       & 308,179.81   & 0.00     & 1 & 0  & 38  &       \\ \midrule
            25U-L1                     & 423,528.15   & 98,853.64  & 35,333.69  & 289,340.83   & 1,464.97 & 2 & 3  & 51  & 81.39 \\
            25U-L2                     & 369,691.77   & 40,332.49  & 0.00       & 329,359.28   & 0.00     & 1 & 0  & 39  &       \\
            25U-L3                     & 425,922.44   & 40,332.49  & 0.00       & 385,589.95   & 0.00     & 1 & 0  & 53  &       \\
            25U-L4                     & 369,691.77   & 40,332.49  & 0.00       & 329,359.28   & 0.00     & 1 & 0  & 39  &       \\ \midrule
            40U-L1                     & 459,954.90   & 63,320.02  & 18,452.88  & 378,182.00   & 1,441.66 & 2 & 3  & 59  & 80.09 \\
            40U-L2                     & 424,777.22   & 63,320.02  & 18,452.88  & 343,004.33   & 1,453.92 & 2 & 3  & 46  & 80.77 \\
            40U-L3                     & 465,508.65   & 63,320.02  & 25,629.00  & 376,559.63   & 1,508.31 & 2 & 5  & 59  & 94.27 \\
            40U-L4                     & 430,719.72   & 63,320.02  & 25,629.00  & 341,770.71   & 1,508.31 & 2 & 5  & 46  & 94.27 \\ \midrule
            50U-L1                     & 483,848.26   & 96,315.70  & 35,198.38  & 352,334.19   & 1,387.50 & 2 & 3  & 63  & 77.08 \\
            50U-L2                     & 520,570.44   & 96,315.70  & 35,198.38  & 389,056.37   & 1,387.50 & 2 & 3  & 57  & 77.08 \\
            50U-L3                     & 495,933.46   & 96,315.70  & 48,886.64  & 350,731.13   & 1,495.32 & 2 & 5  & 63  & 93.46 \\
            50U-L4                     & 532,409.02   & 96,315.70  & 48,886.64  & 387,206.69   & 1,495.32 & 2 & 5  & 57  & 93.46 \\ \midrule
            75U-L1                     & 608,275.69   & 124,218.94 & 58,375.88  & 425,680.87   & 2,072.25 & 3 & 6  & 86  & 57.56 \\
            75U-L2                     & 655,180.78   & 124,218.94 & 58,375.88  & 472,585.95   & 2,072.25 & 3 & 6  & 80  & 57.56 \\
            75U-L3                     & 612,805.00   & 111,886.08 & 51,520.45  & 449,398.46   & 2,279.96 & 3 & 8  & 86  & 89.06 \\
            75U-L4                     & 660,107.60   & 124,218.94 & 62,412.50  & 473,476.15   & 2,018.48 & 3 & 8  & 80  & 78.85 \\ \midrule
            100U-L1                    & 810,993.33   & 135,801.12 & 80,709.13  & 594,483.08   & 1,750.28 & 3 & 6  & 110 & 48.62 \\
            100U-L2                    & 854,636.49   & 135,801.12 & 80,709.13  & 638,126.23   & 1,750.28 & 3 & 6  & 102 & 48.62 \\
            100U-L3                    & 813,893.82   & 135,801.12 & 83,609.62  & 594,483.08   & 1,750.28 & 3 & 8  & 110 & 68.37 \\
            100U-L4                    & 857,536.97   & 135,801.12 & 83,609.62  & 638,126.23   & 1,750.28 & 3 & 8  & 102 & 68.37 \\ \midrule
            150U-L1                    & 900,848.18   & 105,159.67 & 103,668.22 & 692,020.29   & 2,481.92 & 4 & 12 & 158 & 34.47 \\
            150U-L2                    & 992,952.58   & 105,159.67 & 103,668.22 & 784,124.69   & 2,479.64 & 4 & 12 & 151 & 34.44 \\
            150U-L3                    & 892,931.33   & 105,159.67 & 92,527.64  & 695,244.01   & 2,450.68 & 4 & 13 & 158 & 58.91 \\
            150U-L4                    & 985,531.68   & 105,159.67 & 92,527.64  & 787,844.37   & 2,450.68 & 4 & 13 & 151 & 58.91 \\ \midrule
            200U-L1                    & 1,202,840.91 & 104,030.28 & 125,968.59 & 972,842.04   & 2,329.75 & 4 & 12 & 207 & 32.36 \\
            200U-L2                    & 1,325,387.08 & 106,407.66 & 133,326.50 & 1,085,652.92 & 2,421.86 & 4 & 12 & 200 & 33.64 \\
            200U-L3                    & 1,186,732.32 & 104,030.28 & 104,973.83 & 977,728.21   & 2,274.67 & 4 & 12 & 207 & 59.24 \\
            200U-L4                    & 1,310,489.14 & 142,252.10 & 171,591.30 & 996,645.74   & 2,637.12 & 5 & 20 & 199 & 41.20 \\ \bottomrule
        \end{tabular}%
    }
\end{table}

\section{Instances with Stochastic Demands}
\label{sec:stoch}

\autoref{tab:perf-comp-cap-ss} and \autoref{tab:perf-comp-u-ss} show the detailed performance measures of BC for capacitated and uncapacitated stochastic instances under five scenarios, respectively.
Values of different cost components and strategic/operational decisions of the stochastic solutions are listed in Tables \ref{tab:val-comp-cap-ss} and \ref{tab:val-comp-u-ss}.

\begin{table}
    \centering
    \caption{BC performance in solving stochastic capacitated instances.\label{tab:perf-comp-cap-ss}}
    {\footnotesize \begin{tabular}{@{}lrrrrr@{}}
            \toprule
            \multicolumn{1}{c}{Inst}                   &
            \multicolumn{1}{c}{\#BNodes}               &
            \multicolumn{1}{c}{\#Cuts}                 &
            \multicolumn{1}{c}{CPU\textsuperscript{r}} &
            \multicolumn{1}{c}{CPU\textsuperscript{c}} &
            \multicolumn{1}{c@{}}{CPU (\%Gap)}                                                           \\ \midrule
            25T-L1                                     & 4.52e+2 & 295    & 1.07     & 0.59   & 1.48     \\
            25T-L2                                     & 1.05e+3 & 467    & 9.07     & 0.91   & 9.93     \\
            25T-L3                                     & 2.39e+3 & 777    & 2.61     & 1.50   & 7.44     \\
            25T-L4                                     & 1.17e+4 & 1,168  & 11.95    & 2.38   & 21.73    \\ \midrule
            25L-L1                                     & 5.38e+2 & 285    & 1.41     & 0.75   & 2.10     \\
            25L-L2                                     & 1.78e+2 & 88     & 0.79     & 0.38   & 1.03     \\
            25L-L3                                     & 5.08e+3 & 637    & 4.81     & 1.58   & 6.61     \\
            25L-L4                                     & 1.45e+3 & 346    & 4.09     & 0.87   & 5.40     \\ \midrule
            40T-L1                                     & 6.38e+4 & 3,847  & 65.76    & 11.90  & 191.50   \\
            40T-L2                                     & 1.69e+7 & 2,593  & 27.49    & 8.59   & (4.75)   \\
            40T-L3                                     & 3.50e+6 & 15,208 & 62.79    & 53.55  & (3.44)   \\
            40T-L4                                     & 1.86e+6 & 12,733 & 28.27    & 47.05  & (6.55)   \\ \midrule
            40L-L1                                     & 6.84e+3 & 1,566  & 30.88    & 6.04   & 60.40    \\
            40L-L2                                     & 1.03e+4 & 765    & 25.19    & 3.00   & 30.72    \\
            40L-L3                                     & 2.05e+4 & 3,105  & 66.59    & 11.21  & 122.32   \\
            40L-L4                                     & 8.48e+3 & 2,053  & 53.09    & 7.43   & 77.01    \\ \midrule
            50T-L1                                     & 6.60e+6 & 6,407  & 24.39    & 31.27  & (4.97)   \\
            50T-L2                                     & 5.89e+6 & 6,999  & 16.05    & 35.62  & (6.23)   \\
            50T-L3                                     & 2.23e+6 & 12,753 & 24.75    & 69.52  & (2.67)   \\
            50T-L4                                     & 2.85e+6 & 15,519 & 65.72    & 84.11  & (4.59)   \\ \midrule
            50L-L1                                     & 1.45e+4 & 1,176  & 13.23    & 7.48   & 43.67    \\
            50L-L2                                     & 2.96e+4 & 1,515  & 13.12    & 8.67   & 49.20    \\
            50L-L3                                     & 1.41e+6 & 2,495  & 64.75    & 16.50  & 1,818.67 \\
            50L-L4                                     & 6.57e+5 & 2,634  & 104.50   & 16.00  & 802.71   \\ \midrule
            75T-L1                                     & 5.06e+5 & 17,601 & 1,868.20 & 177.01 & (9.03)   \\
            75T-L2                                     & 4.80e+5 & 18,860 & 148.83   & 197.63 & (15.43)  \\
            75T-L3                                     & 1.87e+5 & 26,427 & 596.87   & 275.09 & (9.95)   \\
            75T-L4                                     & 2.80e+5 & 29,857 & 390.38   & 313.68 & (13.57)  \\ \midrule
            75L-L1                                     & 5.24e+6 & 6,282  & 78.15    & 65.45  & (2.88)   \\
            75L-L2                                     & 3.56e+6 & 8,305  & 74.42    & 91.43  & (2.75)   \\
            75L-L3                                     & 2.95e+5 & 20,655 & 180.81   & 230.93 & (2.89)   \\
            75L-L4                                     & 3.79e+4 & 25,867 & 189.82   & 278.77 & (8.57)   \\ \bottomrule
        \end{tabular}%
    }
\end{table}

\begin{table}
    \centering
    \caption{BC performance in solving stochastic ucapacitated instances.\label{tab:perf-comp-u-ss}}
    {\footnotesize \begin{tabular}{@{}lrrrrr@{}}
            \toprule
            \multicolumn{1}{c}{Inst}                   &
            \multicolumn{1}{c}{\#BNodes}               &
            \multicolumn{1}{c}{\#Cuts}                 &
            \multicolumn{1}{c}{CPU\textsuperscript{r}} &
            \multicolumn{1}{c}{CPU\textsuperscript{c}} &
            \multicolumn{1}{c@{}}{CPU (\%Gap)}                                                         \\ \midrule
            25U-L1                                     & 4.27e+2 & 210    & 1.23     & 0.57   & 1.97   \\
            25U-L2                                     & 2.30e+1 & 55     & 0.49     & 0.20   & 0.60   \\
            25U-L3                                     & 6.01e+2 & 300    & 1.65     & 0.65   & 2.14   \\
            25U-L4                                     & 2.40e+1 & 153    & 2.26     & 0.39   & 2.30   \\ \midrule
            40U-L1                                     & 1.69e+4 & 1,992  & 39.67    & 8.54   & 74.11  \\
            40U-L2                                     & 1.40e+3 & 723    & 22.68    & 2.76   & 24.29  \\
            40U-L3                                     & 4.69e+4 & 1,909  & 54.20    & 7.12   & 98.31  \\
            40U-L4                                     & 4.77e+3 & 2,259  & 128.25   & 8.04   & 136.29 \\ \midrule
            50U-L1                                     & 2.31e+5 & 1,367  & 93.37    & 9.12   & 264.06 \\
            50U-L2                                     & 4.99e+4 & 1,795  & 70.26    & 11.22  & 150.13 \\
            50U-L3                                     & 1.39e+5 & 3,107  & 229.64   & 17.05  & 412.01 \\
            50U-L4                                     & 5.67e+5 & 3,093  & 105.50   & 18.33  & 672.07 \\ \midrule
            75U-L1                                     & 4.41e+6 & 4,940  & 81.45    & 55.90  & (4.62) \\
            75U-L2                                     & 4.52e+6 & 6,172  & 571.78   & 66.03  & (4.27) \\
            75U-L3                                     & 1.86e+6 & 10,870 & 1,072.15 & 121.01 & (3.88) \\
            75U-L4                                     & 1.38e+6 & 29,514 & 103.32   & 301.37 & (3.89) \\ \bottomrule
        \end{tabular}%
    }
\end{table}

\begin{table}
    \centering
    \caption{Solution details of the stochastic capacitated instances.\label{tab:val-comp-cap-ss}}
    {\footnotesize \begin{tabular}{@{}lrrrrrrrrr@{}}
            \toprule
            \multicolumn{1}{c}{Inst}   &
            \multicolumn{1}{c}{TC}     &
            \multicolumn{1}{c}{LC}     &
            \multicolumn{1}{c}{HC}     &
            \multicolumn{1}{c}{DC}     &
            \multicolumn{1}{c}{HFlow}  &
            \multicolumn{1}{c}{\#Hubs} &
            \multicolumn{1}{c}{\#Veh1} &
            \multicolumn{1}{c}{\#Veh2} &
            \multicolumn{1}{c@{}}{\%VUtil}                                                                                     \\ \midrule
            25T-L1                     & 395,264.34 & 98,034.61  & 30,723.11  & 266,506.62 & 1,228.80 & 2 & 2.80  & 41 & 73.14 \\
            25T-L2                     & 366,191.96 & 98,034.61  & 30,723.11  & 237,434.25 & 1,228.80 & 2 & 2.80  & 31 & 73.14 \\
            25T-L3                     & 406,602.63 & 98,034.61  & 42,061.40  & 266,506.62 & 1,228.80 & 2 & 4.60  & 41 & 83.48 \\
            25T-L4                     & 377,530.25 & 98,034.61  & 42,061.40  & 237,434.25 & 1,228.80 & 2 & 4.60  & 31 & 83.48 \\ \midrule
            25L-L1                     & 381,937.67 & 114,600.55 & 31,084.69  & 236,252.43 & 1,254.40 & 2 & 2.80  & 41 & 74.67 \\
            25L-L2                     & 339,249.83 & 62,876.95  & 25,741.49  & 250,631.38 & 825.60   & 2 & 2.00  & 31 & 68.80 \\
            25L-L3                     & 393,147.18 & 68,115.09  & 14,642.09  & 310,390.00 & 1,008.00 & 2 & 4.20  & 42 & 75.00 \\
            25L-L4                     & 348,719.00 & 62,876.95  & 36,467.11  & 249,374.94 & 857.40   & 2 & 3.40  & 31 & 78.81 \\ \midrule
            40T-L1                     & 580,503.96 & 138,434.32 & 74,996.13  & 367,073.51 & 2,346.20 & 3 & 6.00  & 65 & 65.17 \\
            40T-L2                     & 542,719.38 & 138,434.32 & 74,996.13  & 329,288.93 & 2,346.20 & 3 & 6.00  & 49 & 65.17 \\
            40T-L3                     & 611,117.13 & 179,627.61 & 90,801.15  & 340,688.38 & 2,810.80 & 4 & 12.60 & 64 & 69.71 \\
            40T-L4                     & 568,371.24 & 179,627.61 & 91,531.99  & 297,211.66 & 2,817.40 & 4 & 12.80 & 47 & 68.78 \\ \midrule
            40L-L1                     & 501,505.40 & 63,320.02  & 23,373.65  & 414,811.73 & 1,780.20 & 2 & 3.80  & 67 & 78.08 \\
            40L-L2                     & 455,346.14 & 63,320.02  & 23,373.65  & 368,652.48 & 1,780.20 & 2 & 3.80  & 50 & 78.08 \\
            40L-L3                     & 509,126.60 & 63,320.02  & 30,754.80  & 415,051.79 & 1,738.60 & 2 & 6.00  & 67 & 90.55 \\
            40L-L4                     & 463,004.29 & 63,320.02  & 30,754.80  & 368,929.47 & 1,738.60 & 2 & 6.00  & 50 & 90.55 \\ \midrule
            50T-L1                     & 632,538.04 & 239,836.08 & 60,576.55  & 332,125.42 & 1,887.40 & 3 & 6.00  & 63 & 52.43 \\
            50T-L2                     & 638,805.43 & 239,836.08 & 60,576.55  & 338,392.80 & 1,887.40 & 3 & 6.00  & 54 & 52.43 \\
            50T-L3                     & 633,528.79 & 239,836.08 & 59,628.14  & 334,064.58 & 1,923.20 & 3 & 7.60  & 63 & 79.08 \\
            50T-L4                     & 645,546.49 & 214,563.92 & 83,351.22  & 347,631.36 & 2,335.00 & 4 & 12.40 & 53 & 58.85 \\ \midrule
            50L-L1                     & 519,399.98 & 93,859.16  & 26,498.01  & 399,042.81 & 1,078.20 & 2 & 2.40  & 64 & 74.88 \\
            50L-L2                     & 527,339.56 & 93,859.16  & 26,498.01  & 406,982.38 & 1,078.20 & 2 & 2.40  & 55 & 74.88 \\
            50L-L3                     & 532,303.70 & 93,859.16  & 36,802.80  & 401,641.74 & 1,012.00 & 2 & 4.00  & 64 & 79.06 \\
            50L-L4                     & 538,954.06 & 112,527.43 & 21,042.71  & 405,383.92 & 1,164.60 & 2 & 4.40  & 55 & 82.71 \\ \midrule
            75T-L1                     & 800,332.99 & 201,590.46 & 92,821.78  & 505,920.75 & 3,438.80 & 4 & 12.00 & 95 & 47.76 \\
            75T-L2                     & 863,739.63 & 196,870.68 & 198,986.31 & 467,882.64 & 3,759.00 & 5 & 20.00 & 81 & 31.32 \\
            75T-L3                     & 824,989.34 & 227,110.99 & 135,204.83 & 462,673.53 & 3,737.40 & 5 & 20.80 & 94 & 56.15 \\
            75T-L4                     & 848,962.94 & 196,870.68 & 177,932.19 & 474,160.07 & 3,655.80 & 5 & 21.60 & 81 & 52.89 \\ \midrule
            75L-L1                     & 647,231.00 & 111,886.08 & 54,340.37  & 481,004.55 & 2,767.80 & 3 & 7.00  & 96 & 65.90 \\
            75L-L2                     & 667,341.85 & 111,886.08 & 54,340.37  & 501,115.40 & 2,767.80 & 3 & 7.00  & 83 & 65.90 \\
            75L-L3                     & 659,756.14 & 111,886.08 & 65,173.62  & 482,696.43 & 2,769.20 & 3 & 10.80 & 96 & 80.13 \\
            75L-L4                     & 706,050.89 & 153,758.36 & 126,942.54 & 425,350.00 & 3,004.80 & 4 & 15.60 & 82 & 60.19 \\ \bottomrule
        \end{tabular}%
    }
\end{table}

\begin{table}
    \centering
    \caption{Solution details of the stochastic uncapacitated instances.\label{tab:val-comp-u-ss}}
    {\footnotesize \begin{tabular}{@{}lrrrrrrrrr@{}}
            \toprule
            \multicolumn{1}{c}{Inst}   &
            \multicolumn{1}{c}{TC}     &
            \multicolumn{1}{c}{LC}     &
            \multicolumn{1}{c}{HC}     &
            \multicolumn{1}{c}{DC}     &
            \multicolumn{1}{c}{HFlow}  &
            \multicolumn{1}{c}{\#Hubs} &
            \multicolumn{1}{c}{\#Veh1} &
            \multicolumn{1}{c}{\#Veh2} &
            \multicolumn{1}{c@{}}{\%VUtil}                                                                                    \\ \midrule
            25U-L1                     & 367,402.09 & 85,182.42  & 32,611.65 & 249,608.02 & 1,254.40 & 2 & 2.80  & 42 & 74.67 \\
            25U-L2                     & 323,541.65 & 40,332.49  & 0.00      & 283,209.16 & 0.00     & 1 & 0.00  & 32 &       \\
            25U-L3                     & 368,791.49 & 40,332.49  & 0.00      & 328,458.99 & 0.00     & 1 & 0.00  & 44 &       \\
            25U-L4                     & 323,541.65 & 40,332.49  & 0.00      & 283,209.16 & 0.00     & 1 & 0.00  & 32 &       \\ \midrule
            40U-L1                     & 500,827.50 & 63,320.02  & 23,373.65 & 414,133.84 & 1,786.60 & 2 & 3.80  & 67 & 78.36 \\
            40U-L2                     & 452,358.23 & 66,796.48  & 18,341.02 & 367,220.74 & 1,423.20 & 2 & 3.00  & 50 & 79.07 \\
            40U-L3                     & 508,996.02 & 66,796.48  & 26,492.58 & 415,706.96 & 1,489.80 & 2 & 5.20  & 67 & 89.53 \\
            40U-L4                     & 458,951.47 & 66,796.48  & 26,492.58 & 365,662.41 & 1,489.80 & 2 & 5.20  & 50 & 89.53 \\ \midrule
            50U-L1                     & 511,753.76 & 111,688.92 & 34,705.24 & 365,359.60 & 1,315.00 & 2 & 3.00  & 64 & 73.06 \\
            50U-L2                     & 523,764.08 & 111,688.92 & 37,018.93 & 375,056.23 & 1,433.60 & 2 & 3.20  & 55 & 74.67 \\
            50U-L3                     & 523,832.46 & 111,688.92 & 48,201.73 & 363,941.82 & 1,375.40 & 2 & 5.00  & 64 & 85.96 \\
            50U-L4                     & 533,190.15 & 97,154.21  & 29,264.83 & 406,771.11 & 1,368.00 & 2 & 5.00  & 55 & 85.50 \\ \midrule
            75U-L1                     & 652,115.96 & 124,218.94 & 63,955.10 & 463,941.91 & 2,498.80 & 3 & 7.00  & 96 & 59.50 \\
            75U-L2                     & 665,520.33 & 123,463.64 & 65,863.22 & 476,193.47 & 2,758.00 & 3 & 7.00  & 83 & 65.67 \\
            75U-L3                     & 664,294.83 & 124,218.94 & 72,641.08 & 467,434.81 & 2,603.80 & 3 & 10.20 & 96 & 79.77 \\
            75U-L4                     & 679,638.91 & 111,886.08 & 63,313.88 & 504,438.94 & 2,766.60 & 3 & 10.40 & 83 & 83.13 \\ \bottomrule
        \end{tabular}%
    }
\end{table}

We also compare the solution of the stochastic problem with that of the expected value problem (EVP). The EVP is solved as a deterministic problem, where the flow values $w_{ij}$ are replaced by the average of flows in five scenarios. 
Tables \ref{tab:vals-ssva-cap} and \ref{tab:vals-ssva-u} list the results for the capacitated and uncapacitated instances, respectively. We compare the solution in terms of the objective function value (TC) and selected hubs (Hubs). We also study a case where one makes the main location and allocation decisions based on the EVP before real demand uncertainty is observed. The `EVP feasible' column show that whether such a strategy would result in a feasible stochastic solution.

Tables \ref{tab:vals-ssva-cap} and \ref{tab:vals-ssva-u} indicate that the stochastic solution often suggests a hub topology to handle uncertainty in demand at the minimum possible cost.
For some instances (25L-L1, 25U-L1, and 40U-L1) the stochastic problem provided a better solution than the optimal EVP solution.
For other instances, the total cost of the stochastic solution is commonly equal or larger than that of the EVP.
This is due to the flow distribution in AP dataset. The amount of flow originating at each node in the AP instances is highly variable. All instances have a small set of nodes for which the outgoing flow is significantly large compared to the other nodes. As a result, these few nodes can have a large impact on the total transportation cost and, hence, bias the optimal location of the hubs.
Even though the location and allocation decisions made in the EVP solution yields a better total cost in most cases, they may result in an infeasible solution under demand variability. \autoref{tab:vals-ssva-cap} shows that when hubs are capacitated, one usually needs to solve the stochastic problem to make a reliable location/allocation decision.

\begin{table}
    \centering
    \caption{Comparison of the expected value problem and stochastic problem solutions for the capacitated instances.\label{tab:vals-ssva-cap}}
    {\footnotesize \begin{tabular}{@{}lrlrll@{}}
            \toprule
            \multicolumn{1}{c}{}     & \multicolumn{2}{c}{EVP solution} & \multicolumn{3}{c}{Stochastic solution}                                                                                           \\ \cmidrule(lr){2-3}\cmidrule(l){4-6}
            \multicolumn{1}{c}{Inst} & \multicolumn{1}{c}{TC}           & \multicolumn{1}{c}{Hubs}                & \multicolumn{1}{c}{TC} & \multicolumn{1}{c}{Hubs} & \multicolumn{1}{c@{}}{EVP feasible} \\ \midrule
            25T-L1                   & 392,649.19                       & 9 14                                    & 395,264.34             & 9 25                     & No                                  \\
            25T-L2                   & 353,656.86                       & 9 12                                    & 366,191.96             & 9 25                     & No                                  \\
            25T-L3                   & 396,235.86                       & 9 14                                    & 406,602.63             & 9 25                     & No                                  \\
            25T-L4                   & 358,374.68                       & 13 14                                   & 377,530.25             & 9 25                     & No                                  \\ \midrule
            25L-L1                   & 382,488.84                       & 9 13                                    & 381,937.67             & 9 23                     & No                                  \\
            25L-L2                   & 336,548.26                       & 11 14                                   & 339,249.83             & 11 14                    & Yes                                 \\
            25L-L3                   & 384,194.77                       & 14                                      & 393,147.18             & 9 14                     & No                                  \\
            25L-L4                   & 337,949.82                       & 14                                      & 348,719.00             & 11 14                    & No                                  \\ \midrule
            40T-L1                   & 549,604.29                       & 14 19 40                                & 580,503.96             & 14 25 40                 & No                                  \\
            40T-L2                   & 504,256.39                       & 14 19 40                                & 542,719.38             & 14 25 40                 & No                                  \\
            40T-L3                   & 573,355.17                       & 14 19 40                                & 611,117.13             & 10 14 19 38              & No                                  \\
            40T-L4                   & 528,325.86                       & 14 19 40                                & 568,371.24             & 10 14 19 38              & No                                  \\ \midrule
            40L-L1                   & 501,148.58                       & 14 19                                   & 501,505.40             & 14 19                    & No                                  \\
            40L-L2                   & 454,438.22                       & 19 22                                   & 455,346.14             & 14 19                    & No                                  \\
            40L-L3                   & 507,299.54                       & 14 19                                   & 509,126.60             & 14 19                    & No                                  \\
            40L-L4                   & 459,418.76                       & 19 22                                   & 463,004.29             & 14 19                    & No                                  \\ \midrule
            50T-L1                   & 574,216.50                       & 6 26 48                                 & 632,538.04             & 6 26 32                  & No                                  \\
            50T-L2                   & 581,310.90                       & 12 26 48                                & 638,805.43             & 6 26 32                  & No                                  \\
            50T-L3                   & 585,050.80                       & 6 26 48                                 & 633,528.79             & 6 26 32                  & No                                  \\
            50T-L4                   & 593,979.61                       & 6 26 48                                 & 645,546.49             & 12 26 28 48              & No                                  \\ \midrule
            50L-L1                   & 514,226.77                       & 3 24                                    & 519,399.98             & 3 24                     & No                                  \\
            50L-L2                   & 521,638.23                       & 3 26                                    & 527,339.56             & 3 24                     & No                                  \\
            50L-L3                   & 525,279.11                       & 3 24                                    & 532,303.70             & 3 24                     & No                                  \\
            50L-L4                   & 532,498.63                       & 3 24                                    & 538,954.06             & 15 24                    & No                                  \\ \midrule
            75T-L1                   & 738,010.36                       & 27 33 53 74                             & 800,332.99             & 27 33 38 53              & No                                  \\
            75T-L2                   & 754,232.54                       & 27 33 53 74                             & 863,739.63             & 17 27 33 53 74           & No                                  \\
            75T-L3                   & 742,446.44                       & 27 33 53 74                             & 824,989.34             & 27 33 38 53 74           & No                                  \\
            75T-L4                   & 761,060.60                       & 27 33 53 74                             & 848,962.94             & 17 27 33 53 74           & No                                  \\ \midrule
            75L-L1                   & 647,231.00                       & 33 40 53                                & 647,231.00             & 33 40 53                 & Yes                                 \\
            75L-L2                   & 667,341.85                       & 33 40 53                                & 667,341.85             & 33 40 53                 & Yes                                 \\
            75L-L3                   & 653,939.16                       & 33 40 53                                & 659,756.14             & 33 40 53                 & Yes                                 \\
            75L-L4                   & 674,085.26                       & 33 40 53                                & 706,050.89             & 3 33 40 53               & Yes                                 \\ \bottomrule
        \end{tabular}
    }
\end{table}

\begin{table}
    \centering
    \caption{Comparison of the expected value problem and stochastic problem solutions for the uncapacitated instances.\label{tab:vals-ssva-u}}
    {\footnotesize \begin{tabular}{@{}lrlrll@{}}
            \toprule
            \multicolumn{1}{c}{}     & \multicolumn{2}{c}{EVP solution} & \multicolumn{3}{c}{Stochastic solution}                                                                                           \\ \cmidrule(lr){2-3}\cmidrule(l){4-6}
            \multicolumn{1}{c}{Inst} & \multicolumn{1}{c}{TC}           & \multicolumn{1}{c}{Hubs}                & \multicolumn{1}{c}{TC} & \multicolumn{1}{c}{Hubs} & \multicolumn{1}{c@{}}{EVP feasible} \\ \midrule
            25U-L1                   & 368,791.49                       & 13                                      & 367,402.09             & 9 24                     & Yes                                 \\
            25U-L2                   & 323,541.65                       & 13                                      & 323,541.65             & 13                       & Yes                                 \\
            25U-L3                   & 368,791.49                       & 13                                      & 368,791.49             & 13                       & Yes                                 \\
            25U-L4                   & 323,541.65                       & 13                                      & 323,541.65             & 13                       & Yes                                 \\ \midrule
            40U-L1                   & 500,844.45                       & 19 22                                   & 500,827.50             & 14 19                    & Yes                                 \\
            40U-L2                   & 450,799.91                       & 19 22                                   & 452,358.23             & 19 22                    & Yes                                 \\
            40U-L3                   & 507,299.54                       & 14 19                                   & 508,996.02             & 19 22                    & Yes                                 \\
            40U-L4                   & 457,932.53                       & 19 22                                   & 458,951.47             & 19 22                    & Yes                                 \\ \midrule
            50U-L1                   & 509,695.46                       & 15 48                                   & 511,753.76             & 15 48                    & Yes                                 \\
            50U-L2                   & 521,450.39                       & 15 48                                   & 523,764.08             & 15 48                    & Yes                                 \\
            50U-L3                   & 523,191.94                       & 15 48                                   & 523,832.46             & 15 48                    & Yes                                 \\
            50U-L4                   & 532,498.63                       & 3 24                                    & 533,190.15             & 17 24                    & Yes                                 \\ \midrule
            75U-L1                   & 647,231.00                       & 33 40 53                                & 652,115.96             & 17 40 53                 & Yes                                 \\
            75U-L2                   & 665,520.33                       & 27 33 53                                & 665,520.33             & 27 33 53                 & Yes                                 \\
            75U-L3                   & 653,939.15                       & 33 40 53                                & 664,294.83             & 17 40 53                 & Yes                                 \\
            75U-L4                   & 674,085.29                       & 33 40 53                                & 679,638.91             & 33 40 53                 & Yes                                 \\ \bottomrule
        \end{tabular}
    }
\end{table}

\section{Instances with General Inter-hub Network Structure}
\label{sec:gen}

In cases where the inter-hub network does not follow a complete-network topology, the structure of the hub-level network becomes a part of the decision process.
In the incomplete networks, two hubs can be connected indirectly, which makes some OD paths to visit more than two hub nodes.
Under such circumstances, the quadratic constraint (4b) and inequality (16c) are no longer valid (see the original paper).
Therefore, to solve the problem using Benders decomposition, one needs to incorporate a flow conservation constraint for each hub node in the subproblem to ensure that the flows are correctly routed through the network.
We use the same flow variable $\bm{z}$ and, for any given solution $\paren{\overline{\bm{x}}, \overline{\bm{y}}}$ of the BMP (eq. 12 in the original paper), and define the Benders subproblem for the general network structures (G-BSP) as:

\begin{IEEEeqnarray}{llll}
    \IEEEyesnumber\label{eq:gbsp}
    \IEEEyessubnumber*
    \text{G-BSP}\paren{\overline{\bm{x}}, \overline{\bm{y}}}:\ \ & \min\ & \IEEEeqnarraymulticol{2}{l}{0} \label{eq:gpsp.0} \\ [1ex]
    & \text{s.t.:}\ & \sum_{i\in N} z_{ihk} \leq Q_{hk}\, y_{hk},\ h, k\in H, \\
    && \sum_{k\in H} z_{ihk} - \sum_{k\in H} z_{ikh} = O_i \overline{x}_{hi} - \sum_{j\in N} w_{ij} \overline{x}_{hj},\  & h\in H, i \in N. \IEEEeqnarraynumspace \label{eq:gpsp.00} \\
    && \bm{z} \in \mathbb{R}^{\card{N}\times\card{H}\times\card{H}}_{\geq 0}.
\end{IEEEeqnarray}

Constraint \eqref{eq:gpsp.00} is the flow conservation constraints, which can also be obtained by subtracting equation (9d) from equation (9c) in the BSP for complete inter-hub networks (see the original paper).

Let $\bm{\lambda}$ and $\bm{\eta}$ be the dual variables. Then, the dual of the G-BSP, namely the G-DSP, is defined by \eqref{eq:gdsp} given below.

\begin{IEEEeqnarray}{llll}
    \IEEEyesnumber\label{eq:gdsp}
    \IEEEyessubnumber*
    \text{G-DSP}\paren{\overline{\bm{x}}, \overline{\bm{y}}}:\ \ & \max\ & \IEEEeqnarraymulticol{2}{l}{\sum_{h \in H}\sum_{k \in H} Q_{hk}\overline{y}_{hk}\, \lambda_{hk} + \sum_{h \in H}\sum_{i \in N}\paren{O_i \overline{x}_{hi} - \sum_{j\in N} w_{ij} \overline{x}_{hj}} \eta_{hi}} \label{eq:gdsp.0} \IEEEeqnarraynumspace\\ [1ex]
    & \text{s.t.:}\ & \lambda_{hk} + \eta_{hi} - \eta_{ki} \leq 0,\  & h, k\in H, i \in N, \label{eq:gdsp.1} \\
    && \bm{\lambda} \in \mathbb{R}^{\card{H}\times\card{H}}_{\leq 0}, \bm{\eta} \in \mathbb{R}^{\card{H}\times\card{N}}. \label{eq:gdsp.00}
\end{IEEEeqnarray}

Please note that there exists at least one feasible solution to the G-DSP (e.g., the origin point). Let $\paren{\overline{\bm{\lambda}}, \overline{\bm{\eta}}}$ be the unbounded ray of the G-DSP when the G-BSP is infeasible. Then, one can apply the following feasibility cut to the BMP and continue to the next step of the Benders decomposition algorithm.

\begin{IEEEeqnarray}{l}
    \sum_{h \in H}\sum_{k \in H} Q_{hk}\overline{\lambda}_{hk}\,y_{hk} + \sum_{h \in H}\sum_{i \in N} \overline{\eta}_{hi}\paren{ O_i \, x_{hi}  - \sum_{j \in N} w_{ij}\, x_{hj}} \leq 0. \label{eq:gfcut}
\end{IEEEeqnarray}

\autoref{tab:perf-incomp} shows the performance measures of the proposed Benders decomposition approach (BD) in solving small and medium-size instances with general inter-hub network structure.
Our approach is able to solve capacitated and uncapacitated instances with up to 50 nodes in 12,000 seconds (3 hrs and 20 mins).
\autoref{tab:val-incomp} provides the solution details.

\begin{table}
    \centering
    \caption{BD performance in solving instances with general network structure.\label{tab:perf-incomp}}
    {\footnotesize \begin{tabular}{@{}lrrrrr@{}}
            \toprule
            \multicolumn{1}{c}{Inst}                   &
            \multicolumn{1}{c}{\#BNodes}               &
            \multicolumn{1}{c}{\#Cuts}                 &
            \multicolumn{1}{c}{CPU\textsuperscript{r}} &
            \multicolumn{1}{c}{CPU\textsuperscript{c}} &
            \multicolumn{1}{c@{}}{CPU (\%Gap)}                                                          \\ \midrule
            20T-L1                                     & 6.58e+3 & 86    & 3.78   & 4.10     & 8.90     \\
            20T-L2                                     & 4.40e+3 & 101   & 2.07   & 4.75     & 8.18     \\
            20T-L3                                     & 2.98e+4 & 64    & 1.61   & 2.30     & 5.91     \\
            20T-L4                                     & 3.28e+3 & 114   & 1.94   & 4.90     & 8.80     \\ \midrule
            20L-L1                                     & 7.48e+4 & 142   & 1.95   & 8.36     & 19.84    \\
            20L-L2                                     & 1.08e+4 & 53    & 1.74   & 1.82     & 3.60     \\
            20L-L3                                     & 1.38e+4 & 117   & 4.61   & 4.12     & 8.17     \\
            20L-L4                                     & 2.60e+3 & 86    & 3.89   & 2.22     & 5.33     \\ \midrule
            20U-L1                                     & 7.18e+4 & 134   & 2.80   & 6.96     & 12.31    \\
            20U-L2                                     & 6.67e+3 & 73    & 1.61   & 3.03     & 5.74     \\
            20U-L3                                     & 7.52e+3 & 101   & 3.95   & 5.04     & 10.11    \\
            20U-L4                                     & 3.26e+3 & 46    & 1.03   & 1.31     & 2.13     \\ \midrule
            25T-L1                                     & 1.03e+5 & 173   & 4.37   & 22.67    & 33.94    \\
            25T-L2                                     & 4.00e+4 & 183   & 8.12   & 26.59    & 38.82    \\
            25T-L3                                     & 7.85e+4 & 161   & 4.76   & 19.95    & 29.54    \\
            25T-L4                                     & 6.46e+4 & 200   & 4.64   & 29.34    & 38.03    \\ \midrule
            25L-L1                                     & 1.20e+4 & 120   & 3.48   & 11.85    & 14.09    \\
            25L-L2                                     & 9.37e+3 & 77    & 3.35   & 5.91     & 8.69     \\
            25L-L3                                     & 9.56e+3 & 76    & 3.92   & 5.41     & 8.15     \\
            25L-L4                                     & 3.92e+3 & 65    & 2.95   & 5.04     & 6.44     \\ \midrule
            25U-L1                                     & 1.52e+4 & 107   & 4.84   & 11.23    & 15.14    \\
            25U-L2                                     & 3.98e+3 & 65    & 3.22   & 3.60     & 4.64     \\
            25U-L3                                     & 1.21e+4 & 64    & 3.40   & 5.79     & 8.54     \\
            25U-L4                                     & 3.39e+3 & 58    & 2.91   & 3.59     & 5.12     \\ \midrule
            40T-L1                                     & 8.50e+6 & 281   & 43.67  & 253.58   & 1,844.98 \\
            40T-L2                                     & 7.49e+7 & 349   & 31.68  & 319.70   & (1.30)   \\
            40T-L3                                     & 9.47e+6 & 110   & 26.58  & 75.17    & 1,309.78 \\
            40T-L4                                     & 5.86e+7 & 324   & 113.64 & 318.75   & (3.55)   \\ \midrule
            40L-L1                                     & 3.08e+4 & 177   & 8.67   & 110.04   & 127.10   \\
            40L-L2                                     & 4.86e+2 & 64    & 16.23  & 16.86    & 21.03    \\
            40L-L3                                     & 2.35e+5 & 98    & 21.22  & 59.43    & 96.40    \\
            40L-L4                                     & 1.11e+3 & 63    & 20.02  & 27.15    & 34.77    \\ \midrule
            40U-L1                                     & 4.58e+4 & 174   & 56.06  & 104.77   & 124.73   \\
            40U-L2                                     & 2.66e+2 & 56    & 17.84  & 16.49    & 20.40    \\
            40U-L3                                     & 5.64e+4 & 109   & 5.91   & 63.94    & 79.48    \\
            40U-L4                                     & 6.21e+4 & 55    & 11.52  & 14.18    & 24.66    \\ \midrule
            50T-L1                                     & 2.19e+7 & 628   & 581.83 & 2,547.89 & (5.05)   \\
            50T-L2                                     & 3.00e+7 & 495   & 626.04 & 1,584.53 & (3.66)   \\
            50T-L3                                     & 9.27e+6 & 1,114 & 441.57 & 5,232.04 & (8.85)   \\
            50T-L4                                     & 1.09e+7 & 990   & 417.81 & 4,267.35 & (7.56)   \\ \midrule
            50L-L1                                     & 2.14e+7 & 56    & 24.00  & 82.43    & 2,344.54 \\
            50L-L2                                     & 5.73e+5 & 64    & 54.92  & 104.51   & 191.52   \\
            50L-L3                                     & 6.36e+7 & 197   & 166.53 & 392.47   & (0.84)   \\
            50L-L4                                     & 6.20e+7 & 199   & 120.96 & 411.69   & (2.88)   \\ \midrule
            50U-L1                                     & 1.56e+5 & 111   & 140.76 & 213.65   & 253.24   \\
            50U-L2                                     & 1.41e+6 & 290   & 253.03 & 944.41   & 1,299.17 \\
            50U-L3                                     & 5.36e+7 & 245   & 598.74 & 586.54   & (4.59)   \\
            50U-L4                                     & 3.62e+7 & 292   & 32.81  & 979.80   & (5.12)   \\ \bottomrule
        \end{tabular}%
    }
\end{table}


\begin{table}
    \centering
    \caption{Solution details of the instances with general network structure.\label{tab:val-incomp}}
    {\footnotesize \begin{tabular}{@{}lrrrrrrrr@{}}
            \toprule
            \multicolumn{1}{c}{Inst}   &
            \multicolumn{1}{c}{TC}     &
            \multicolumn{1}{c}{LC}     &
            \multicolumn{1}{c}{HC}     &
            \multicolumn{1}{c}{DC}     &
            \multicolumn{1}{c}{\#Hubs} &
            \multicolumn{1}{c}{\#Veh1} &
            \multicolumn{1}{c}{\#Veh2} &
            \multicolumn{1}{c@{}}{\%VUtil}                                                                      \\ \midrule
            20T-L1                     & 393,933.67 & 75,763.01  & 28,123.55 & 290,047.10 & 2 & 4  & 48 & 73.36 \\
            20T-L2                     & 358,634.95 & 75,763.01  & 28,123.55 & 254,748.38 & 2 & 4  & 36 & 73.36 \\
            20T-L3                     & 400,964.56 & 75,763.01  & 35,154.44 & 290,047.10 & 2 & 6  & 48 & 91.70 \\
            20T-L4                     & 365,665.84 & 75,763.01  & 35,154.44 & 254,748.38 & 2 & 6  & 36 & 91.70 \\ \midrule
            20L-L1                     & 370,221.96 & 77,781.69  & 33,217.39 & 259,222.88 & 2 & 3  & 48 & 79.09 \\
            20L-L2                     & 334,821.03 & 77,781.69  & 33,217.39 & 223,821.95 & 2 & 3  & 36 & 79.09 \\
            20L-L3                     & 375,980.02 & 28,436.87  & 0.00      & 347,543.15 & 1 & 0  & 50 &       \\
            20L-L4                     & 336,616.68 & 28,436.87  & 0.00      & 308,179.81 & 1 & 0  & 38 &       \\ \midrule
            20U-L1                     & 370,221.96 & 77,781.69  & 33,217.39 & 259,222.88 & 2 & 3  & 48 & 79.09 \\
            20U-L2                     & 334,821.03 & 77,781.69  & 33,217.39 & 223,821.95 & 2 & 3  & 36 & 79.09 \\
            20U-L3                     & 375,980.02 & 28,436.87  & 0.00      & 347,543.15 & 1 & 0  & 50 &       \\
            20U-L4                     & 336,616.68 & 28,436.87  & 0.00      & 308,179.81 & 1 & 0  & 38 &       \\ \midrule
            25T-L1                     & 462,716.46 & 135,626.35 & 34,539.64 & 292,550.48 & 3 & 5  & 50 & 72.05 \\
            25T-L2                     & 417,631.78 & 121,863.97 & 35,851.40 & 259,916.42 & 3 & 5  & 37 & 71.85 \\
            25T-L3                     & 472,623.79 & 135,626.35 & 45,053.60 & 291,943.84 & 3 & 8  & 50 & 84.85 \\
            25T-L4                     & 427,906.07 & 135,626.35 & 45,053.60 & 247,226.13 & 3 & 8  & 37 & 84.85 \\ \midrule
            25L-L1                     & 434,368.05 & 128,271.77 & 32,736.50 & 273,359.79 & 2 & 3  & 51 & 81.39 \\
            25L-L2                     & 387,027.79 & 114,600.55 & 33,305.03 & 239,122.22 & 2 & 3  & 37 & 77.81 \\
            25L-L3                     & 441,156.96 & 114,600.55 & 46,256.98 & 280,299.42 & 2 & 5  & 52 & 98.27 \\
            25L-L4                     & 394,911.22 & 114,600.55 & 46,256.98 & 234,053.69 & 2 & 5  & 37 & 98.27 \\ \midrule
            25U-L1                     & 423,528.15 & 98,853.64  & 35,333.69 & 289,340.83 & 2 & 3  & 51 & 81.39 \\
            25U-L2                     & 369,691.77 & 40,332.49  & 0.00      & 329,359.28 & 1 & 0  & 39 &       \\
            25U-L3                     & 425,922.44 & 40,332.49  & 0.00      & 385,589.95 & 1 & 0  & 53 &       \\
            25U-L4                     & 369,691.77 & 40,332.49  & 0.00      & 329,359.28 & 1 & 0  & 39 &       \\ \midrule
            40T-L1                     & 511,901.10 & 131,166.19 & 45,923.73 & 334,811.19 & 3 & 5  & 58 & 73.90 \\
            40T-L2                     & 484,728.44 & 108,415.85 & 50,290.22 & 326,022.37 & 3 & 5  & 45 & 76.64 \\
            40T-L3                     & 520,285.33 & 131,166.19 & 59,967.56 & 329,151.59 & 3 & 8  & 58 & 87.18 \\
            40T-L4                     & 496,223.51 & 108,415.85 & 65,425.68 & 322,382.03 & 3 & 8  & 45 & 91.15 \\ \midrule
            40L-L1                     & 459,954.90 & 63,320.02  & 18,452.88 & 378,182.00 & 2 & 3  & 59 & 80.09 \\
            40L-L2                     & 424,777.22 & 63,320.02  & 18,452.88 & 343,004.33 & 2 & 3  & 46 & 80.77 \\
            40L-L3                     & 465,508.65 & 63,320.02  & 25,629.00 & 376,559.63 & 2 & 5  & 59 & 94.27 \\
            40L-L4                     & 430,719.72 & 63,320.02  & 25,629.00 & 341,770.71 & 2 & 5  & 46 & 94.27 \\ \midrule
            40U-L1                     & 459,954.90 & 63,320.02  & 18,452.88 & 378,182.00 & 2 & 3  & 59 & 80.09 \\
            40U-L2                     & 424,777.22 & 63,320.02  & 18,452.88 & 343,004.33 & 2 & 3  & 46 & 80.77 \\
            40U-L3                     & 465,508.65 & 63,320.02  & 25,629.00 & 376,559.63 & 2 & 5  & 59 & 94.27 \\
            40U-L4                     & 430,719.72 & 63,320.02  & 25,629.00 & 341,770.71 & 2 & 5  & 46 & 94.27 \\ \midrule
            50T-L1                     & 566,348.12 & 216,636.58 & 59,496.23 & 290,215.31 & 4 & 7  & 61 & 59.97 \\
            50T-L2                     & 585,187.18 & 216,636.58 & 59,496.23 & 309,054.42 & 4 & 7  & 55 & 60.60 \\
            50T-L3                     & 589,489.09 & 216,636.58 & 85,826.22 & 287,026.29 & 4 & 11 & 61 & 70.65 \\
            50T-L4                     & 607,270.69 & 216,636.58 & 89,532.31 & 301,101.81 & 4 & 10 & 55 & 79.24 \\ \midrule
            50L-L1                     & 493,572.33 & 144,265.22 & 52,197.68 & 297,109.44 & 3 & 5  & 62 & 64.67 \\
            50L-L2                     & 519,870.45 & 144,265.22 & 52,197.68 & 323,407.56 & 3 & 5  & 56 & 64.67 \\
            50L-L3                     & 516,022.37 & 93,859.16  & 36,802.80 & 385,360.42 & 2 & 4  & 63 & 84.94 \\
            50L-L4                     & 539,289.51 & 93,859.16  & 36,802.80 & 408,627.56 & 2 & 4  & 57 & 84.94 \\ \midrule
            50U-L1                     & 481,048.54 & 134,880.74 & 49,529.37 & 296,638.43 & 3 & 5  & 62 & 66.39 \\
            50U-L2                     & 511,300.54 & 134,880.74 & 49,529.37 & 326,890.43 & 3 & 5  & 56 & 66.39 \\
            50U-L3                     & 495,933.46 & 96,315.70  & 48,886.64 & 350,731.13 & 2 & 5  & 63 & 93.46 \\
            50U-L4                     & 531,845.53 & 134,880.74 & 66,417.71 & 330,547.12 & 3 & 7  & 56 & 83.76 \\ \bottomrule
        \end{tabular}%
    }
    {}
\end{table}